\documentclass[journal ]{new-aiaa}
\usepackage[utf8]{inputenc}
\usepackage{textcomp}
\usepackage{xcolor}
\usepackage{graphicx}
\usepackage{amsmath}
\usepackage{cancel}
\usepackage[version=4]{mhchem}
\usepackage{siunitx}
\usepackage{longtable,tabularx}
\usepackage{subcaption}
\usepackage{threeparttable}
\usepackage[normalem]{ulem}
\title{Numerical Dissipation Based Error Estimators and Grid Adaptation for Large Eddy Simulation}

\author{Yao Jiang \footnote{Graduate Student, Student Member AIAA.} and Siva Nadarajah \footnote{Associate Professor, AIAA Associate Fellow.}}
\affil{McGill University, Montreal, QC, H3A 0G4}

\begin{document}

\maketitle

\begin{abstract}
Grid adaptation for implicit Large Eddy Simulation (LES) is a non-trivial challenge due to the inherent coupling of the modeling and numerical errors. An attempt to address the challenge first requires a comprehensive assessment and then the development of error estimators to highlight regions that require refinement. Following the work of Schranner et al.~\cite{schranner2015assessing}, a novel approach to estimate the numerical dissipation of the turbulent kinetic energy (TKE) equations is proposed. The presented approach allows the computation of the local numerical dissipation for arbitrary curvilinear grids through a post-processing procedure. This method, as well as empirical and kinetic-energy-based~\cite{castiglioni2014numerical,sun2018implicit} approaches, are employed to estimate the inherent numerical TKE. We incorporate the numerical TKE to evaluate an effective eddy viscosity, an effective Kolmogorov length scale, and an effective TKE to build a family of Index Quality (IQ) based error estimators. The proposed IQ based estimators are then assessed and utilized to show their effectiveness through an application of grid adaptation for the periodic hill test case and transitional flow over the SD 7003 airfoil. Numerical results are validated through a comparison against reference LES and experimental data. Flow over the adapted grids appear better abled to capture pertinent flow features and integrated functions, such as the lift and drag coefficients.
\end{abstract}

\section*{Nomenclature}

\noindent

{\renewcommand\arraystretch{1.0}
\noindent\begin{longtable*}{@{}l @{\quad=\quad} l@{}}
$\epsilon$  & dissipation \\
$\nu$ &    viscosity \\
$\eta$& Kolmogorov length scale \\
$k$ & turbulent kinetic energy \\
$\epsilon$ & dissipation rate \\
$s_{ij}$ & string rate tensor\\
$\Delta$ & filter size \\
$h$   & grid cell size \\
$\Delta t$ & time step \\
$\overline{\bullet}$ & time-averaged value \\
$f, g$   & generic functions\\
$\alpha, m, n$   & generic constants
\end{longtable*}}

\section{Introduction}

There are multitude of approaches in CFD to model and capture turbulence. Direct Numerical Simulation (DNS), for example, resolves all scales down to the grid scale. However, since the range increases rapidly with Reynolds number, the computational cost of DNS limits its current application to low-Reynolds number flows over flat-plate type geometries, precluding industrially relevant problems. In industry, the Reynolds-averaged Navier-Stokes (RANS) method is the default approach to modeling turbulent flow. By using the Reynolds decomposition, the RANS method provides an approximation to the time-averaged solutions and models the influence of all scales of turbulence motion. Despite its widespread use in industry, RANS often fails to capture unsteady flow physics accurately. This is due in large part to the lack of a clear separation of scales between the time scales associated with the mean flow variation and the turbulent fluctuations.

Large Eddy Simulation (LES) has thus emerged as an alternative approach to provide a compromise between computational expense and the extent to which scales are resolved. In LES, a filter is employed for the state variables in space in order to separate the large and small scales. The unsteady large scales are captured while the unresolved small scales are modeled. Two filtering processes, namely implicit and explicit filtering, have been developed by the CFD community in the past two decades. When using an integral based finite-volume discretization, if a built-in filter (commonly of the grid cell size~\cite{denaro2011does}) is implicitly applied~\cite{denaro2011does,lund2003} with no explicit implementation of a filtering process, then the approach is deemed implicit filtering. Unlike implicit filtering, explicit filtering techniques are applied directly to the Navier-Stokes (NS) equation or only on the non-linear terms as an independent filter, whose width is fixed manually throughout the grid and decoupled with the grid cell size. As a consequence, a grid independent LES is achieved. The advantages and disadvantages of both techniques have been exclusively studied~\cite{vasilyev1997general,lund2003,brandt2006,deconinck2008,denaro2011does,bull2016explicit} and is beyond the scope of this article. To date implicit filtering has generally been adopted by commercial CFD softwares for industrial problems thanks to its flexibility and simplicity. 

{\color{black} Although LES is less computationally expensive compared to DNS, as the Reynolds number of the flow increases, the range of scales to be resolved also increase, which is especially true in the near-wall region where the length scales are much smaller than the boundary layer thickness ($\eta \ll \delta$, typically $\Delta x\Delta y\Delta z \sim 2500\delta^3$ for wall-resolved LES~\cite{chapman1979}). It is crucial to wisely refine the computational grid in critical regions in order to minimize error as well as to keep the computational cost affordable. However, it is not obvious how to evaluate the total error in LES, as both grid resolution and model accuracy play an equally important role nor are the two sources of error mutually exclusive. A commonly accepted concept~\cite{vermeire2016implicit,GHOSAL1996,celik2009assessment} to define the total error is to split the error in LES into a numerical error, $\epsilon_{\text{num}}$ defined as the difference between a fine and coarse grid LES, and a modeling error, $\epsilon_{\text{mod}}$ which represents the departure of the fine grid LES from a filtered~\cite{celik2009assessment} or unfiltered~\cite{meyers2011error} DNS result. Assuming that the errors arising from the filtering and the numerical discretization are unrelated, we can represent the total error for an arbitrary variable $\mathbf{u}$ as,
\begin{equation}
    \begin{aligned}
    &\epsilon_{\text{tot}}=\epsilon_{\text{num}}+\epsilon_{\text{mod}}=\mathbf{u}_{\text{DNS}}-\mathbf{u}_{\text{LES}},
    \end{aligned}
\end{equation}
where $\epsilon_{\text{num}}=\mathbf{u}_{\text{fine LES}}-\mathbf{u}_{\text{LES}}$ and $\epsilon_{\text{mod}}=\mathbf{u}_{\text{DNS}}-\mathbf{u}_{\text{fine LES}}$. In the context of explicit filtering, through the definitions of the filtering width and the grid size, the two error contributions can be dissociated to achieve grid independent solutions. However, in implicit filtering, the common practice of using the grid size within the SGS model strongly couples the errors and renders grid independence intractable. In the latter approach grid refinement not only leads to a reduction of the truncation error in the numerical discretization and the capture of smaller scales but modifies the inherent dissipation of the SGS model. If we are to further organize the sub-contributors to the numerical error into dissipation, dispersion and aliasing; then the primary contributor is dissipation, as noted by Celik et al.~\cite{celik2009assessment}, as it alters the apparent Reynolds number of the flow. 

The challenge of separating numerical from modeling errors further confounds the problem in creating a general framework to develop error estimators.} Existing error estimators and grid adaptation procedures for RANS solvers could serve as references~\cite{hauser2015adaptive} but cannot be directly applied to LES due to the nature of the turbulence model. Despite the number of publications devoted to the assessment of LES quality~{\cite{john2003analysis,klein2005attempt,freitag2006improved,klein2008assessment,xing2015general,ries2018evaluating}}, few of them could be directly applied to grid adaptation for practical LES applications, due to the following limitations:

\begin{enumerate}
    \item Different from RANS, where the grid size mainly affects the discretization error; in LES, both modeling and numerical errors are implicitly dependent on each other and are non-trivial to estimate;
    \item LES is intrinsically unsteady and chaotic, such that the error estimators should include some averaging process in space and/or time;
    \item Different from detailed LES quality assessment, the error estimation for LES grid adaptation only allows for single-grid estimators for practical reasons.
\end{enumerate}

Geurts and Frohlich proposed the activity parameter~\cite{geurts2002framework} for LES quality assessment based on the estimation of turbulent and viscous dissipation. Celik et al.~\cite{celik2005index,celik2009assessment} introduced a family of Index Qualities for LES error estimation, based on the effective Kolmogorov scale $\eta_{\text{eff}}$ or the eddy viscosity $\nu_{\text{eff}}$ which incorporates the contribution from the SGS model and the numerical dissipation. The estimator was extended using the proportion of resolved to total TKE based on Pope's suggestion~\cite{pope2001turbulent} that a ``good'' LES approach should resolve at least 80\% of the total TKE. Instead of using the proportion of resolved TKE, Antepara et al.~\cite{antepara2015parallel} used the residual velocity magnitude (without scaling) as an error indicator. {\color{black} In the approaches above, there is at least one parameter~\cite{xing2015general} which requires an estimation based on empirical equations. 
  For instance, the Index Quality of Celik et al.~\cite{celik2005index,celik2009assessment} relies on an empirical estimation of the numerical TKE for single grid estimators. In addition, the Index Quality based on $\eta_{\text{eff}}$ and $\nu_{\text{eff}}$ call for proper tuning of coefficients to shape the Index Quality curve such that a value of 0.8 corresponds to a sufficient LES grid. Lastly, in~\cite{antepara2015parallel}, grid refinement is initiated when the local residual velocity magnitude reaches a threshold that depends on a user-defined averaging function and case-dependent constants.}

To avoid the use of empirical equations to estimate the amount of numerical dissipation, multiple grid estimates were introduced, where Celik and Karatekin~\cite{celik1997numerical} proposed a Richardson extrapolation based approach. Klein~\cite{klein2005attempt} extended the Richardson extrapolation method of~\cite{celik1997numerical} using multiple grids to evaluate the Index Quality, where he assumed that both modeling and numerical errors for implicit filtering can be combined and represented as a unique function of grid size. The multiple grid Index Quality approach was further extended in~\cite{freitag2006improved,xing2015general} to separate the modeling and numerical errors, in addition to the coupling error through systematic independent grid refinements and model variations. {\color{black} 
However neither single nor multiple grid estimators provide for a directional approach to adapt the filter or the grid size. Toosi and Larsson~\cite{toosi2017anisotropic,toosi2018grid,toosi2020towards} proposed a small energy density indicator depending on the residual velocity magnitude which included directional information. A second anisotropic test filter~\cite{toosi2017anisotropic} is then applied in the post-processing stage to estimate the directional error to study the solution sensitivity to a change in filter size. However, the approach is based on the assumption that the application of a larger test filter in post-processing is similar to the application of a larger filter during LES.} 


The quality of the error estimator hinges on an appropriate estimation of the inherent numerical dissipation, regardless of the LES approach. At present, multiple grid estimators do not offer a cost effective approach to provide an accurate estimate of both the modeling and numerical errors. The development of an efficient single-grid estimator relies on the study of crucial flow features to quantify the impact of numerical dissipation. {\color{black} Domaradzki et al.~\cite{domaradzki2003effective} proposed an approach to evaluate the numerical dissipation of kinetic energy in implicit LES in spectral space. The evolution of kinetic energy is crucial for accurate simulations of flow~\cite{jameson2008formulation}, where there is an energy cascade between the different eddy scales. Currently, kinetic energy preserving schemes~\cite{jameson2008formulation,schranner2015assessing} have contributed towards this objective; however, the rationale proposed by Domaradzki et al.~\cite{domaradzki2003effective} still allows for a quantification and perhaps further refines the estimate of numerical dissipation present in the flow. The approach was extended to the physical space by Schranner et al.~\cite{schranner2015assessing}, where the authors employed the numerical dissipation of the kinetic energy equation as an LES quality assessment criterion and presented results for the Taylor-Green vortex. Castiglioni and Domaradzki~\cite{castiglioni2015numerical} confirmed the findings of~\cite{schranner2015assessing} for a realistic flow over a NACA0012 airfoil at a Mach number of 0.4 and Reynolds number 50,000. However, the authors~\cite{castiglioni2015numerical} highlighted a limitation of the approach, where negative numerical dissipation was observed in regions of laminar flow with weak velocity gradients. Nevertheless, Cadieux et al.~\cite{cadieux2017effects} further promoted the approach for explicit LES {\color{black} and Sun and Domaradzki~\cite{sun2018implicit}} for adaptive implicit filtering; where the authors performed adaptive grid refinement using a block-based ratio of the numerical and physical dissipation rates as a criterion.

The present work extends~\cite{schranner2015assessing,cadieux2017effects,sun2018implicit} to evaluate the numerical dissipation of the turbulent kinetic energy (TKE) equation and employ it within {\it a posteriori} single-grid error estimators of Celik et al.~\cite{celik2009assessment}. Our extension differs from the original approach in four aspects. First, instead of using the numerical dissipation as the error estimator, we further develop the method to fit within the framework of Index Quality through the derivation of the numerical TKE. Second, unlike the block-based evaluation approach in~\cite{sun2018implicit}, we employ cell-based error estimation and provide direct indication of the local truncation error. Third, instead of the kinetic energy equation, we utilize the turbulent kinetic energy equation to estimate the numerical TKE. Fourth, we engage an automatic grid adaptation approach. The proposed error estimators are easily implemented within a post-processing procedure without modifying the underlying flow solver. Our aim is to focus on developing cost effective single grid estimators, provide a comprehensive assessment of them, and demonstrate the effectiveness through a single grid adaptation cycle. A conventional mesh adaptation approach that repeatedly adapts the grid is currently infeasible for LES, and hence our objective is to present an approach where an acceptable but not sufficient grid is improved in a single cycle. Thus our contributions are: First, a novel error estimator that appropriately evaluates the total error, which includes contributions from both numerical and modeling errors in LES for a given grid, especially focusing on the impact of numerical dissipation. Second, a comparison of several error estimators on the performance of estimating the grid quality for LES. Third, a proposal for a complete pipeline for LES grid adaptation and demonstration of the effectiveness of the LES error estimator in providing significant improvement in engineering interest.} The paper is organized as follows: Section~\ref{sec:err_estim} details the error estimators for LES based on Index of Quality and the classical approach to estimate the numerical TKE. Section~\ref{sec:num_tke} presents novel forms for evaluating numerical TKE. In Section~\ref{sec:num_res}, the numerical validation of error estimators are performed through a complete cycle of grid adaptation via the periodic hill and SD7003 airfoil test cases.

\section{Error Estimation}
\label{sec:err_estim}

{\color{black} Index of Quality is a family of error estimators for LES which allows for single-grid error estimation~\cite{celik2005index} and this concept is adopted in our grid adaptation process. In this section the formulas and the implementation of three classical Index of Quality error estimators are presented.}
    
\subsection{Index Quality $IQ_{\nu}$}
The Index Quality based on effective eddy viscosity takes the form of
    \begin{equation}
    IQ_{\nu}=\frac{1}{1+\alpha_{\nu}(\frac{\nu_{\text{eff}}}{\nu})^n},
    \label{iq_nu}
    \end{equation}
where the eddy viscosity $\nu_{\text{eff}}$ incorporates the contribution from the SGS model and the numerical dissipation. The implementation of $ IQ_{\nu} $ can be directly fulfilled by using
    \begin{equation}
    \nu_{\text{eff}}=\nu_{\text{num}}+\nu_{\text{sgs}},
    \label{eq:nu_eff}
    \end{equation} 
where $\nu_{\text{num}}$ is estimated from
    \begin{equation}
    \nu_{\text{num}}=\text{sgn}(k_{\text{num}})C_{\nu}\Delta\sqrt{\text{abs}(k_{\text{num}})},
    \label{nu_num}
    \end{equation}
by analogy to the SGS model formula~\cite{yoshizawa1982statistically}, {\color{black}while $k_{\text{num}}$ is an estimate of the numerical TKE due to numerical dissipation inherent in the discretization of the governing equations. The inclusion of the numerical TKE was proposed by Celik et al.~\cite{celik2009assessment} with an empirical formula which will be discussed in Section~\ref{sec:k_num_emp}. In Section~\ref{sec:num_tke} we extend the estimator to two novel versions depending on the approach to evaluate $k_{\text{num}}$.} 
The $IQ_{\nu}$ estimator will target regions where either the SGS eddy viscosity or the numerical eddy viscosity is large compared to the kinematic viscosity $\nu$.
\subsection{Index Quality $IQ_{\eta}$}
The Index Quality based on the effective Kolmogorov scale is defined as,
    \begin{equation}
    IQ_{\eta}=\frac{1}{1+\alpha_{\eta}(\frac{h}{\eta_{\text{eff}}})^m},
    \label{iq_eta}
    \end{equation}
where the effective Kolmogorov length scale $\eta_{\text{eff}}$ is defined by
    \begin{equation}
    \eta_{\text{eff}}=(\frac{\nu^3}{\epsilon})^{1/4}.
    \end{equation}
The dissipation $\epsilon$ is estimated from
    \begin{equation}
    s=(\frac{\epsilon}{\nu_{\text{eff}}})^{1/2},
    \label{epsilon}
    \end{equation}
where $ s=\sqrt{2s_{ij}s_{ij}} $ is the square root of the double inner product of the mean strain rate tensor $ s_{ij} $ and the effective eddy viscosity, $\nu_{\text{eff}}$ is evaluated from Eq.~\ref{eq:nu_eff}. Lastly, the grid cell length scale $h$ is given by
    $h=V^{1/3}$,
where $V$ is the volume of the cell and $\alpha_{\eta}$ and $m$ are empirical constants. The estimator targets regions where the grid size $h$ is large compared to the effective Kolmogorov length scale. The estimator is related to $IQ_{\nu}$ by involving the same effective eddy viscosity $\nu_{\text{eff}}$, while the $IQ_{\eta}$ estimator mainly targets regions where the grid cell size is not sufficiently small to resolve flow structures on the length scale of the effective Kolmogorov.

\subsection{Index Quality $IQ_k$}

The Index Quality was extended using the proportion of resolved to total TKE,
    \begin{equation}
    IQ_{k}=\frac{k_{\text{res}}}{k_{\text{res}}+k_{\text{eff}}},
    \label{iq_k}
    \end{equation}
where the resolved TKE, $k_{\text{res}}$, is obtained by adding the contribution from the diagonal terms in the Reynolds stress tensor, obtained by subtracting the instantaneous flow field from the time-averaged solution,
    \begin{equation}
    k_{\text{res}}=\frac{1}{2}\sum_{i=1}^{3} u_i'u_i',
    \end{equation}
    with
    $u_i'=u_i-\bar{u}_i$,
where $u_i$ and $\bar{u}_i$ are the instantaneous and time-averaged velocities. For a single-grid approach, $ k_{\text{eff}} $ is the sum of the modeled TKE, $ k_{\text{sgs}} $, and numerical TKE, $ k_{\text{num}} $. The modeled TKE can be derived from the SGS eddy viscosity using the formula
\begin{equation}
    k_{\text{sgs}}=(\frac{\nu_{\text{sgs}}}{C_{\nu}\Delta})^2.
\end{equation}
The evaluation of $ k_{\text{num}} $ is nontrivial and the classical approach relies on an empirical formula~\cite{celik2009assessment}. We proposed different formulas for the derivation of $ k_{\text{num}} $ which will be presented in Section~\ref{sec:num_tke}. The estimator targets regions where insufficient percentage of TKE is resolved by the flow solver.

{\color{black} The above mentioned Index Quality error estimators follow the same conception and the constants in each formula are tuned~\cite{celik2009assessment} such that the value of all $IQ$-based error estimators are in a comparable range. However, the estimators focus on different scales within the flow field. $IQ_{\nu}$ focuses on the effective eddy viscosity, thus the product of a velocity and length scale in the dissipative range. {\color{black} $IQ_{\eta}$ focuses on the Kolmogorov length scale and the indicator is related to $IQ_{\nu}$ with the inclusion of the local strain rate, $s$, as an additional factor, such that the $IQ_{\eta}$ estimator also targets regions with high strain rate (e.g. the turbulent mixing layer and boundary layer).} $IQ_k$ focuses on the resolution of TKE. An empirical laminar flow correction factor proposed by Celik et al.~\cite{celik2009assessment} was applied for all three $IQ$ estimators to avoid an incorrect estimation in near wall regions, where the molecular viscosity is much higher than the eddy viscosity.}

\subsection{Empirical Evaluation of the Numerical TKE}
\label{sec:k_num_emp}
{\color{black} All of the above formulas depend on an estimation of numerical TKE, $k_{\text{num}}$; however, its assessment is nontrivial for single-grid simulations. Considering that $k_{\text{sgs}}$ and $k_{\text{num}}$ are related to the length scale $\Delta$ and $h$ respectively, Celik et al.~\cite{celik2009assessment} proposed a direct way of evaluating $k_{\text{num}}$ on a single-grid by assuming a linear relationship between $k_{\text{num}}$ and $k_{\text{sgs}}$, and the coefficient reflects the ratio of $h/\Delta$,}
    \begin{equation}
    k_{\text{num}}=C_n(\frac{h}{\Delta})^2 k_{\text{sgs}}.
    \label{k_num}
    \end{equation}
    
In the case of implicit filtering, the filter size is directly related to the grid size, such that the only parameter to tune is $C_n$, {\color{black} which is properly tuned following~\cite{jiang2020assessment}.}
Although the uniform turbulence may show a linear relation between $k_{\text{sgs}}$ and $k_{\text{num}}$, the linear assumption has two primary drawbacks because of the nature of the error even if the coefficient $C_n$ is tuned:
    \begin{itemize}
        \item The SGS model primarily contributes to the modeled TKE in the high wavenumber region while the numerical dissipation contributes to the numerical TKE in a much wider wavenumber range;
        \item When approaching the wall, the modeled TKE tends to zero, while clearly $k_{\text{num}}$ is expected to be non-zero. Therefore $k_{\text{sgs}}$ should not serve as an indicator to estimate $k_{\text{num}}$ in the near-wall region.
    \end{itemize}


{\color{black} An alternative way to study the impact of numerical dissipation is through the framework proposed by Schranner et al.~\cite{schranner2015assessing,cadieux2017effects}. The method estimates the numerical dissipation in the kinetic energy equation and further derives the numerical eddy viscosity, $\nu_{\text{num}}$. The application to LES was performed without~\cite{schranner2015assessing} and with~\cite{cadieux2017effects} an explicit SGS model. The extension of this approach to estimate the numerical TKE forms the fundamental contribution of this article and is presented in the following Section}

\section{Evaluation of Numerical TKE}
\label{sec:num_tke}
{\color{black} In this section, we focus on two novel approaches to estimate the numerical TKE. Starting from the LES governing equation, we detail the procedure of~\cite{cadieux2017effects} to estimate the numerical eddy viscosity, $\nu_{\text{num}}$, and extend the approach to evaluate the numerical TKE, $k_{\text{num}}$, in Section~\ref{sec:ke_diss}. In Section~\ref{sec:tke_diss} we present an alternative novel approach to estimate the numerical TKE. Two evaluation methods are implemented in the Index Quality, leading to two {novel error estimator families} for LES.}

\subsection{Governing Equation}
For LES to correctly resolve the flow, sub-filtered terms in the Navier-Stokes equation need to be modeled by a subgrid-scale (SGS) model. In this work, we solve the Navier-Stokes equations using implicit filtering and a wall-adapting local eddy-viscosity turbulence model~\cite{nicoud1999subgrid}. Using Einstein notation, the governing momentum and energy equations can be obtained,
\begin{equation}
\begin{aligned}
    \frac{\partial \bar{u}_{i}}{\partial t}+\frac{\partial \bar{u}_i \bar{u}_{i}}{\partial x_i}&=-\frac{1}{\rho}\frac{\partial \bar{p}}{\partial x_i}+\nu \frac{\partial \bar{\tau}_{ij}}{\partial x_j}-\frac{\partial \tau_{ij}^{SGS}}{\partial x_j},\\
    \frac{\partial \bar{e}}{\partial t}+\frac{\partial \bar{u}_{i} \bar{e}}{\partial x_i}&=-\frac{1}{\rho}\frac{\partial {u}_{i}\bar{p}}{\partial x_i}+\nu \frac{\partial {u}_{j}\bar{\tau}_{ij}}{\partial x_j}-\frac{\partial {u}_{j}\tau_{ij}^{SGS}}{\partial x_j}-\frac{\partial}{\partial x_i}\left( k\frac{\partial \bar{T}}{\partial x_i} \right),
\end{aligned}
\label{eq:governing_eq}
\end{equation}
with
\begin{equation}
    \bar{\tau}_{ij}=\frac{\partial \bar{u}_{i}}{\partial x_j}+\frac{\partial \bar{u}_{j}}{\partial x_i},
\end{equation}
where $\bar{u}_{i}$, $\bar{p}$, $\bar{e}$, $\bar{T}$ and $\bar{\tau}_{ij}$ are the filtered velocity, pressure, energy, temperature and stress tensor, and the filtering operation is defined as a convolution integral with the kernel function $G$,
\begin{equation}
    \overline{\bar{u}_{i}(\mathbf{x})}=\iiint G(\mathbf{x}-\mathbf{r},\Delta)u_{i}(\mathbf{r})d\mathbf{r}.
\end{equation}
For any finite-volume discretization employed, an implicit filter with a top-hat kernel function are implicitly applied to the flow variables~\cite{denaro2011does}, and the filter size is assumed to be equal to or proportional to the grid space $h$. The subgrid scale stress tensor is employed to model the unclosed sub-filtered part
\begin{equation}
    \tau_{ij}^{SGS}=\overline{u_iu_j}-\bar{u}_i\bar{u}_j.
\end{equation}
When an eddy viscosity model is employed, the Boussinesq hypothesis is applied to the subgrid scale stress tensor with the introduction of an eddy viscosity, $\nu_t$,
\begin{equation}
    \tau_{ij}^{SGS}=\nu_t(\frac{\partial \bar{u}_{i}}{x_j}+\frac{\partial \bar{u}_{j}}{x_i}).
\end{equation}
A Wall-Adapting Local Eddy-viscosity (WALE) model~\cite{nicoud1999subgrid} was adopted in the current study with the expression of eddy viscosity
\begin{equation}
    \nu_{t} = \Delta _s^2 \frac{(S_{ij}^{d} S_{ij}^{d})^{3/2}}{(\overline{S}_{ij} \overline{S}_{ij})^{5/2} + (S_{ij}^{d} S_{ij}^{d})^{5/4}},
\end{equation}
where
\begin{equation}
    \begin{aligned}
    \Delta _s &= C_w V^{1/3}, \\
    S_{ij}^{d} &= \frac{1}{2} \left( \overline{g}_{ij}^{2} + \overline{g}_{ji}^{2}  \right) - \frac{1}{3} \delta_{ij} \overline{g}_{kk}^{2}, \\ 
    \overline{g}_{ij} &= \frac{\partial \overline{u_i}}{\partial x_{j}}, \\
    \overline{g}_{ij}^{2} &= \overline{g}_{ik} \overline{g}_{kj},\\
    \bar S_{ij}  &= \frac{1}{2}\left( {\frac{{\partial \bar u_i }}{{\partial x_j }} + \frac{{\partial \bar u_j }}{{\partial x_i }}} \right),
    \end{aligned}
\end{equation}
and the constant $ C_w = 0.325 $. {\color{black} The model is based on the commonly held believe that the turbulent eddy viscosity is assumed to be proportional to the resolved velocity gradients and hence the impact of the unresolved scales are to dissipate the larger resolved scales.}

\subsection{Evaluation from Kinetic Energy Numerical Dissipation {\color{black} (KE-based approach)}}
\label{sec:ke_diss}
    
{\color{black} Following the procedure described by Schranner et al.~\cite{schranner2015assessing}, we present briefly the derivation of the numerical dissipation of kinetic energy, {\color{black} or KE numerical dissipation for short}. The transport equation for the kinetic energy $e_{\text{kin}}$ is formed by separating the transport equation for internal energy from the conservation of total energy. To simplify our notation, we choose to omit the filter symbol~($\;\bar{ }\;$) and assume all variables obtained in the following sections are LES resolved values. Assuming incompressible flow condition, the final form can be defined as},
    \begin{equation}
    \frac{\partial e_{\text{kin}}}{\partial t}+\frac{\partial (u_i e_{\text{kin}})}{\partial x_i}=-\frac{u_i}{\rho}\frac{\partial p}{\partial x_i}+(\nu+\nu_t) u_j\frac{\partial \tau_{ij}}{\partial x_i},
    \label{eq:trans_k}
    \end{equation}
    where $\nu_t$ is the subgrid‐scale eddy viscosity, and
    \begin{equation}
    \tau_{ij}=\frac{\partial u_i}{\partial x_j}+\frac{\partial u_j}{\partial x_i}.
    \end{equation}
    We then employ the same numerical scheme that was used to discretize the Navier-Stokes solver to discretize the transport of kinetic energy to estimate $e_{\text{kin}}$. In the following equation, the left-hand-side represents the discretization of Eq.~\ref{eq:trans_k}; while the non-zero right-hand-side provides the estimate of the numerical dissipation of kinetic energy $\epsilon_n$,
    \begin{equation}
    (\frac{\partial e_{\text{kin}}}{\partial t})_d+(NS)_d e_{\text{kin}}=\underbrace{(\frac{\partial e_{\text{kin}}}{\partial t})+(NS) e_{\text{kin}}}_{=0}+\epsilon_n,
    \label{eq:discretization_trans_ke}
    \end{equation}
    {\color{black} where $(NS)_d$ is the discretized Navier–Stokes operator. Schranner et al.~\cite{schranner2015assessing} employed both high-order and second-order finite volume schemes and demonstrated that although a high-order scheme is preferred to minimize the numerical errors, the method could be applied to second-order schemes and the difference was minimal. Castiglioni et al.~\cite{castiglioni2014numerical,castiglioni2015numerical} and Cadieux et al.~\cite{cadieux2017effects} confirmed that using the same discretization order as in the flow solver provided for sufficient accuracy for various test cases.} {\color{black} The residual $\epsilon_n$ serves as an estimate of the numerical dissipation of kinetic energy. The residual is largely a function of the truncation error of the numerical scheme and includes both dissipative and dispersive errors. For the purposes of this work where $\epsilon_n$ will be employed towards a grid adaptation indicator, an accurate approximate is not as critical as a consistent trend in the distribution of the numerical dissipation.} Next, we rewrite Eq.~\ref{eq:discretization_trans_ke} for the finite volume framework employed in this work as, 
    \begin{equation}
    -\epsilon_n=e_{\text{kin},t}+F_{e_{\text{kin}}}+F_{ac}+F_{\nu}+\epsilon_{\nu},
    \label{eq:diss_fv}
    \end{equation}
    where the right hand side includes the total time-rate of change of kinetic energy:
    \begin{itemize}
        \item $e_{\text{kin},t}=\frac{\partial e_{\text{kin}}}{\partial t}$,
        \item $F_{e_{\text{kin}}}=\frac{\partial (u_j e_{\text{kin}})}{\partial x_j}$,
        \item $F_{ac}=\frac{1}{\rho}\frac{\partial (p u_i)}{\partial x_i}$,
        \item $F_{\nu}=-(\nu+\nu_t) \frac{\partial (\tau_{ij}u_i) }{\partial x_j}$,
        \item $\epsilon_{\nu}=(\nu+\nu_t)\tau_{ij}\frac{\partial u_i}{\partial x_j}$.
    \end{itemize}

\noindent {\color{black} The first term on the right-hand-side, $e_{\text{kin},t}$ denotes the time rate of change of kinetic energy within the control volume, which we calculate using a second-order central difference using  $e_{\text{kin}}$ at the $t_{n-1}$ and $t_{n+1}$ time steps. The term averages to be zero for a statistically steady regime. The following term on the right-hand-side, $F_{e_{\text{kin}}}$ represents the convective flux term and tracks the net rate of change of kinetic energy through the control volume. The term $F_{ac}$ represents the acoustic flux driven by pressure work. The viscous flux, $F_{\nu}$, represents the transport of kinetic energy by Reynolds stresses; while, $\epsilon_{\nu}$ is the kinetic energy dissipation term. Each cell is treated as a control volume and flux terms are evaluated on each cell surface.}

In order to achieve a local cell-based error estimation, Eq.~\ref{eq:diss_fv} is evaluated locally in each computational cell, providing local $\epsilon_n$ values. {\color{black} Eq.~\ref{eq:diss_fv} is then evaluated at each time interval $\Delta t$ (in the order of $0.05h_{\text{ref}}/v_{\text{ref}}$ where $h_{\text{ref}}$ and $v_{\text{ref}}$ are the reference length and velocity scales), which lead to time and spanwise averaged values of $\bar{\epsilon}_{n}$. Following the concept proposed in~\cite{schneider2010}, a temporal and spatial along the $z$ (spanwise) averaged numerical eddy viscosity $\nu_{\text{num}}$ is estimated by analogy to the formula for the viscous dissipation term,}
    \begin{equation}
    \bar{\epsilon}_{n}=\nu_{\text{num}}\overline{\tau_{ij}\frac{\partial u_i}{\partial x_j}},
    \label{num-vis}
    \end{equation}
    such that
    \begin{equation}
    \nu_{\text{num}}=\bar{\epsilon}_{n}/\overline{\tau_{ij}\frac{\partial u_i}{\partial x_j}}.
    \label{eq:ratio}
    \end{equation}
    {\color{black}A temporal and spanwise averaged numerical TKE, $k_{\text{num}}$, is estimated from $\nu_{\text{num}}$ through Eq.~\ref{nu_num}.} The process of $\epsilon_n$ evaluation is similar to the residual evaluation procedure in the flow solver. The implementation of the approach could easily be carried out for an arbitrary flow solver. {However, the approach shows two drawbacks when applied for grid adaptation for LES:}
    
    \begin{itemize}
        \item Due to the nature of the term, the dissipation $\epsilon_{\nu}$ is always positive; the unsteady term $e_{\text{kin},t}$ tends to zero after time-averaging for a statistically steady case, while the flux terms, $F_{e_{\text{kin}}}$, $F_{ac}$ and $F_{\nu}$, have the possibility of showing negative values, such that the total numerical dissipation could result in negative values in certain regions where the dissipative effect is negligible compared to the convective effect. For example, $F_{e_{\text{kin}}}$ and $F_{ac}$ could become the leading terms in \color{black} both steady or laminar regions with high mean flow velocities and low turbulent intensities as well as low-velocity gradients; where it leads to negative numerical dissipation ~\cite{castiglioni2015numerical}. Therefore in practice, we only consider positive values of the {\color{black} numerical dissipation,
        \begin{equation}
            \bar{\epsilon}_{n,pos}=\max(\bar{\epsilon},0),
            \label{eq:pos_diss}
        \end{equation}
        to avoid steady or laminar regions with negative values from being targeted for grid refinement. Since the motivation is to approximate the truncation error, an alternate procedure would be to take the absolute value of the numerical dissipation. However, in practice, we observed that laminar regions with large convective acceleration yielded larger $\epsilon_n$ values. Since our objective is to perform grid adaptation in turbulent regions of the flow, Eq.~\ref{eq:pos_diss} proved to be the more appropriate choice.}
        \item {\color{black} The current approach is focused on the kinetic energy numerical dissipation which incorporates both mean flow kinetic energy and turbulent kinetic energy numerical dissipations, where both parts contribute to the evaluation of $\nu_{\text{num}}$. However, the concept of the $IQ_k$ error estimator is based on the percentage of captured turbulent kinetic energy instead of the kinetic energy, thus ideally $IQ_k$ should only be formulated based on the numerical dissipation of the turbulent kinetic energy.}
    \end{itemize}

    \begin{figure}[hbt!]
        \begin{center}
            \includegraphics[clip=true, trim= 1.0cm 3.0cm 1.0cm 0.5cm,width=.6\columnwidth]{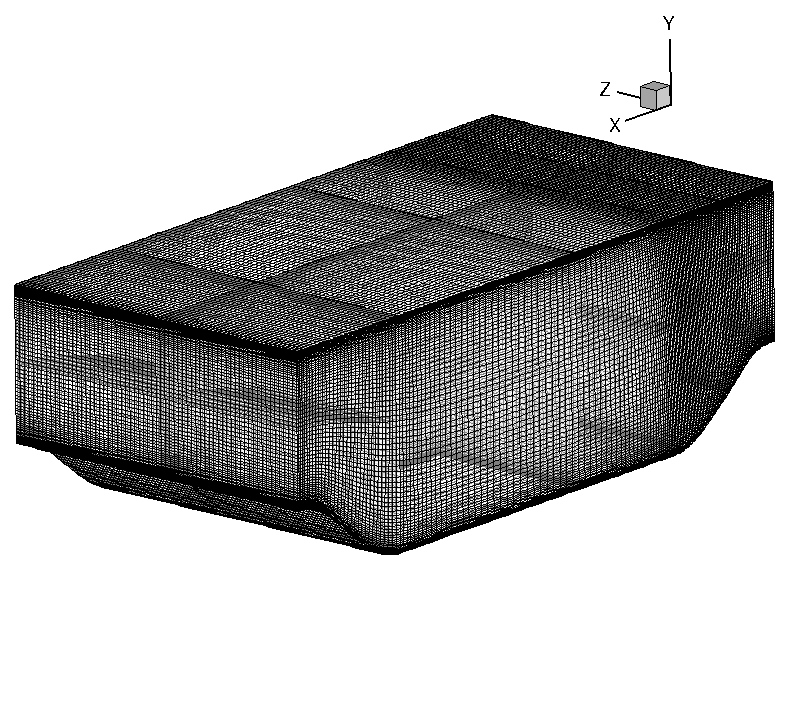}
        \end{center}
        \caption{Periodic hill grid.}
        \label{figure:perhill_geometry} 
    \end{figure}

\subsection{Evaluation from the Turbulent Kinetic Energy Numerical Dissipation {\color{black} (TKE-based approach)}}
\label{sec:tke_diss}
{\color{black} In this subsection, we present the novel aspect of this work. We believe that the most appropriate way to estimate numerical TKE, or the amount of turbulent kinetic energy dissipated by the discretization error, is to extend the approach of~\cite{schranner2015assessing} to the transport equation for the TKE. A generalization of Eq.~\ref{eq:discretization_trans_ke} for the turbulent kinetic energy equation can be expressed as,
    \begin{equation}
    (\frac{\partial k}{\partial t})_d+(NS)_d k=\underbrace{(\frac{\partial k}{\partial t})+(NS) k}_{=0}+\epsilon_n.
    \label{eq:discretization_trans_tke}
    \end{equation}
In LES, Eq.~\ref{eq:discretization_trans_tke} refers to the resolved TKE evolution and $\epsilon_n$ estimates the influence of numerical dissipation on the resolved TKE, or TKE numerical dissipation. Since most of the SGS model is based on eddy viscosity, we focus on the derivation of the discretized resolved turbulent kinetic energy equation for LES with the inclusion of $\nu_t$. Since we aim to achieve a \textit{posteriori} error estimators through a post-processing approach, as employed for the numerical dissipation of kinetic energy we engage the same numerical discretization as the flow solver to assess the TKE numerical dissipation. We begin by expressing the instantaneous momentum equation for LES,}
\begin{equation}
\frac{\partial u_i}{\partial t}+u_k\frac{\partial u_i}{\partial x_k}=-\frac{1}{\rho}\frac{\partial p}{\partial x_i}+(\nu+\nu_t) \frac{\partial^2 u_i}{\partial x_k \partial x_k},
\label{eq:inst}
\end{equation}
where the instantaneous resolved velocity $u_i$ could be decomposed into a mean and fluctuation part,
\begin{equation}
u_i=\bar{u_i}+u'_i,
\end{equation}
where the mean part is the time-averaging over the entire simulation period,
\begin{equation}
\bar{u_i}=\frac{1}{T}\int_T u_i dt
\end{equation}
and the time-averaged resolved momentum equation for incompressible LES leads to
\begin{equation}
\frac{\partial \bar{u_i}}{\partial t}+\bar{u_k}\frac{\partial \bar{u_i}}{\partial x_k}+\frac{\partial}{\partial x_k}\overline{u'_iu'_k}-\cancelto{0}{\overline{u'_i\frac{\partial u'_k}{\partial x_k}}}=-\frac{1}{\rho}\frac{\partial \bar{p}}{\partial x_i}+\nu \frac{\partial^2 \bar{u_i}}{\partial x_k \partial x_k}+\bar{\nu_t} \frac{\partial^2 \bar{u_i}}{\partial x_k \partial x_k}+\overline{\nu'_t \frac{\partial^2 u'_i}{\partial x_k \partial x_k}}.
\label{eq:avg}
\end{equation}
We can subtract Eq.~\ref{eq:inst} by Eq.~\ref{eq:avg} to acquire the instantaneous fluctuation equation
\begin{equation}
\begin{aligned}
&\frac{\partial u'_i}{\partial t}+\bar{u_k}\frac{\partial u'_i}{\partial x_k}+u'_k\frac{\partial \bar{u_i}}{\partial x_k}+u_k'\frac{\partial u'_i}{\partial x_k}-\frac{\partial}{\partial x_k}\overline{u'_iu'_k}=\\
&-\frac{1}{\rho}\frac{\partial p'}{\partial x_i}+\nu \frac{\partial^2 u'_i}{\partial x_k \partial x_k}+\bar{\nu_t} \frac{\partial^2 u'_i}{\partial x_k \partial x_k}+\nu_t' \frac{\partial^2 \bar{u_i}}{\partial x_k \partial x_k}+\nu_t' \frac{\partial^2 u'_i}{\partial x_k \partial x_k}-\overline{\nu'_t \frac{\partial^2 u'_i}{\partial x_k \partial x_k}}.
\label{eq:fluc}
\end{aligned}
\end{equation}
Next we premultiply Eq.~\ref{eq:fluc} by $u'_j$,
\begin{equation}
\begin{aligned}
&u'_j\frac{\partial u'_i}{\partial t}+u'_j\bar{u_k}\frac{\partial u'_i}{\partial x_k}+u'_ju'_k\frac{\partial \bar{u_i}}{\partial x_k}+u'_ju'_k\frac{\partial u'_i}{\partial x_k}-u'_j\frac{\partial}{\partial x_k}\overline{u'_iu'_k}=\\
&-\frac{u'_j}{\rho}\frac{\partial p'}{\partial x_i}+\nu u'_j\frac{\partial^2 u'_i}{\partial x_k\partial x_k}+\bar{\nu_t} u'_j\frac{\partial^2 u'_i}{\partial x_k\partial x_k}+\nu'_t u'_j\frac{\partial^2 \bar{u_i}}{\partial x_k\partial x_k}+\nu'_t u'_j\frac{\partial^2 u'_i}{\partial x_k\partial x_k}-u'_j\overline{\nu'_t \frac{\partial^2 u'_i}{\partial x_k. \partial x_k}}.
\end{aligned}
\label{eq:reynolds1}
\end{equation}
Switching $i$ and $j$ in Eq.~\ref{eq:reynolds1} and adding itself to Eq.~\ref{eq:reynolds1} yields,
\begin{equation}
\begin{aligned}
&\frac{\partial u'_iu'_j}{\partial t}+\bar{u_k}\frac{\partial u'_iu'_j}{\partial x_k}+u'_iu'_k\frac{\partial \bar{u_j}}{\partial x_k}+u'_ju'_k\frac{\partial \bar{u_i}}{\partial x_k}+\frac{\partial}{\partial x_k}u'_iu'_ju'_k-\cancelto{0}{u'_iu'_j\frac{\partial u'_k}{\partial x_k}}-\overbrace{u'_j\frac{\partial}{\partial x_k}\overline{u'_iu'_k}-u'_i\frac{\partial}{\partial x_k}\overline{u'_ju'_k}}^{\text{disappear when averaging in time}}=\\
&-\frac{u'_i}{\rho}\frac{\partial p'}{\partial x_j}-\frac{u'_j}{\rho}\frac{\partial p'}{\partial x_i}+\nu u'_i\frac{\partial^2 u'_j}{\partial x_k\partial x_k}+\nu u'_j\frac{\partial^2 u'_i}{\partial x_k\partial x_k}\\
&+\bar{\nu_t} u'_j\frac{\partial^2 u'_i}{\partial x_k\partial x_k}+\nu'_t u'_j\frac{\partial^2 \bar{u_i}}{\partial x_k\partial x_k}+\nu'_t u'_j\frac{\partial^2 u'_i}{\partial x_k\partial x_k}-u'_j\overline{\nu'_t \frac{\partial^2 u'_i}{\partial x_k \partial x_k}}\\
&+\bar{\nu_t} u'_i\frac{\partial^2 u'_j}{\partial x_k\partial x_k}+\nu'_t u'_i\frac{\partial^2 \bar{u_j}}{\partial x_k\partial x_k}+\nu'_t u'_i\frac{\partial^2 u'_j}{\partial x_k\partial x_k}-u'_i\overline{\nu'_t \frac{\partial^2 u'_j}{\partial x_k \partial x_k}}=\\
&-\frac{1}{\rho}(u'_i\frac{\partial p'}{\partial x_j}+u'_j\frac{\partial p'}{\partial x_i})
+\nu \frac{\partial^2 u'_iu'_j}{\partial x_k\partial x_k}-2\nu \frac{\partial u'_i}{\partial x_k}\frac{\partial u'_j}{\partial x_k}
+\bar{\nu_t} \frac{\partial^2 u'_iu'_j}{\partial x_k\partial x_k}-2\bar{\nu_t} \frac{\partial u'_i}{\partial x_k}\frac{\partial u'_j}{\partial x_k}\\
&+\nu'_t \frac{\partial^2 u'_iu'_j}{\partial x_k\partial x_k}-2\nu'_t \frac{\partial u'_i}{\partial x_k}\frac{\partial u'_j}{\partial x_k}+\nu'_t u'_i\frac{\partial^2 \bar{u_j}}{\partial x_k\partial x_k}+\nu'_t u'_j\frac{\partial^2 \bar{u_i}}{\partial x_k\partial x_k}\underbrace{-u'_i\overline{\nu'_t \frac{\partial^2 u'_j}{\partial x_k \partial x_k}}-u'_j\overline{\nu'_t \frac{\partial^2 u'_i}{\partial x_k \partial x_k}}}_{\text{disappear when averaging in time}},\\
\end{aligned}
\label{eq:reynolds}
\end{equation}
{\color{black}which is the transport equation for the Reynolds stress tensor $\tau'_{ij}=u'_iu'_j$.} We now average Eq.~\ref{eq:reynolds}
\begin{equation}
\begin{aligned}
&\frac{\partial \overline{u'_iu'_j}}{\partial t}+\bar{u_k}\frac{\partial \overline{u'_iu'_j}}{\partial x_k}+\overline{u'_iu'_k}\frac{\partial \bar{u_j}}{\partial x_k}+\overline{u'_ju'_k}\frac{\partial \bar{u_i}}{\partial x_k}+\frac{\partial}{\partial x_k}\overline{u'_iu'_ju'_k}=\\
&-\frac{1}{\rho}(\overline{u'_i\frac{\partial p'}{\partial x_j}}+\overline{u'_j\frac{\partial p'}{\partial x_i}})+\overline{(\nu+\nu_t) \frac{\partial^2 u'_iu'_j}{\partial x_k\partial x_k}}-2\overline{(\nu+\nu_t) \frac{\partial u'_i}{\partial x_k}\frac{\partial u'_j}{\partial x_k}}+\overline{\nu'_t u'_i}\frac{\partial^2 \bar{u_j}}{\partial x_k\partial x_k}+\overline{\nu'_t u'_j}\frac{\partial^2 \bar{u_i}}{\partial x_k\partial x_k}.
\end{aligned}
\label{eq:reynolds_avg}
\end{equation}
The contraction operation to the Reynolds stress tensor leads to the definition of TKE,
\begin{equation}
contr(\tau'_{ij})=\delta_{ij}\tau'_{ij}=u'_iu'_i,
\end{equation}
using Einstein summation convention. We are able to apply the contraction operation to Eq.~\ref{eq:reynolds_avg} to obtain the equation for the TKE budget,
\begin{equation}
\begin{aligned}
&\frac{\partial \overline{\frac{1}{2}u'_iu'_i}}{\partial t}+\bar{u_k}\frac{\partial \overline{\frac{1}{2}u'_iu'_i}}{\partial x_k}+\overline{u'_iu'_k}\frac{\partial \bar{u_i}}{\partial x_k}+\frac{\partial}{\partial x_k}\frac{1}{2}\overline{u'_iu'_iu'_k}=\\
&-\frac{1}{\rho}(\overline{u'_i\frac{\partial p'}{\partial x_i}})+\overline{(\nu+\nu_t) \frac{\partial^2 \frac{1}{2}u'_iu'_i}{\partial x_k\partial x_k}}-\overline{(\nu+\nu_t) \frac{\partial u'_i}{\partial x_k}\frac{\partial u'_i}{\partial x_k}}+\overline{\nu'_t u'_i}\frac{\partial^2 \bar{u_i}}{\partial x_k\partial x_k},
\end{aligned}
\end{equation}
which can be re-written as the time averaged resolved TKE equation:
\begin{equation}
\begin{aligned}
&\frac{\partial \overline{k}}{\partial t}+\bar{u_k}\frac{\partial \overline{k}}{\partial x_k}+\overline{u'_iu'_k}\frac{\partial \bar{u_i}}{\partial x_k}+\frac{\partial}{\partial x_k}\frac{1}{2}\overline{u'_iu'_iu'_k}=\\
&-\frac{1}{\rho}(\overline{u'_i\frac{\partial p'}{\partial x_i}})+\overline{(\nu+\nu_t) \frac{\partial^2 \frac{1}{2}u'_iu'_i}{\partial x_k\partial x_k}}-\overline{(\nu+\nu_t) \frac{\partial u'_i}{\partial x_k}\frac{\partial u'_i}{\partial x_k}}+\overline{\nu'_t u'_i}\frac{\partial^2 \bar{u_i}}{\partial x_k\partial x_k},
\end{aligned}
\end{equation}
such that at each time step,
\begin{equation}
\begin{aligned}
&\frac{\partial k}{\partial t}+\bar{u_j}\frac{\partial k}{\partial x_j}+u'_iu'_j\frac{\partial \bar{u_i}}{\partial x_j}+\frac{\partial}{\partial x_j}ku'_j=\\
&-\frac{u'_i}{\rho}(\frac{\partial p'}{\partial x_i})+(\nu+\nu_t) \frac{\partial^2 k}{\partial x_j\partial x_j}-(\nu+\nu_t) \frac{\partial u'_i}{\partial x_j}\frac{\partial u'_i}{\partial x_j}+\nu'_t u'_i\frac{\partial^2 \bar{u_i}}{\partial x_j\partial x_j}\\
&+\text{first-order fluctuation terms},
\end{aligned}
\end{equation}
which can be re-written as
\begin{equation}
\begin{aligned}
&\frac{\partial k}{\partial t}+\frac{\partial}{\partial x_j}(ku_j+\frac{u'_jp'}{\rho})=\cancelto{0}{k\frac{\partial \bar{u_j}}{\partial x_j}}+\cancelto{0}{\frac{p'}{\rho}\frac{\partial u'_j}{\partial x_j}}\\
&-u'_iu'_j\frac{\partial \bar{u_i}}{\partial x_j}+(\nu+\nu_t) \frac{\partial^2 k}{\partial x_j\partial x_j}-(\nu+\nu_t) \frac{\partial u'_i}{\partial x_j}\frac{\partial u'_i}{\partial x_j}+\nu'_t u'_i\frac{\partial^2 \bar{u_i}}{\partial x_j\partial x_j}\\
&+\text{first-order fluctuation terms},
\end{aligned}
\label{eq:tke_inst}
\end{equation}
{\color{black}where all first-order fluctuation terms vanish once time averaging is applied.} Using the same technique as the previous subsection, we express the numerical dissipation of TKE as the discretization of all the terms of Eq.~\ref{eq:tke_inst} in the flow solver,
        \begin{equation}
        \begin{aligned}
        -\epsilon_n=&(\frac{\partial k}{\partial t})_d+(NS)_d k\\=&[\frac{\partial k}{\partial t}+\frac{\partial}{\partial x_j}(ku_j+\frac{u'_jp}{\rho})+u'_iu'_j\frac{\partial \bar{u_i}}{\partial x_j}-(\nu+\bar{\nu_t}) \frac{\partial^2 k}{\partial x_j\partial x_j}\\
        +&(\nu+\bar{\nu_t}) \frac{\partial u'_i}{\partial x_j}\frac{\partial u'_i}{\partial x_j}-\nu'_t u'_i\frac{\partial^2 \bar{u_i}}{\partial x_j\partial x_j}]_d,
        \label{eq:tke_num}
        \end{aligned}
        \end{equation}
        {\color{black} where $(NS)_d$ is the discretized Navier–Stokes operator. It is sufficient to use a numerical discretization scheme consistent with the flow solver to estimate the numerical dissipation of TKE in an arbitrary CFD code. Eq.~\ref{eq:tke_num} is re-written for implementation purpose as},
        \begin{equation}
        -\epsilon_n=k_{t}+F_{k}+F_{ac}+F_{\nu}+P+\epsilon_{\nu}+\epsilon_{inter},
        \label{eq:diss_tke_fv}
        \end{equation}
        where the are defined as:
        \begin{itemize}
            \item $k_t=\frac{\partial k}{\partial t}$,
            \item $F_{k}=\frac{\partial}{\partial x_j}(ku_j)$,
            \item $F_{ac}=\frac{\partial}{\partial x_j}(\frac{u'_jp'}{\rho})$,
            \item $F_{\nu}=-(\nu+\nu_t) \frac{\partial^2 k}{\partial x_j\partial x_j}$,
            \item $P=u'_iu'_j\frac{\partial \bar{u_i}}{\partial x_j}$,
            \item $\epsilon_{\nu}=(\nu+\nu_t) \frac{\partial u'_i}{\partial x_j}\frac{\partial u'_i}{\partial x_j}$,
            \item $\epsilon_{interact}=-\nu'_t u'_i\frac{\partial^2 \bar{u_i}}{\partial x_j\partial x_j}$.
        \end{itemize}
        {\color{black} The first term on the right-hand-side $k_t$ is the time rate of change of TKE within the control volume, which averages to zero for a statistically steady flow. The temporal derivative is approximated using the
        turbulent kinetic energy data from the simulations at time steps $t_{n-1}$ and $t_{n+1}$; $F_{k}$ is the convection of TKE; $F_{ac}$ is the transport of TKE by pressure fluctuation; and $F_{\nu}$ represents the turbulent and viscous transport of TKE, which exhibits a diffusive property. $P$ is the turbulent production of TKE from the interaction between the mean flow and the fluctuation; while,  $\epsilon_{\nu}$ is the dissipation term and $\epsilon_{interact}$ is the interactive term which primarily displays a dissipative behaviour. We adopted a cell-based approach where flux terms are evaluated on all cell surfaces. The numerical dissipation, $\epsilon_n$ is then evaluated at each time interval $\Delta t$ (in the order of $0.05h_{\text{ref}}/v_{\text{ref}}$ where $h_{\text{ref}}$ and $v_{\text{ref}}$ are the reference length and velocity scales), and averaged in the temporal and spanwise directions, resulting in $\bar{\epsilon}_{n}$. The dissipated TKE due to the existence of numerical error, or numerical TKE $k_{\text{num}}$, is then estimated using a local length scale $l_{\text{scale}}$ and velocity scale $u_{\text{scale}}$,
        \begin{equation}
        k_{\text{num}}=\bar{\epsilon}_n\frac{l_{\text{scale}}}{u_{\text{scale}}}.
        \end{equation}
        Considering that the numerical dissipation occurs at the scale of the cell size $h$ and the velocity fluctuation scale, we employ
        \begin{equation}
        l_{\text{scale}}=V^{\frac{1}{3}},\quad \mbox{and}
        \end{equation}
        \begin{equation}
        u_{\text{scale}}=|\overline{u'}|=\sqrt{\frac{2}{3}k_{\text{res}}},
        \end{equation}
        where $V$ is the local cell volume and $k_{\text{res}}$ is the local resolved TKE.}

{\color{black}We can then finally evaluate the Index Qualities. This is accomplished through three possible approaches to evaluate $k_{\text{num}}$ and subsequently the Index Quality error estimators to derive $IQ_{\nu,\eta,k}$; an empirical approach which leads to $IQ_{\nu,\eta,k-emp}$, {\color{black} the KE-based approach} on the evaluation of KE numerical dissipation deemed $IQ_{\nu,\eta,k-ke}$ and the evaluation of TKE numerical dissipation {\color{black} using the TKE-based approach} designated as $IQ_{\nu,\eta,k-tke}$.}

\section{Numerical Results}
\label{sec:num_res}
\subsection{Numerical Setup}
\label{sec:solver}
SYN3D is a finite-volume based multi-block structured flow solver developed at the Computational Aerodynamics Group at McGill University. The solver is designed for three-dimensional Navier-Stokes equations, coupled with the energy equation with the inclusion of various turbulence models including RANS, LES and hybrid RANS/LES models. Parallelization is enabled through a Message Passing Interface (MPI) standard. The governing equations can be rewritten over the computational domain in semi-discrete form as,
\begin{equation}
    V\frac{\partial \mathbf{u}}{\partial t}+\mathbf{R}(\mathbf{u})=0,
\end{equation}
where $V$ is the cell volume and $\mathbf{R}(\mathbf{u})$ is the residual, which is a function of the state variable $\mathbf{u}$. For LES, the solved state variable is the filtered velocity and energy variables,
\begin{equation}
    \mathbf{u}=[\bar{u}_1,\bar{u}_2,\bar{u}_3,\bar{e}].
\end{equation}
The residual term $\mathbf{R}(\mathbf{u})$ contains the convective and dissipative as well as artificial dissipation flux gradients. We employ a cell-centered discretization with an artificial dissipation scheme~\cite{jameson1981numerical} based on a blend of second-and fourth-order differences which result from the sum of additional first-and third-order fluxes. The parameters in the scheme are tuned based on the Isotropic Decaying Turbulence test case~\cite{jiang2020assessment}.

\subsection{Octree data structure}
In order to achieve grid adaptation, a cell-based octree data structure~{\cite{khokhlov1998fully}} is incorporated within SYN3D. When a parent cell is flagged for refinement, eight children cells are generated and a continuous sequence of memory is allocated, such that only the pointer to the first child is stored for the parent cell. The following information is stored in each cell for the data structure:
\begin{itemize}
    \item level of the cell in the tree;
    \item pointer to the parent cell;
    \item pointer to the first child;
    \item pointers to 6 neighbor cells.
\end{itemize}
Given a list of refined cells from the error estimator, the connection pointers are automatically generated and the ranks are rebalanced for parallel efficiency. In the current flow solver only geometry variables exist on face and vertices while all flow variables are cell-center based.
\subsection{Integration process}
Time integration is achieved by using a five-stage Runge-Kutta method. Within each stage, the time increment updates are performed for all levels of cells recursively. The time integration follows the process proposed in~{\cite{khokhlov1998fully}}, where the time steps in refined cells are automatically reduced from that of the parent cell,
\begin{equation}
    \Delta t(l) = 2^{l_{\text{max}}-l} \Delta t,
\end{equation}
where $\Delta t$ is the time step on the root cell and $l$ the level of the current cell. Assuming the flux evaluation and time advancement procedures at level $l$ to be $A(l)$, the operation at level $l$ is

\begin{equation}
    O(l)=\left\{
\begin{array}{lcl}
O(l+1)O(l+1)A(l);       &      & \text{if } l < l_{\text{max}};\\
A(l);   &      & \text{if } l = l_{\text{max}}.
\end{array} \right.
\end{equation}

\subsection{Periodic hill test case}
\label{sec:perhill}
\subsubsection{Introduction}
Our first test case is the periodic hill, whose geometry and grid are shown in figure~\ref{figure:perhill_geometry}, with periodic boundary conditions in both streamwise and spanwise directions, and no-slip conditions are applied at the upper and lower boundaries. The Reynolds number based on the hill height and mean bulk velocity at the hill crest is $Re = 10600$. A pressure forcing term~\cite{jiang2020assessment} is added to the streamwise momentum and energy equations in order to drive the flow to maintain a constant Reynolds number during the simulation. The flow is highly unsteady featuring separation from the continuous surface and the separation point oscillates over a large range on the wall. The mean flow is characterized by a separation bubble with an established separation and reattachment point. The form of the bubble depends on two factors: the position of the separation point and the turbulent intensity on the top of the bubble, which determines the level of energy exchange from the mean flow to the bubble region. The experimental data~\cite{breuer2009flow}, an LES with wall function on a fine grid of 4.6M~\cite{temmerman2003investigation} and a very fine wall-resolved LES on 13.1M grid points~\cite{rapp2011flow} serve as the reference data.
    
Two levels of grids are available for the study. The coarse grid has $ 160\times 160\times 64 $ grid points, while the fine grid is refined from the coarse grid only in the spanwise direction with $ 160\times 160\times 128 $ grid points. The first layer of cells satisfies $y^+\approx 1$, which ensures wall-resolved LES on both the upper and lower surfaces. The goal of the study is to apply the error estimators to the coarse grid based on the temporal and spanwise spatially averaged solution and refine only $5\%$ of cells with the largest error. We have selected $5\%$ as it almost doubles the computational cost of the course grid due to the addition of the refined grids as well as halving the time step in the refined regions.
    \begin{figure}[hbt!]
        \centering
        \begin{subfigure}[b]{0.49\linewidth}            
            \includegraphics[clip=true, trim= 0.1cm 0.1cm 0.1cm 0.1cm, width=\textwidth]{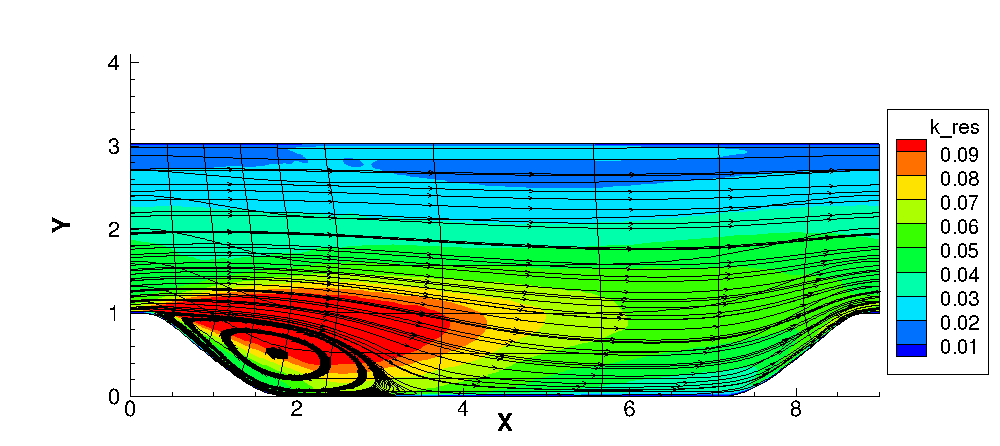}
            \label{fig:bubble_coarse}
            \caption{}
        \end{subfigure}
        \begin{subfigure}[b]{0.49\linewidth}            
            \includegraphics[clip=true, trim= 0.1cm 0.1cm 0.1cm 0.1cm, width=\textwidth]{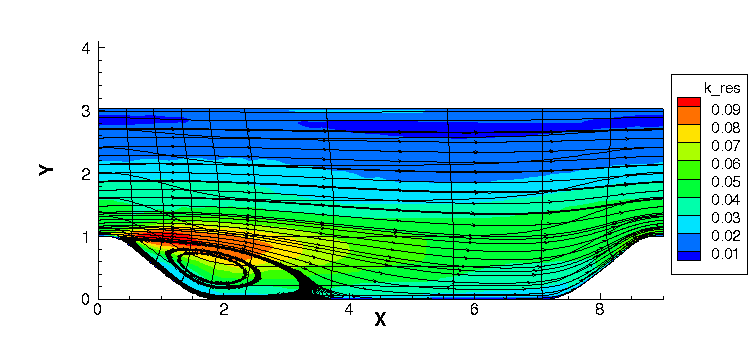}
            \label{fig:bubble_fine}
            \caption{}
        \end{subfigure}
        \caption{Separation bubble and TKE contour on (a): coarse grid; (b): fine grid.}\label{fig:bubble_org}
    \end{figure}

    \begin{figure}[hbt!]
        \begin{center}
            \includegraphics[width=0.49\linewidth]{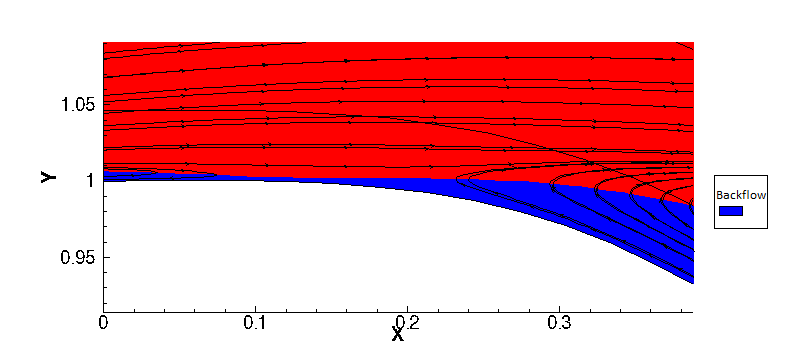}
            \includegraphics[width=0.49\linewidth]{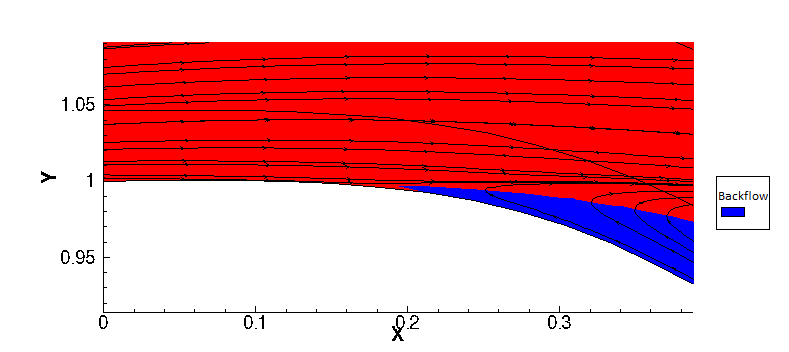}
        \end{center}
        \caption{Backflow on the hill top region, left: coarse grid; right: fine grid.}
        \label{figure:perhill_backflow} 
    \end{figure}
    
    \subsubsection{Flow solution}
    The simulations were launched on coarse and fine grids for 100 flow through periods, and the flow field is averaged over 30 flow through periods in time and the spanwise direction. Fig.~{\ref{fig:bubble_org}}~(a) shows the size of the separation bubble with the contour of the captured TKE. The coarse grid shows an early flow reattachment due to an over-estimation of the TKE level in the region above the separation bubble, the phenomenon is also reported in~{\cite{temmerman2003investigation}} for a coarse grid. A refinement in the spanwise direction allows for a lower TKE level and a better capture of the bubble length as shown in Fig.~{\ref{fig:bubble_org}}~(b). Fig.~\ref{figure:perhill_backflow} provides further evidence that the coarse grid incorrectly captures a backflow layer on the top of the hill which is removed in the fine grid through spanwise refinement. Table~{\ref{tab:adapt_separation}} confirms that the coarse grid incorrectly predicts the separation point by estimating the separation prior to the top of the hill at a negative $ X $ coordinate and predicts an early reattachment of the flow, while the fine grid shows good agreement with the reference LES data, even with a comparatively smaller grid size.
    
    \begin{table}
        \caption{Separation and reattachment points for the periodic hill case.}
        \centering
        \begin{tabular}{c l l l l}
            \hline
            Grid & Size & $ (x/h)_{sep} $ & $ (x/h)_{reat} $ & $ L_{bubble} $\\
            \hline
            Temmerman et al.\cite{temmerman2003investigation} & 4.6M & 0.22 & 4.72 & 4.5\\
            Breuer et al.\cite{breuer2009flow} & 13.1M & 0.19 & 4.69 & 4.5\\
            Coarse & 1.6M & -0.93 & 4.06 & 4.99\\
            Fine & 3.3M & 0.22 & 4.67 & 4.45\\
            $IQ_{\eta-emp}$ adapted & 2.1M & 0.14 & 3.95 & 3.81\\
            $IQ_{k-emp}$ adapted & 2.1M & 0.21 & 4.41 & 4.2\\
            $IQ_{k-ke}$ adapted & 2.1M & -0.93 & 4.99 & 5.92\\
            $IQ_{k-tke}$ adapted & 2.1M & 0.22 & 4.55 & 4.35\\
            \hline
        \end{tabular}
        \label{tab:adapt_separation}
    \end{table}
\subsubsection{Evaluation of numerical TKE}
    \begin{figure}[hbt!]
        \begin{center}
            \includegraphics[clip=true, trim= 0.1cm 0.1cm 0.1cm 0.1cm, width=.49\linewidth]{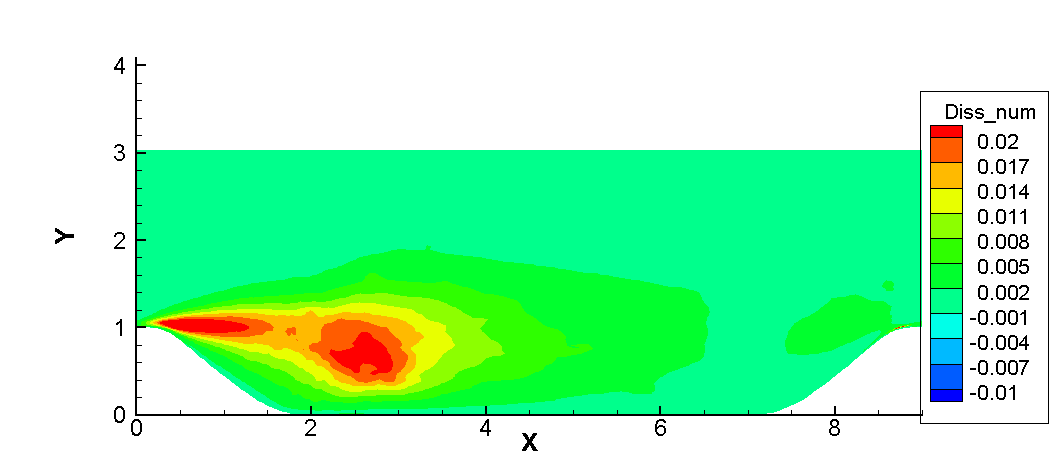}
            \includegraphics[clip=true, trim= 0.1cm 0.1cm 0.1cm 0.1cm, width=.49\linewidth]{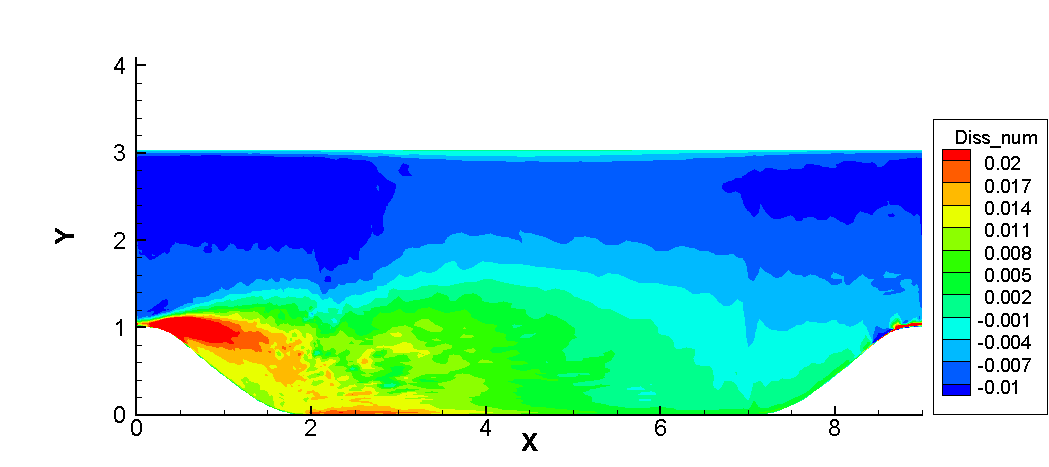}
            \includegraphics[clip=true, trim= 0.1cm 0.1cm 0.1cm 0.1cm, width=.49\linewidth]{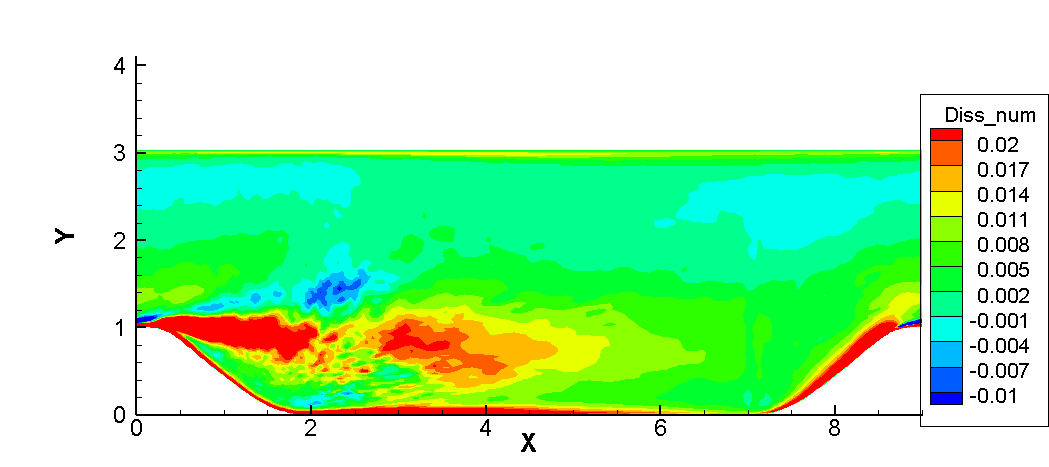}
        \end{center}
        \caption{Numerical dissipation, $\epsilon_n$, based on empirical approach (top left), KE-based approach (top right) and TKE-based approach (bottom).}
        \label{figure:num_diss}
    \end{figure}

   \begin{figure}[hbt!]
        \centering
        \begin{subfigure}[b]{0.49\textwidth}            
            \includegraphics[clip=true, trim= 0.1cm 0.1cm 0.1cm 0.1cm, width=\textwidth]{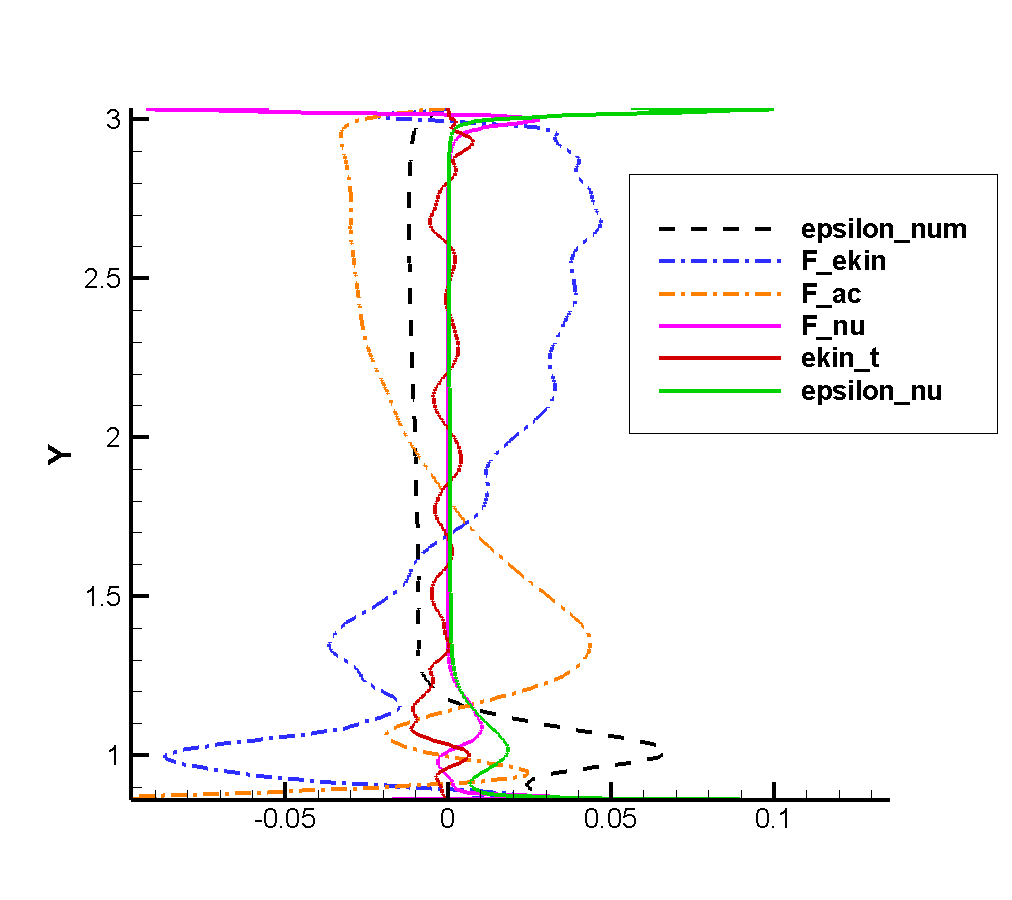}
            \caption{}
        \end{subfigure}
        \begin{subfigure}[b]{0.49\textwidth}            
            \includegraphics[clip=true, trim= 0.1cm 0.1cm 0.1cm 0.1cm, width=\textwidth]{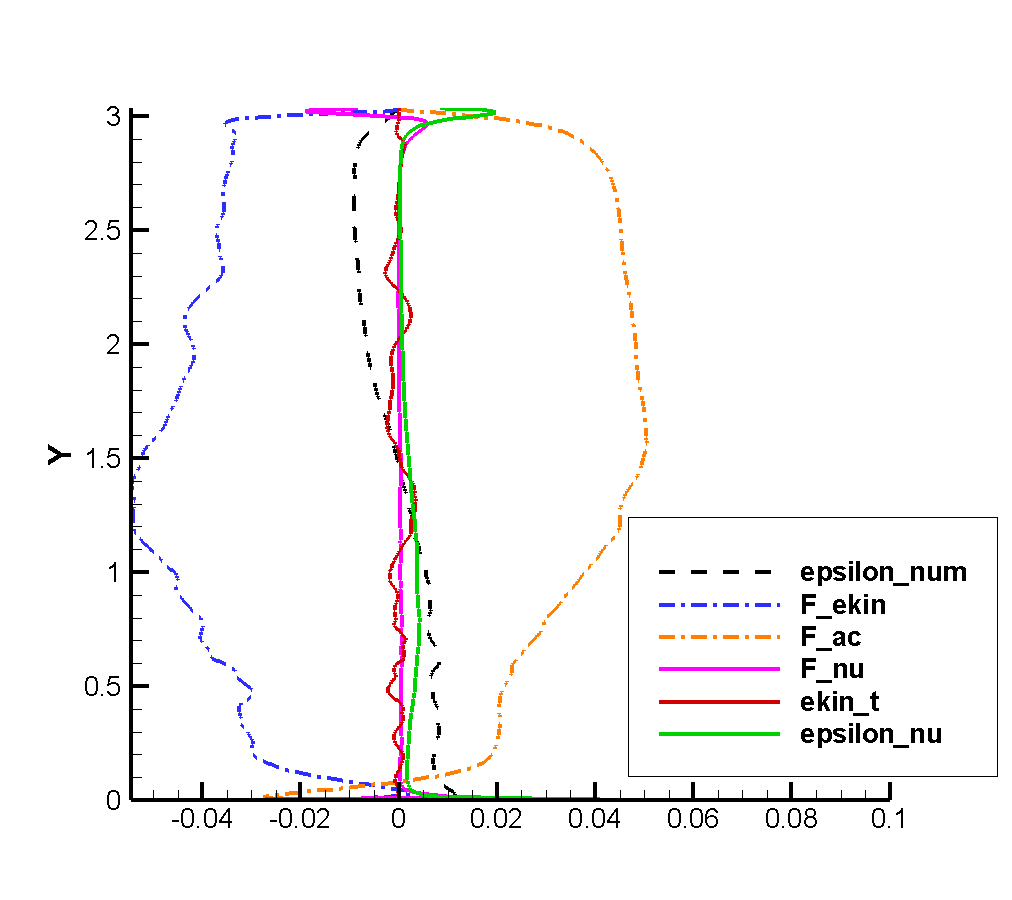}
            \caption{}
        \end{subfigure}
        \caption{Profile of contributions from different terms to KE numerical dissipation {\color{black}(Eq.~\ref{eq:diss_fv}, time-averaged)} at: (a) $X=0.5$; (b) $X=4$.}
        \label{fig:comp_ke}
    \end{figure}
    
   \begin{figure}[hbt!]
        \centering
        \begin{subfigure}[b]{0.49\textwidth}            
            \includegraphics[clip=true, trim= 0.1cm 0.1cm 0.1cm 0.1cm, width=\textwidth]{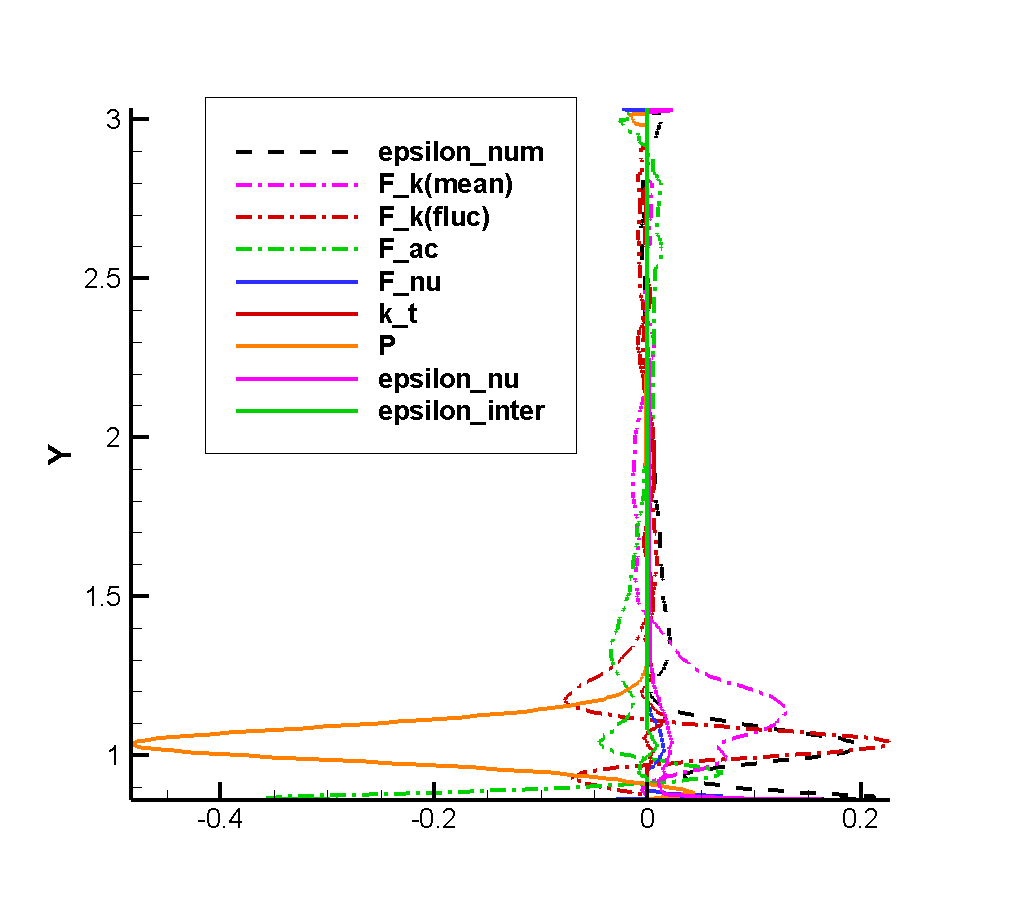}
            \caption{}
        \end{subfigure}
        \begin{subfigure}[b]{0.49\textwidth}            
            \includegraphics[clip=true, trim= 0.1cm 0.1cm 0.1cm 0.1cm, width=\textwidth]{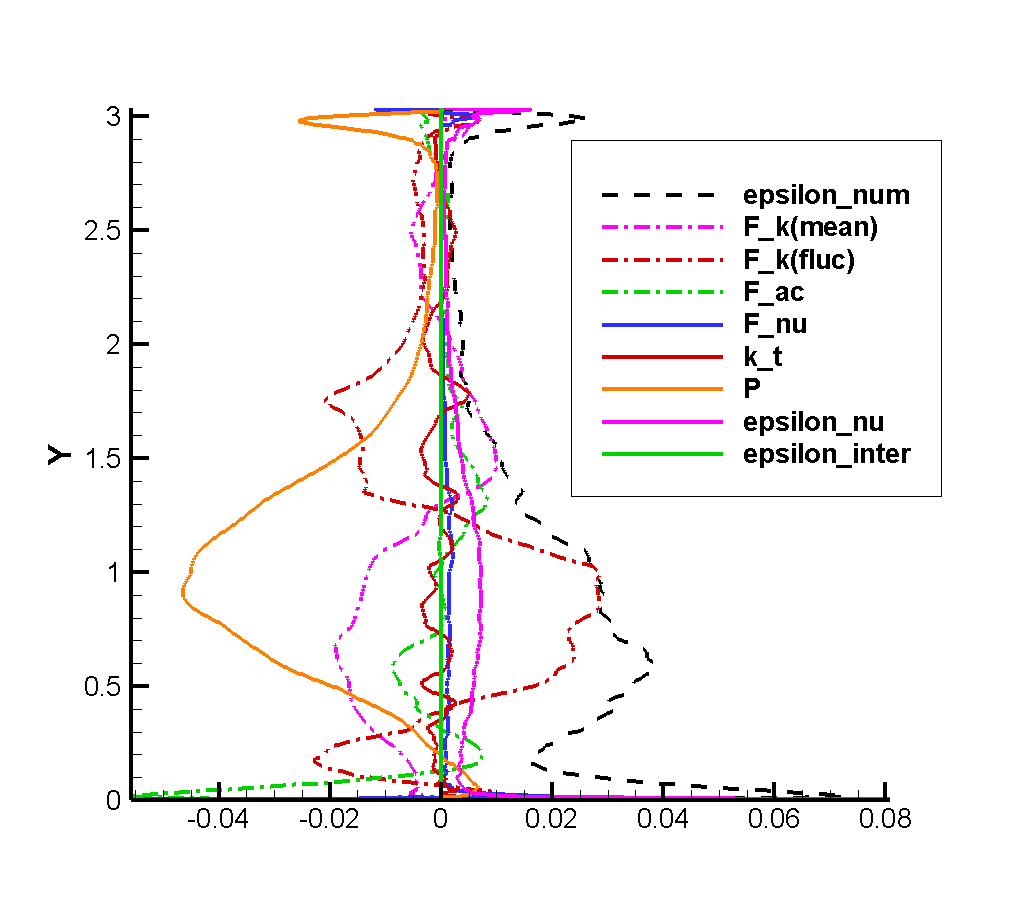}
            \caption{}
        \end{subfigure}
        \caption{Profile of contributions from different terms to TKE numerical dissipation {\color{black}(Eq.~\ref{eq:diss_tke_fv}, time-averaged)} at: (a) $X=0.5$; (b) $X=4$.}
        \label{fig:comp_tke}
    \end{figure}
    
    \begin{figure}[hbt!]
        \begin{center}
            \includegraphics[clip=true, trim= 0.1cm 0.1cm 0.1cm 0.1cm, width=.49\linewidth]{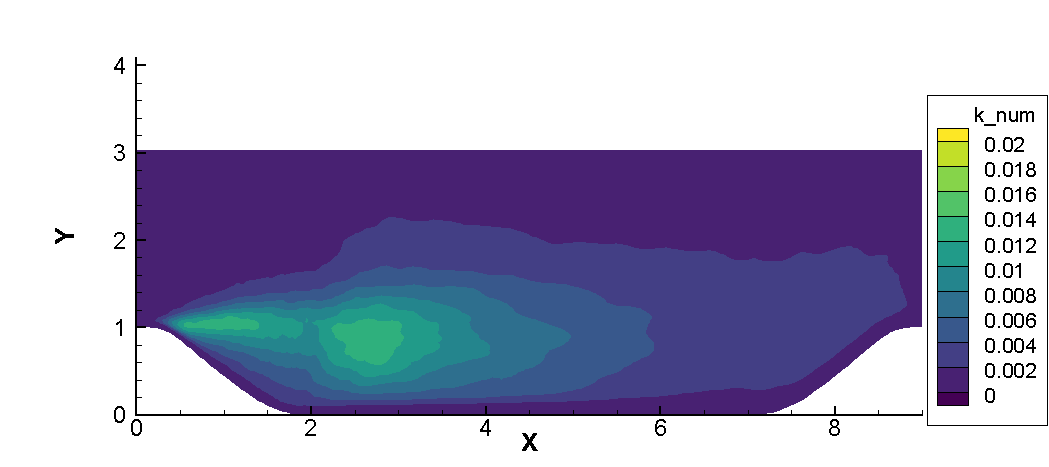}
            \includegraphics[clip=true, trim= 0.1cm 0.1cm 0.1cm 0.1cm, width=.49\linewidth]{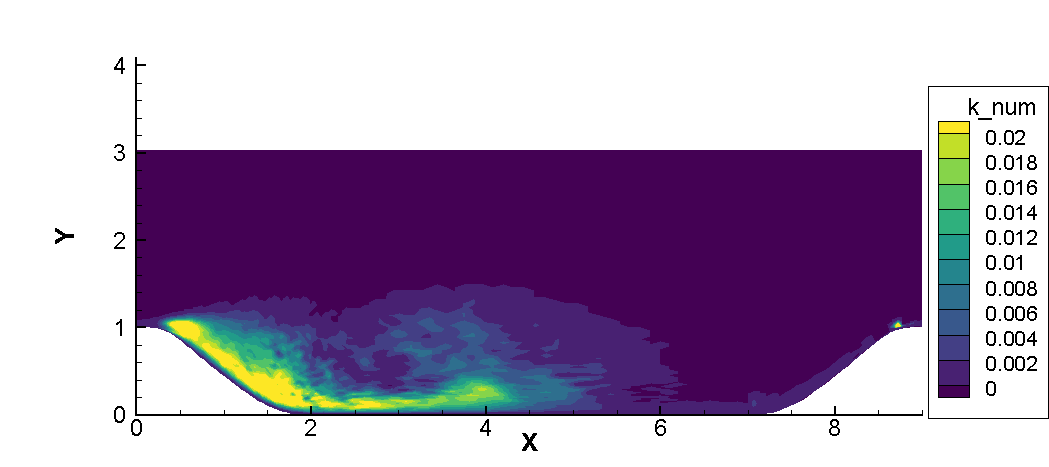}
            \includegraphics[clip=true, trim= 0.1cm 0.1cm 0.1cm 0.1cm, width=.49\linewidth]{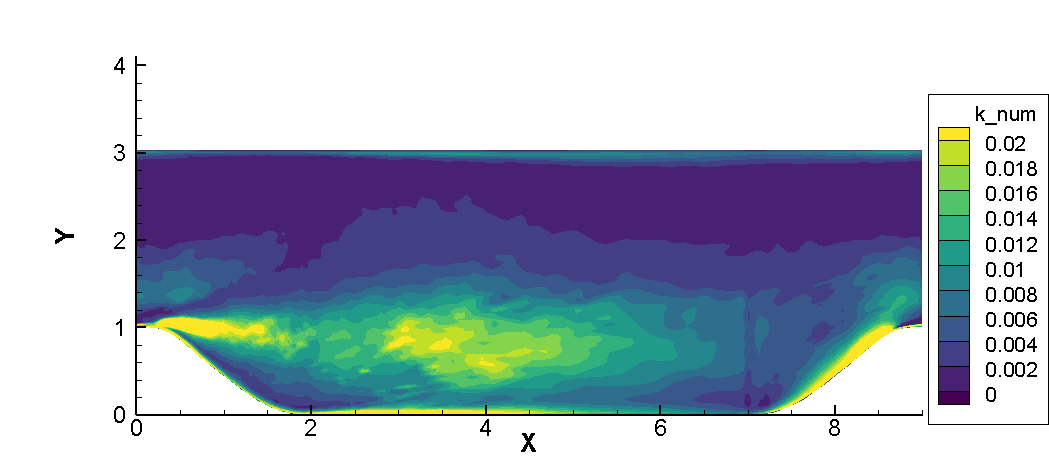}
        \end{center}
        \caption{Numerical TKE, $k_{\text{num}}$, based on empirical formula (top left), KE numerical dissipation (top right) and TKE numerical dissipation (bottom).}
        \label{figure:k_num}
    \end{figure}

    \begin{figure}[hbt!]
        \centering
        \begin{subfigure}[b]{0.49\textwidth}            
            \includegraphics[clip=true, trim= 0.1cm 0.1cm 0.1cm 0.1cm, width=\textwidth]{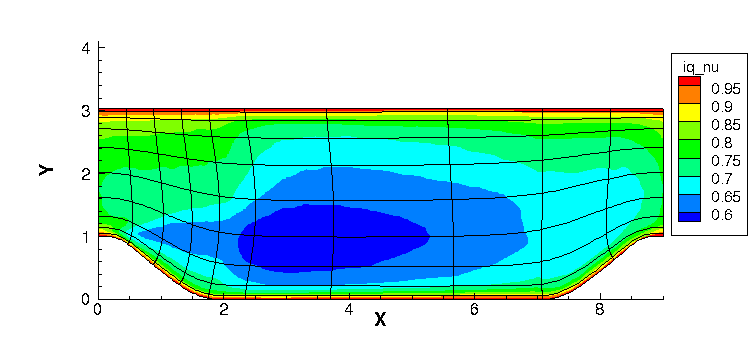}
            \caption{}
        \end{subfigure}
        \begin{subfigure}[b]{0.49\textwidth}            
            \includegraphics[clip=true, trim= 0.1cm 0.1cm 0.1cm 0.1cm, width=\textwidth]{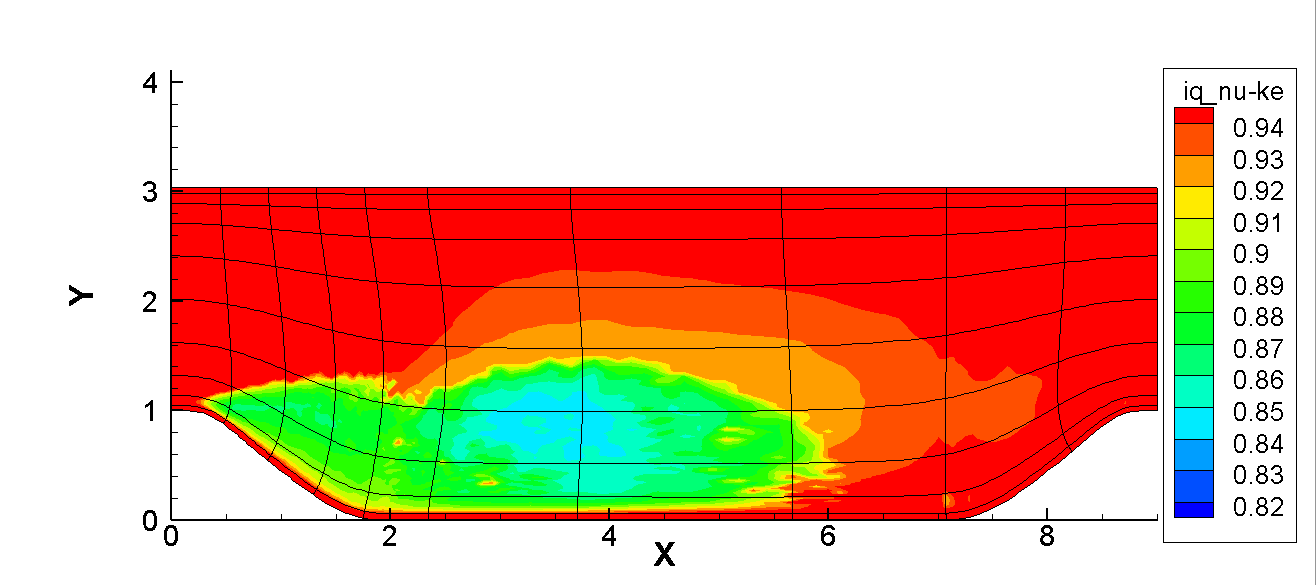}
            \caption{}
        \end{subfigure}
        \begin{subfigure}[b]{0.49\textwidth}            
            \includegraphics[clip=true, trim= 0.1cm 0.1cm 0.1cm 0.1cm, width=\textwidth]{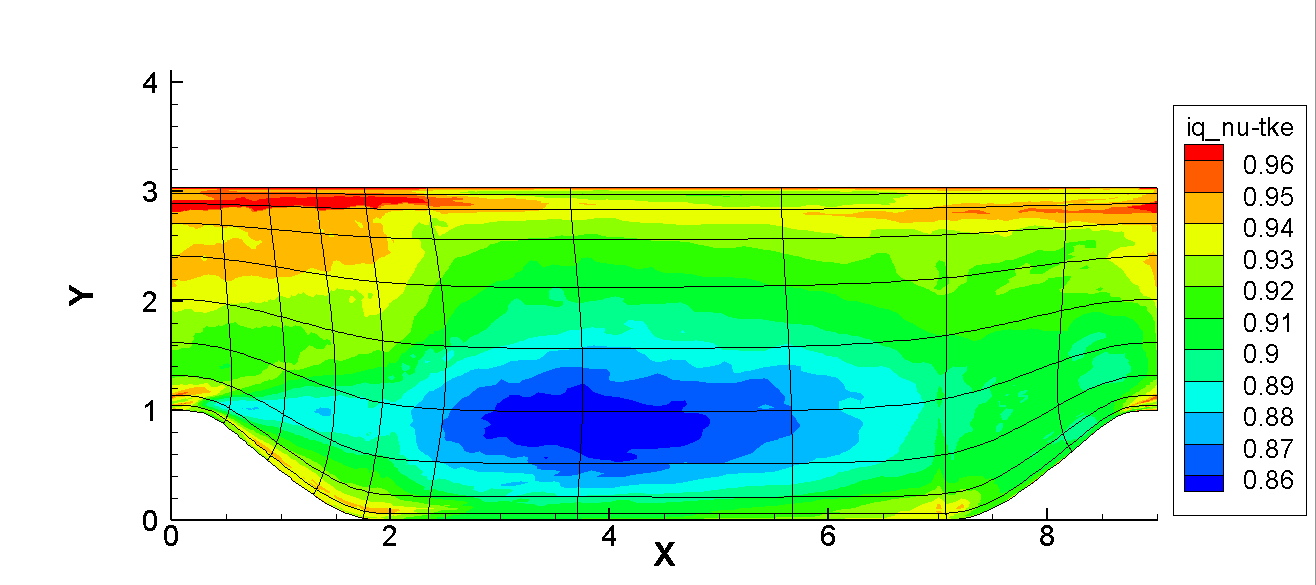}
            \caption{}
        \end{subfigure}
        \caption{Error estimation on coarse grid based on $IQ_{\nu}$: (a) $IQ_{\nu-emp}$; (b) $ IQ_{\nu-ke} $; (c) $IQ_{\nu-tke}$.}\label{fig:iq_nu_3ways}
    \end{figure}

   \begin{figure}[hbt!]
        \centering
        \begin{subfigure}[b]{0.49\textwidth}            
            \includegraphics[clip=true, trim= 0.1cm 0.1cm 0.1cm 0.1cm, width=\textwidth]{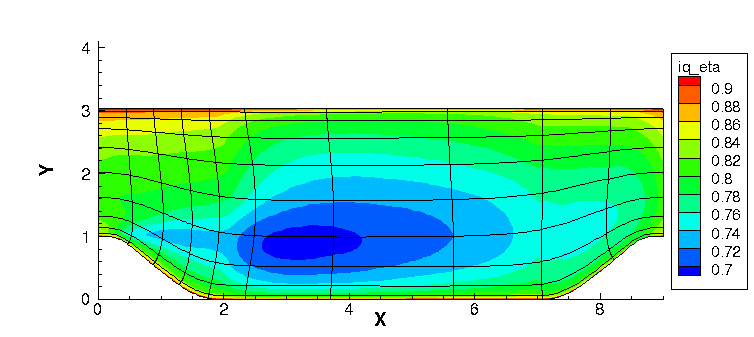}
            \caption{}
        \end{subfigure}
        \begin{subfigure}[b]{0.49\textwidth}            
            \includegraphics[clip=true, trim= 0.1cm 0.1cm 0.1cm 0.1cm, width=\textwidth]{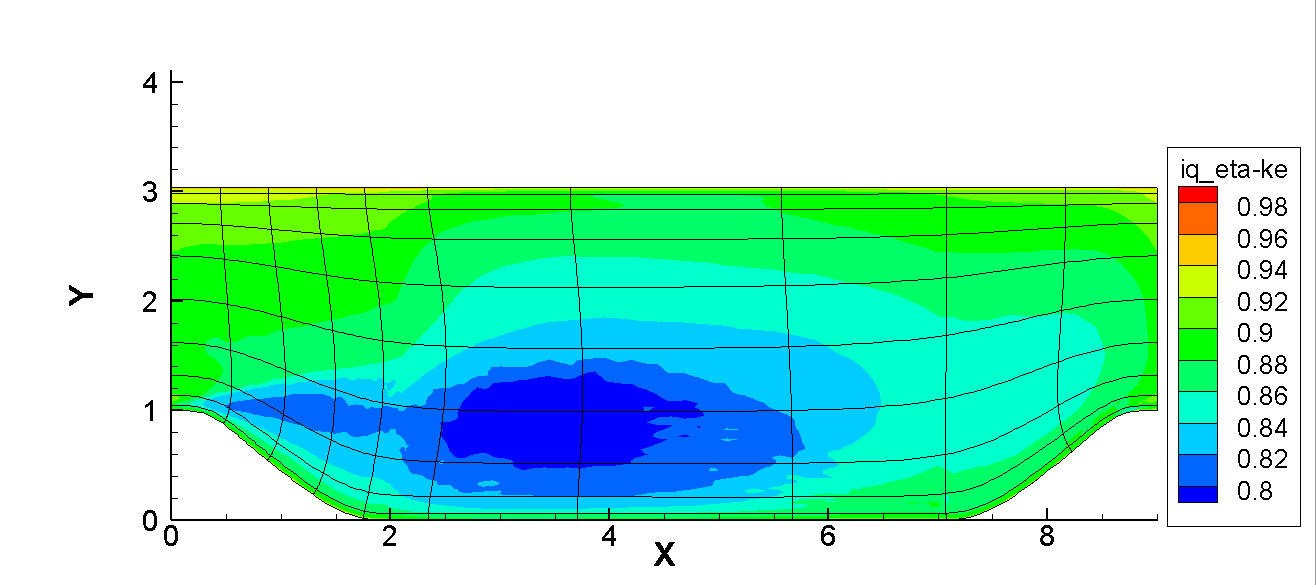}
            \caption{}
        \end{subfigure}
        \begin{subfigure}[b]{0.49\textwidth}            
            \includegraphics[clip=true, trim= 0.1cm 0.1cm 0.1cm 0.1cm, width=\textwidth]{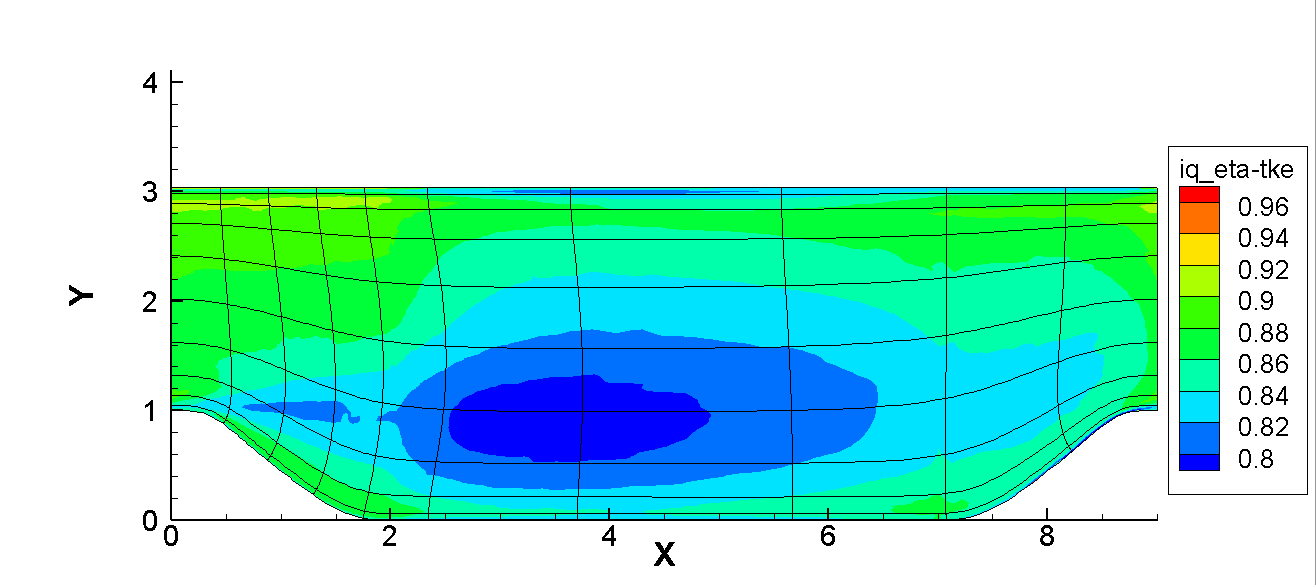}
            \caption{}
        \end{subfigure}
        \caption{Error estimation on coarse grid based on $IQ_{\eta}$: (a) $IQ_{\eta-emp}$; (b) $ IQ_{\eta-ke} $; (c) $IQ_{\eta-tke}$.}\label{fig:iq_eta_3ways}
    \end{figure}
    
    \begin{figure}[hbt!]
        \centering
        \begin{subfigure}[b]{0.49\textwidth}            
            \includegraphics[clip=true, trim= 0.1cm 0.1cm 0.1cm 0.1cm, width=\textwidth]{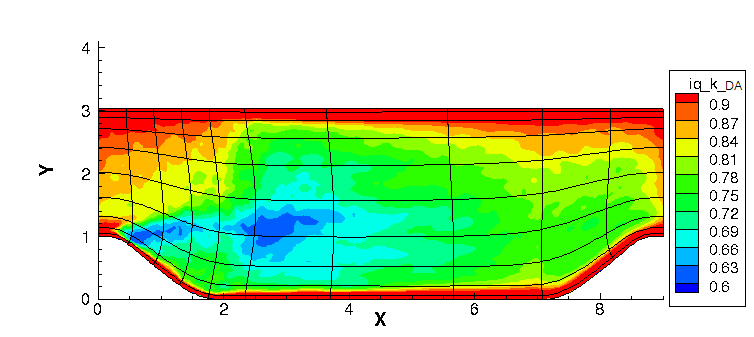}
            \caption{}
        \end{subfigure}
        \begin{subfigure}[b]{0.49\textwidth}            
            \includegraphics[clip=true, trim= 0.1cm 0.1cm 0.1cm 0.1cm, width=\textwidth]{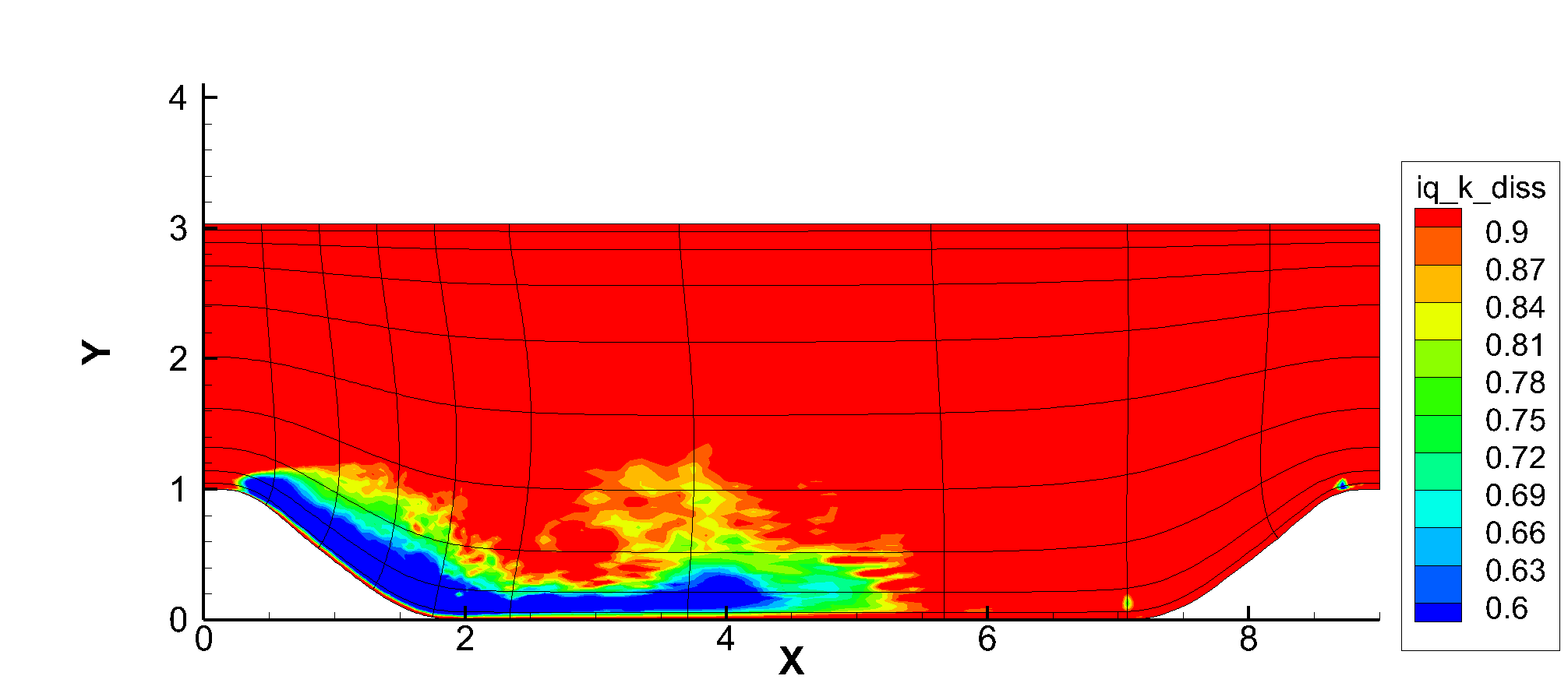}
            \caption{}
        \end{subfigure}
        \begin{subfigure}[b]{0.49\textwidth}            
            \includegraphics[clip=true, trim= 0.1cm 0.1cm 0.1cm 0.1cm, width=\textwidth]{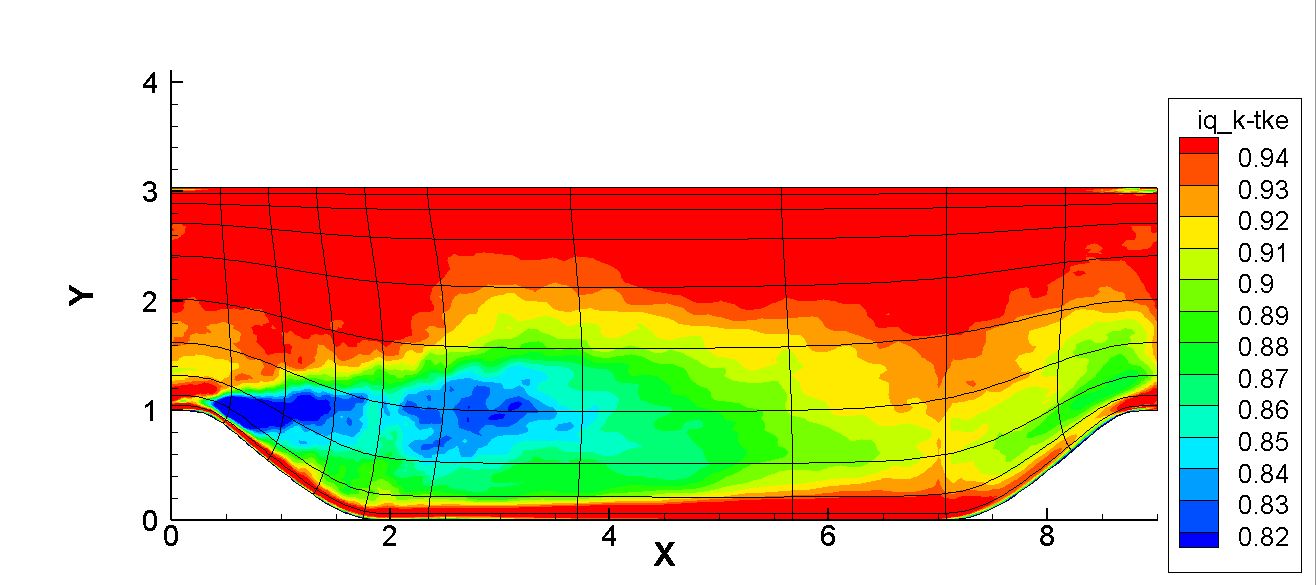}
            \caption{}
        \end{subfigure}
        \caption{Error estimation on coarse grid based on $IQ_{k}$: (a) $IQ_{k-emp}$; (b) $IQ_{k-ke}$; (c) $IQ_{k-tke}$.}\label{fig:error}
    \end{figure}
The numerical TKE $k_{\text{num}}$ is evaluated using all three approaches presented in Section~\ref{sec:num_tke}. First, the {\color{black} time-averaged} numerical dissipation, $\bar{\epsilon}_n$, is evaluated and the results are shown in Fig.~\ref{figure:num_diss}. The results using the empirical formula, shown in the top-left-hand corner, only highlights the mixing layer above the separation bubble. The {\color{black}KE-based approach} of Schranner et al.~\cite{schranner2015assessing} highlights the intensively turbulent region near the separation point, while not targeting the mixing layer between $2<X<3.2$. The approach also shows high dissipation in a thin layer near the wall within the separation bubble $2<X<4$ and in a small near-wall region on the upstream side of the hill top starting from $X \approx 8.5$, where flow separation starts to appear. Generally, the region between the hills with $Y < 1$ shows positive numerical dissipation. Negative dissipation is shown in the region above the hill for $Y>1$, where the convective flux dominates and the turbulence intensity is low, as discussed in Section~\ref{sec:ke_diss}. A similar trend was reported by Castiglioni and Domaradzki~\cite{castiglioni2015numerical} who observed negative numerical dissipation in laminar regions far from the wall. {\color{black} In such regions a combination of low velocity gradients and an absences of a turbulent shear layer lead to a diminished role of the dissipation term, $\epsilon_\nu$ from Eq.~\ref{eq:diss_fv}. In addition, following Eq.~\ref{eq:ratio} and as reported by~\cite{castiglioni2015numerical} we observe a division of two small quantities ($\bar{\epsilon}_n$ and $\overline{\tau_{ij}\frac{\partial u_i}{\partial x_j}}$) which lead to unreasonable values of the numerical viscosity. Moreover the convective and pressure flux terms from Eq.~\ref{eq:diss_fv} dominate over the turbulence dissipation and this eventually results in an estimation of negative numerical dissipation. Thus an application of Eq.~\ref{eq:pos_diss} eliminates the impacted control volumes from adaptation. Lastly, we examine the developed TKE numerical dissipation from Fig.~\ref{figure:num_diss}~(bottom). A quick survey of the contour and a comparison against the KE-based approach, reveals the removal of the negative numerical dissipation in the laminar region and an emphasis along the boundary layer and separated regions. As we traverse down towards the wall, the numerical dissipation begins to show higher value when approaching the mixing layer at the hill top level at $Y \approx 1$, where the production term, $P$, and the convective term, $F_k$, overwhelms the remaining contributions. As we approach the lower wall, the boundary layer is particularly highlighted, where the viscous and turbulent transport terms, $F_{\nu}$, as well as the dissipation terms, $\epsilon_{\nu}$ and $\epsilon_{inter}$, prevail due to high mean velocity gradients.}

{\color{black} To provide a more thorough comparison between the {\color{black}KE-and TKE-based approaches}, we plot and compare in Figs.~\ref{fig:comp_ke} and~\ref{fig:comp_tke} the contributions from all the terms in Eqs.~\ref{eq:diss_fv} and~\ref{eq:diss_tke_fv} at two streamwise locations: at $X=0.5$ shortly after the flow separation and $X=4$ immediately before the flow reattachment. From Fig.~\ref{fig:comp_ke}~(a) at $X=0.5$ as we trace from the lower to the upper wall, in the shear layer region of $Y\approx1$ the KE numerical dissipation, $\epsilon_n$ reaches its maximum value, driven primarily by the convective term. As we move towards the upper wall the KE numerical dissipation is characterized by the convective flux term as well as work due to pressure which results in a constant weak negative value. At $X=4$, the KE numerical dissipation is dominated by the convective and pressure work terms which results in positive values in the lower turbulent region but negative in the upper largely laminar region. Fig.~\ref{fig:comp_tke} shows the TKE numerical dissipation profile and the contributions from different items in Eq.~\ref{eq:diss_tke_fv} at the same positions of $X=0.5$ and $X=4$. It is observed in Fig.~\ref{fig:comp_tke}~(a) that at $X=0.5$ the lower wall shows high TKE numerical dissipation contributed by the convection and production of TKE. The remaining terms play a minor role and all terms diminish as we move towards the upper wall. At the downstream location, $X=4$, an equivalent observation is noted; however, the trends encompass a larger area. Comparing the two Figs.~\ref{fig:comp_ke} and~\ref{fig:comp_tke} the TKE based approach emphasizes the separated region. While the KE {\color{black} numerical dissipation} is driven by the convection and pressure terms, the TKE approach is driven by the production and convection of TKE with a larger concentration on the boundary layer region.}
    
{\color{black} Having established clear differences between the {\color{black}KE- and TKE-based approaches}, we now center our discussion on the numerical TKE $k_{\text{num}}$ that is evaluated based on the {\color{black}empirical approach~(Eq.~\ref{k_num}), KE-based approach~(Eq.~\ref{eq:diss_fv}), and TKE-based approach~(Eq.~\ref{eq:diss_tke_fv})} and we compare their contours in Fig.~\ref{figure:k_num}. As Eq.~\ref{k_num} demonstrates, the empirical approach is based on the modeled TKE and as such highlights the highly turbulent region, while diminishes as it approaches the wall. The empirical based formula only shows high value in the mixing layer above the separation bubble with a peak value of $0.014$, while the KE-based approach mainly targets the near-wall region in the separation bubble with a higher peak value of $0.044$. The comparatively low numerical viscosity in the mixing layer leads to a low numerical TKE since the formula in Eq.~\ref{eq:ratio} relates the numerical eddy viscosity, $\nu_{\text{num}}$, to the ratio between $\bar{\epsilon}_n$ and $\overline{\tau_{ij}\frac{\partial u_i}{\partial x_j}}$. Although the mixing layer figures elevated values of numerical dissipation $\bar{\epsilon}_n$; $\overline{\tau_{ij}\frac{\partial u_i}{\partial x_j}}$ also shows high values in this turbulent region, resulting in largely low values of $k_{\text{num}}$. The TKE-based approach is able to target both the mixing layer and the boundary layer with a peak value of $0.013$. The targeted boundary layer includes three domains, the separation point region, the reattachment point region and a small near-wall region on the upstream side of the hill top region with $7.5< X < 8.6$. The targeted regions are well aligned with the area with high value of numerical dissipation, $\bar{\epsilon}_n$, as shown in Fig.~\ref{figure:num_diss}.}
    
    \begin{figure}[hbt!]
        \begin{center}
            \centering
            \begin{subfigure}[b]{0.49\linewidth}    
                \includegraphics[clip=true, trim= 1.0cm 0.5cm 1.0cm 1.0cm, width=\textwidth]{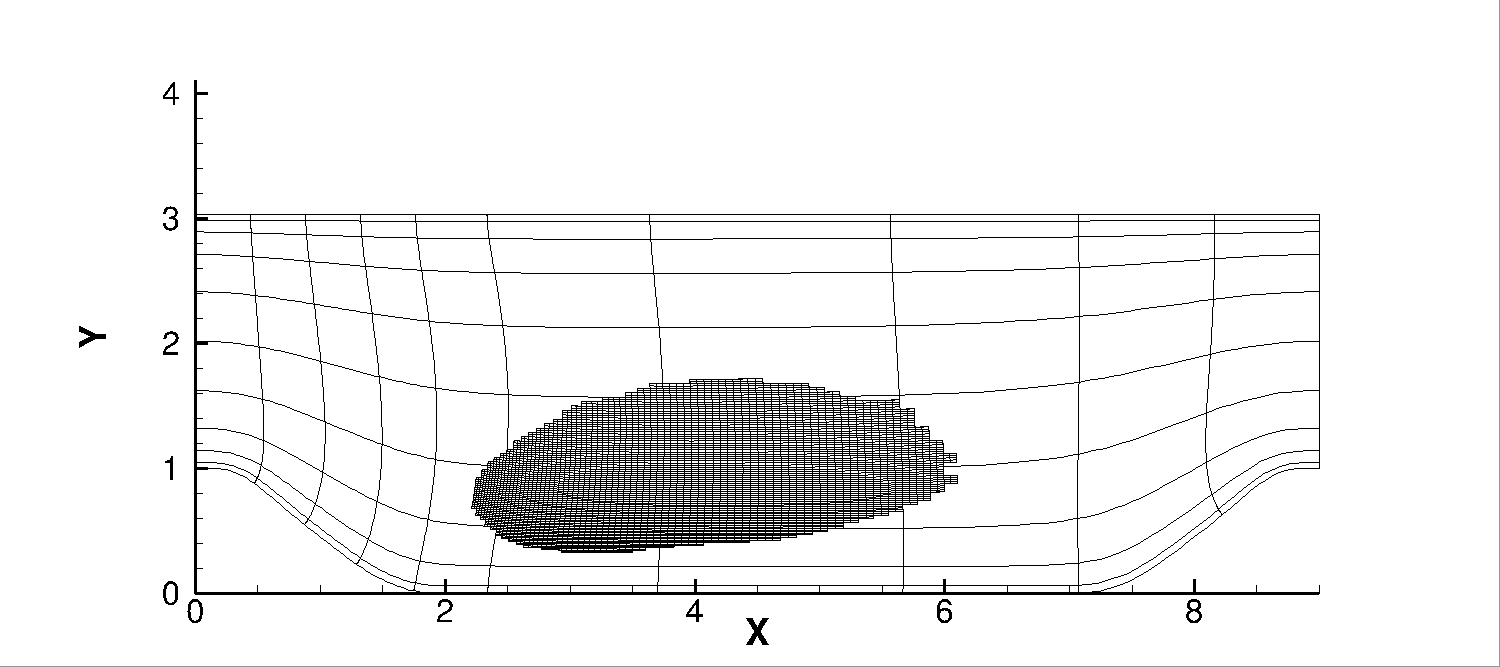}
                \caption{}
            \end{subfigure}
            \begin{subfigure}[b]{0.49\linewidth}            \includegraphics[clip=true, trim= 0.2cm 0.1cm 0.2cm 0.1cm,width=\textwidth]{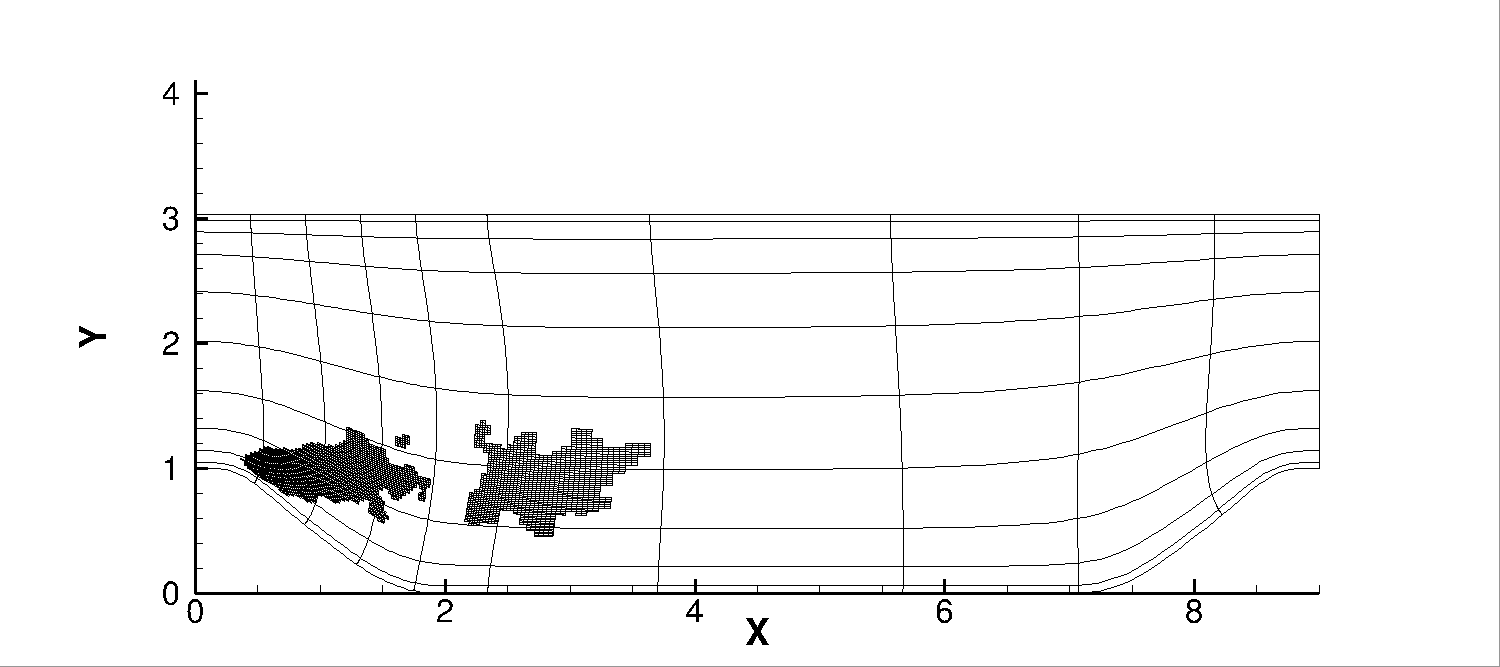}
                \caption{}
            \end{subfigure}
            \begin{subfigure}[b]{0.49\linewidth}            \includegraphics[clip=true, trim= 0.2cm 0.1cm 0.2cm 0.1cm,width=\textwidth]{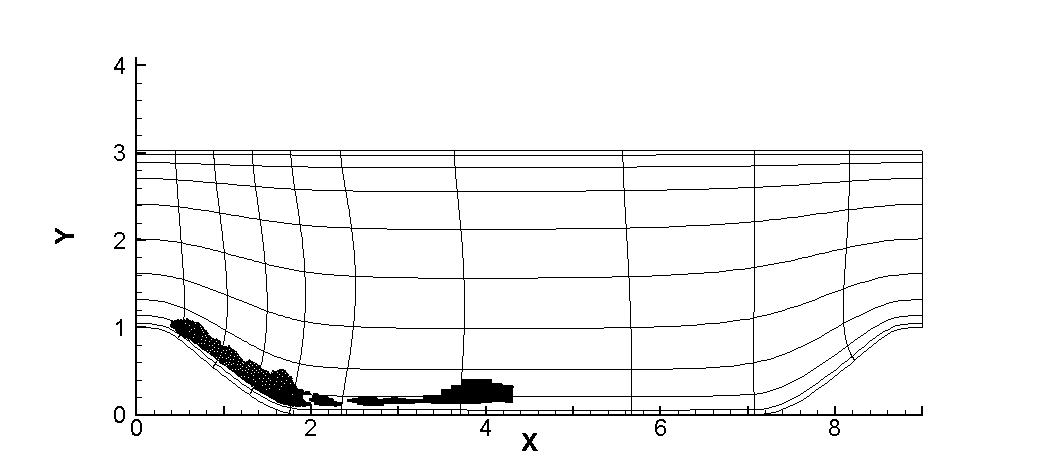}
                \caption{}
            \end{subfigure}
            \begin{subfigure}[b]{0.49\linewidth}            \includegraphics[clip=true, trim= 0.2cm 0.1cm 0.2cm 0.1cm,width=\textwidth]{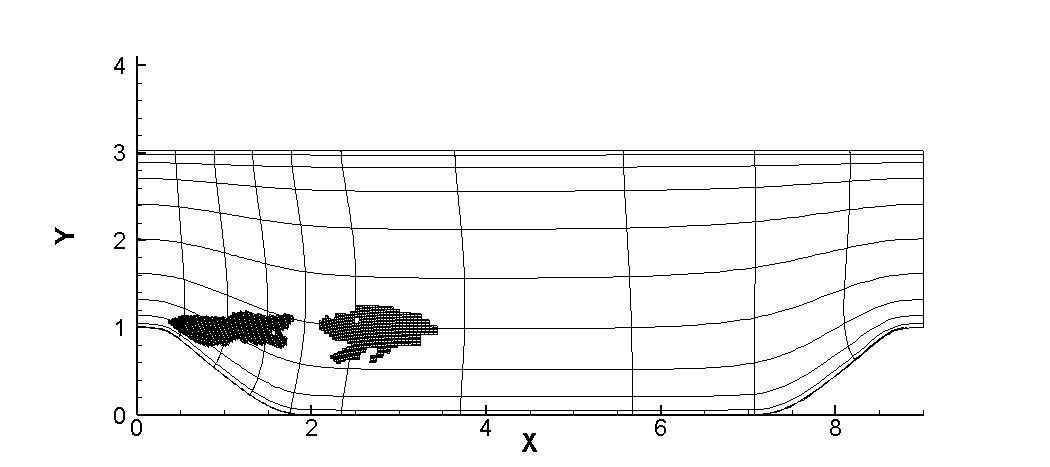}
                \caption{}
            \end{subfigure}            
        \end{center}
        \caption{(a): $ IQ_{\eta}$ adapted grid; (b): $ IQ_{k-emp}$ adapted grid; (c): $ IQ_{k-ke}$ adapted grid; (d): $ IQ_{k-tke}$ adapted grid.}
        \label{figure:mesh_adapt}
    \end{figure}
    
\subsubsection{Error estimation}
{\color{black} In this subsection, we employ the three presented numerical TKE evaluation approaches within the Index Quality based error estimators: $ IQ_{\nu} $, $ IQ_{\eta} $ and $ IQ_{k} ${\color{black}, leading to nine error estimators namely $IQ_{\nu,\eta,k-emp}$, $IQ_{\nu,\eta,k-ke}$ and $IQ_{\nu,\eta,k-tke}$}. The error estimation based on the $IQ_{\nu} $ estimator are shown in Fig.~{\ref{fig:iq_nu_3ways}}. {\color{black}If we are to target the lowest quality region in the computational domain according to the error estimator}, then the central region of the computational domain as well as a small region of the mixing layer near the separation point will require refinement. Three subfigures based on the three versions of $k_{\text{num}}$ lead to minimal differences in  targeted regions. The $ IQ_{\nu} $ estimator is a function of the ratio between $\nu_{\text{eff}}$ and $\nu$ as shown in Eq.~\ref{iq_nu}, thus only highly turbulent regions where $\nu_{\text{eff}}$ has large values will be identified. Without normalization against proper length and velocity scales, neither the turbulent shear layer near the separation point nor the boundary layer is targeted. It should be mentioned that the $IQ_{\nu-ke}$ estimator shows very high value in most of the computational domain outside the separation bubble region.  {\color{black}These regions are mostly dominated by convection and the KE-based approach tends to estimate non-physical negative numerical dissipation, which is clipped by Eq.~\ref{eq:pos_diss}, leading to zero numerical TKE. As a consequence, following Eq.~\ref{eq:nu_eff} and~\ref{nu_num}, $IQ_{\nu-ke}$ will have no contribution from $\nu_{\text{num}}$ in those regions and thus show high values.} Next, we present the $IQ_{\eta}$ error estimation in Fig.~{\ref{fig:iq_eta_3ways}}. The estimator targets domains similar to that of $IQ_{\nu}$; moreover, all versions of $k_{\text{num}}$ lead to very similar results.

Lastly, we present the $IQ_k$ error estimator in Fig.~\ref{fig:error}. The indicator $ IQ_{k-emp} $ with the empirical formula of $k_{\text{num}}$ targets primarily the mixing layer above the separation bubble. The indicator $ IQ_{k-ke} $ based on the evaluation of KE numerical dissipation provides an unexpected outcome, by exhibiting high values (over 90\%) in most of the flow field and targets the near-wall region on the downstream side of the hill. Unlike $ IQ_{k-emp} $,  $ IQ_{k-ke} $ does not target the mixing layer that is present above the separation bubble as a low-quality region, since the region is dominated by high strain rate and shows a low numerical eddy viscosity, $\nu_{\text{num}}$, and numerical TKE, $k_{\text{num}}$. The indicator $ IQ_{k-tke} $ based on the TKE numerical dissipation targets both the mixing layer and the near-wall region, including the separation point and the reattachment point. Since the correct capture of TKE in those regions is essential for the correct prediction of the bubble length as shown in{~\cite{temmerman2003investigation}}, the indicator $IQ_{k-tke}$ is expected to efficiently aid grid adaptation for LES.}

    \begin{figure}[hbt!]
        \centering
        \begin{subfigure}[b]{0.49\linewidth}            
            \includegraphics[clip=true, trim= 0.1cm 0.1cm 0.1cm 0.1cm, width=\textwidth]{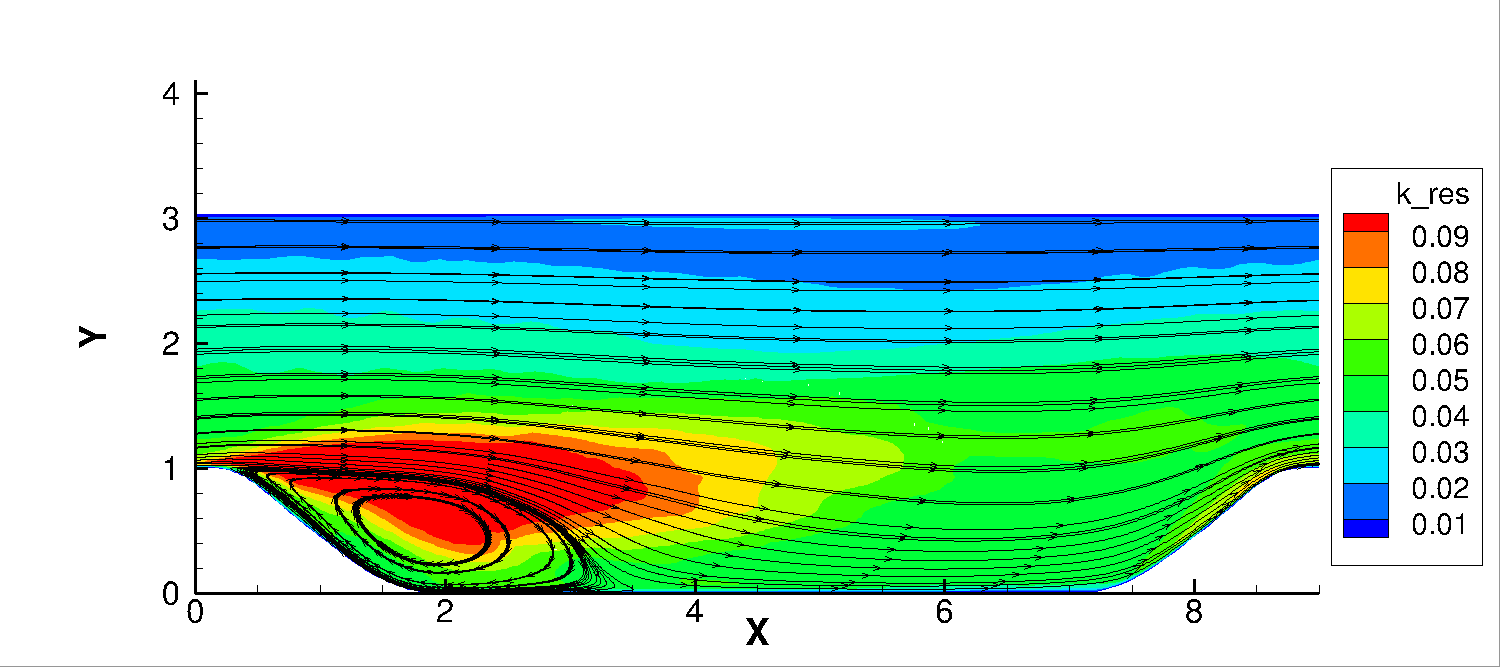}
            \caption{}
        \end{subfigure}
        \begin{subfigure}[b]{0.49\linewidth}
            \includegraphics[clip=true, trim= 0.1cm 0.1cm 0.1cm 0.1cm, width=\textwidth]{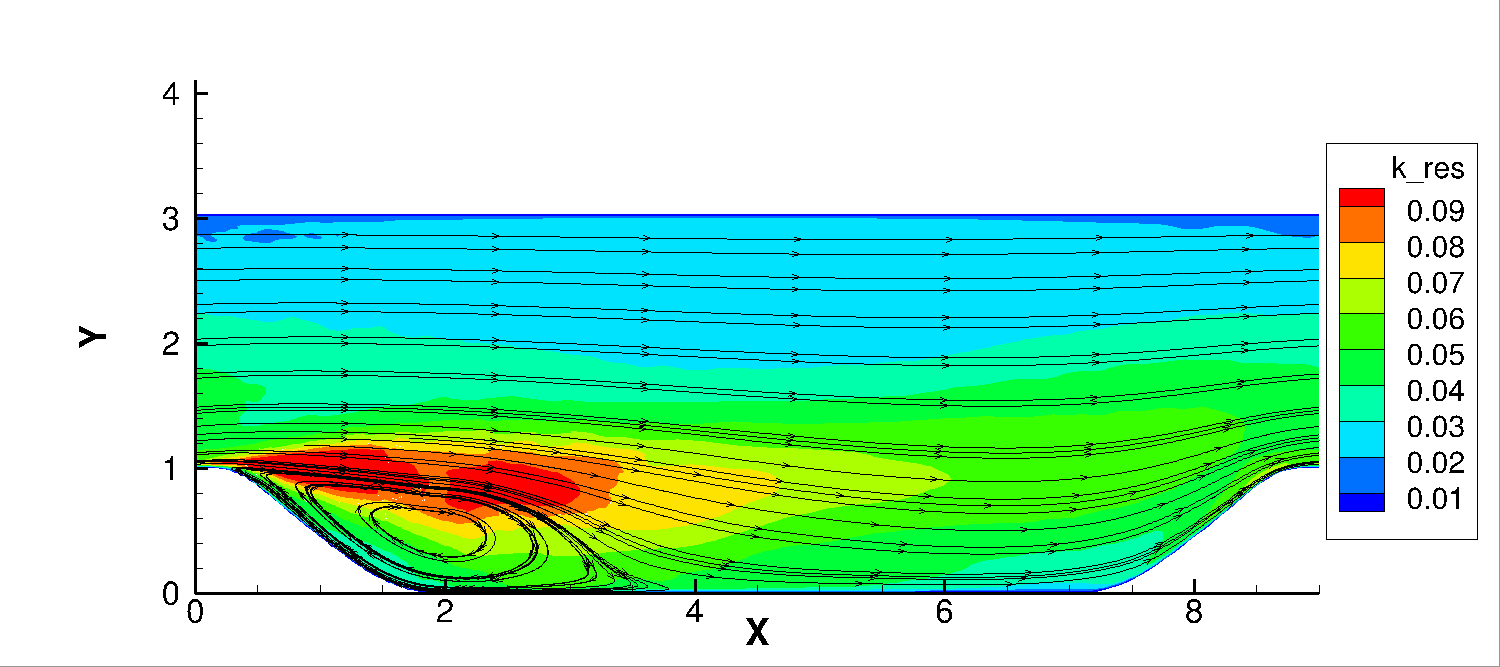}
            \caption{}
        \end{subfigure}
        \begin{subfigure}[b]{0.49\linewidth}
            \includegraphics[clip=true, trim= 0.1cm 0.1cm 0.1cm 0.1cm, width=\textwidth]{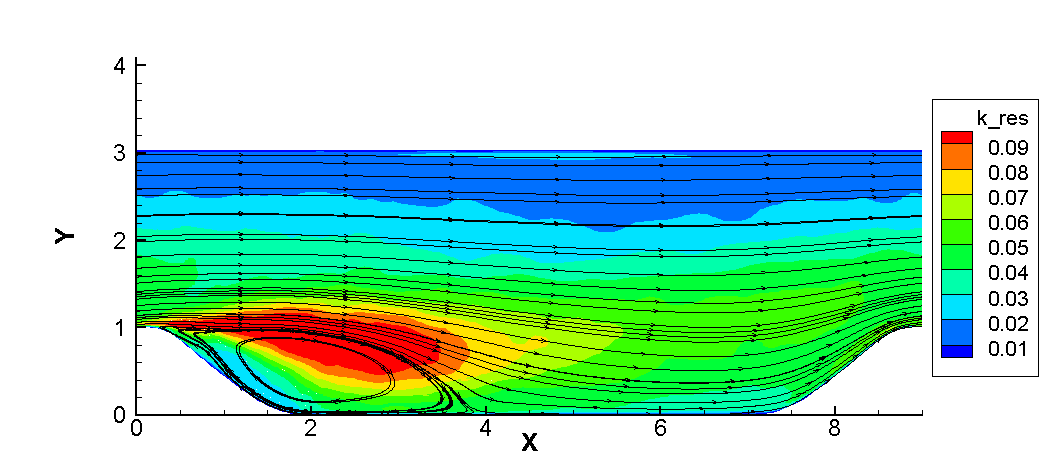}
            \caption{}
        \end{subfigure}
        \begin{subfigure}[b]{0.49\linewidth}
            \includegraphics[clip=true, trim= 0.1cm 0.1cm 0.1cm 0.1cm, width=\textwidth]{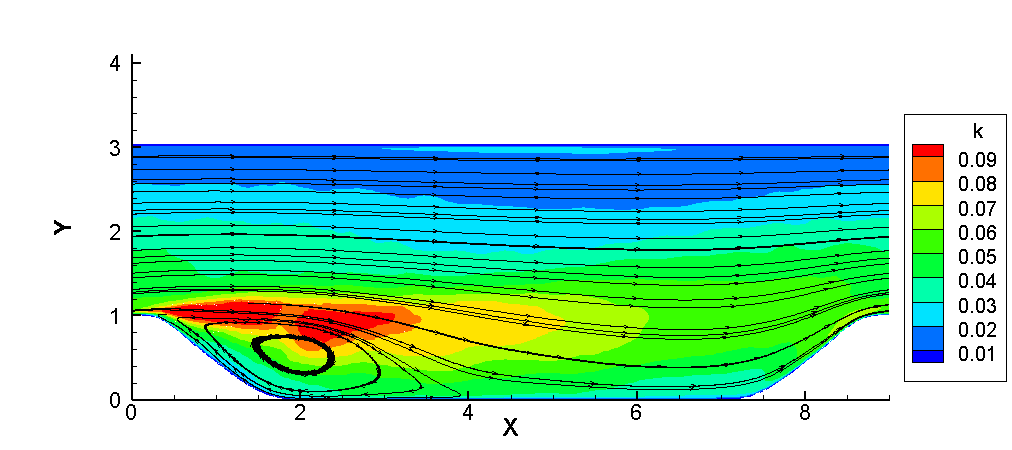}
            \caption{}
        \end{subfigure}
        \caption{Separation bubble and TKE coutour on (a): $ IQ_{\eta} $ adapted grid; (b): $IQ_{k-emp}$ adapted grid; (c): $IQ_{k-ke}$ adapted grid; (d): $IQ_{k-tke}$ adapted grid.}\label{fig:bubble}
    \end{figure}
    
    \begin{figure}[hbt!]
        \centering
        \begin{subfigure}[b]{0.49\linewidth}            
            \includegraphics[clip=true, trim= 0.1cm 0.1cm 0.1cm 0.1cm, width=\textwidth]{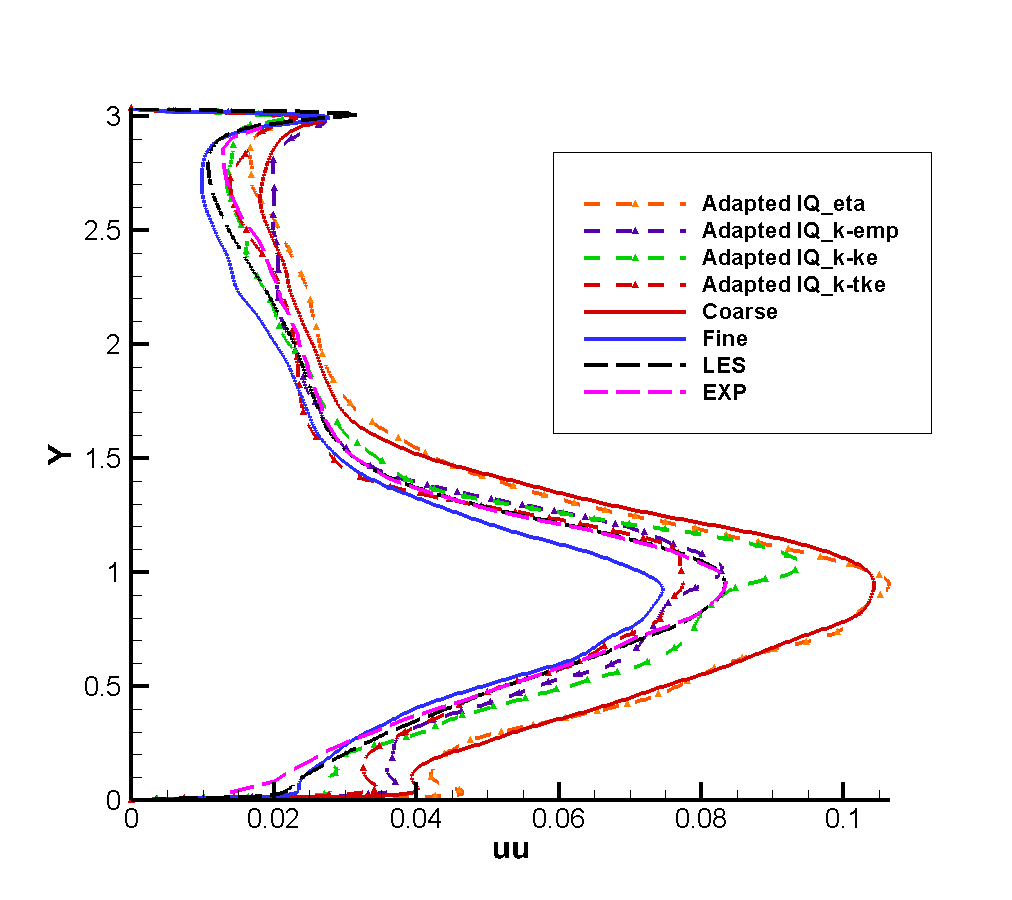}
            \caption{}
        \end{subfigure}
        \begin{subfigure}[b]{0.49\linewidth}
            \includegraphics[clip=true, trim= 0.1cm 0.1cm 0.1cm 0.1cm, width=\textwidth]{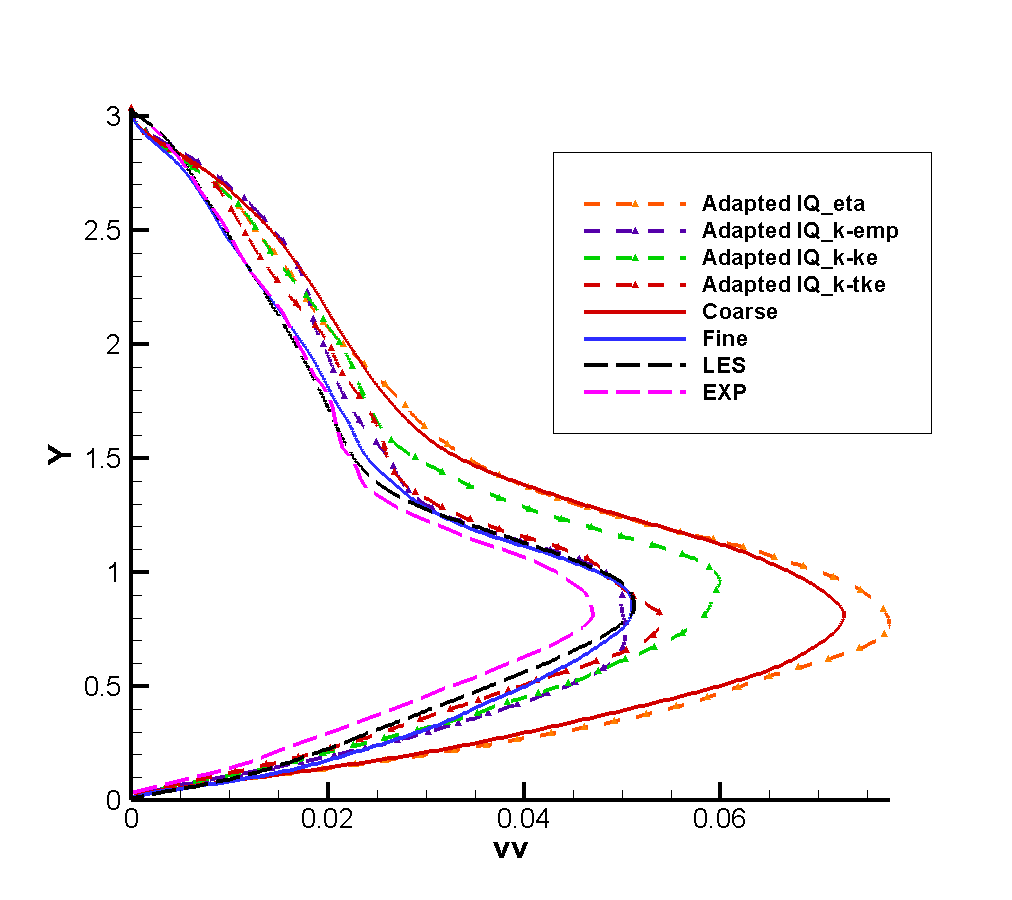}
            \caption{}
        \end{subfigure}
        \begin{subfigure}[b]{0.49\linewidth}
            \includegraphics[clip=true, trim= 0.1cm 0.1cm 0.1cm 0.1cm, width=\textwidth]{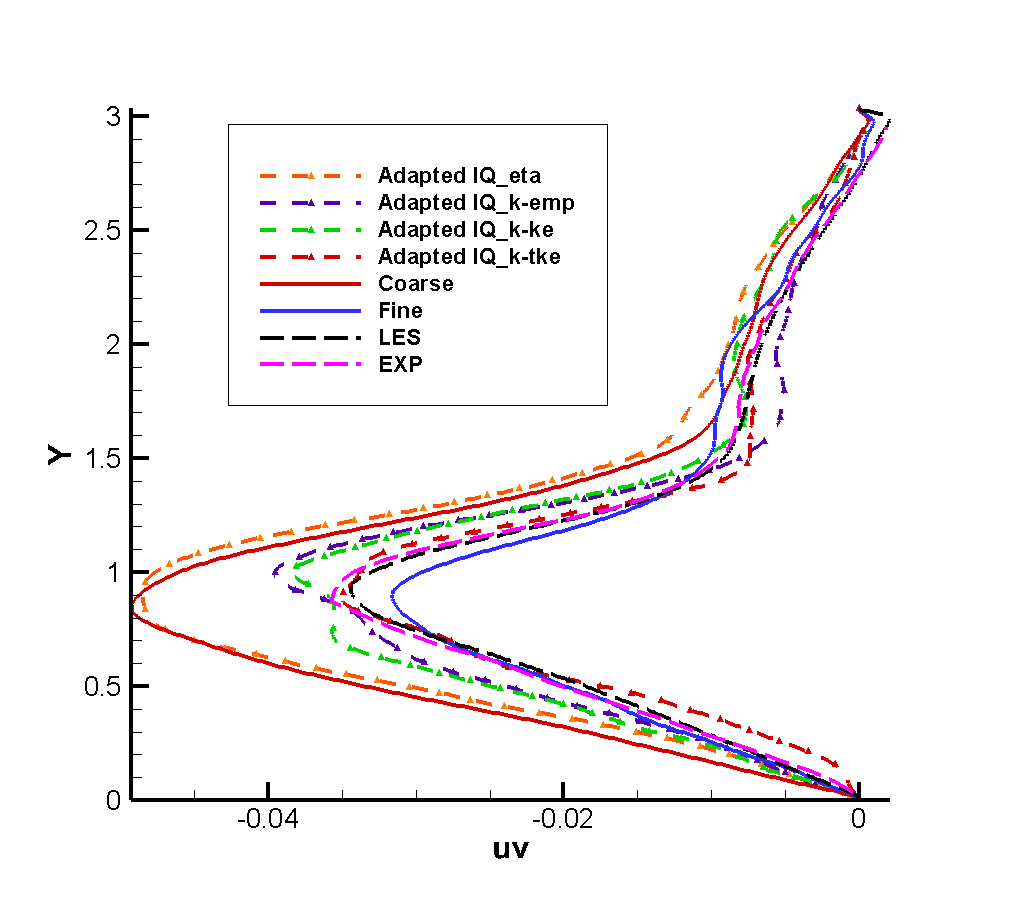}
            \caption{}
        \end{subfigure}
        \begin{subfigure}[b]{0.49\linewidth}
            \includegraphics[clip=true, trim= 0.1cm 0.1cm 0.1cm 0.1cm, width=\textwidth]{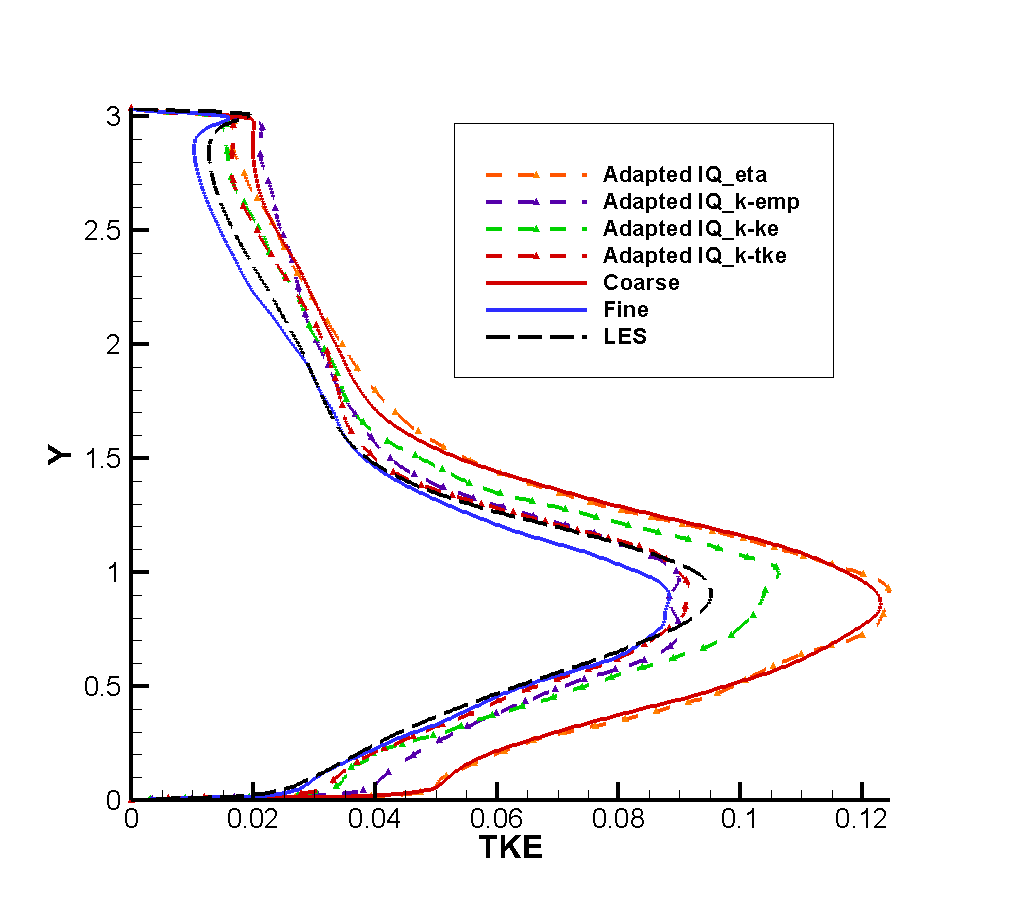}
            \caption{}
        \end{subfigure}
        \caption{Reynolds stress tensor components profiles at $X=2$ (a): $u'u'$ profile; (b): $v'v'$ profile; (c): $u'v'$ profile; (d): TKE profile.}\label{figure:adapt_reynolds_2}
    \end{figure}

    \subsubsection{Simulations on adapted grids}
    
    We have established that the three presented versions of $k_{\text{num}}$ does not lead to a clear difference in the regions targeted for grid adaptation for both $ IQ_{\nu} $ and $ IQ_{\eta} $, while $ IQ_{k-emp} $, $ IQ_{k-ke} $ and $ IQ_{k-tke} $ do target disparate regions. Hence, we decided to carry out a comparative study based on adapted grids generated by $ IQ_{\eta-emp} $, $ IQ_{k-emp} $, $ IQ_{k-ke} $ and $ IQ_{k-tke} $ as shown in Fig.~{\ref{figure:mesh_adapt}}, based on the largest $5\%$ of the error levels of the original cells. The simulation was stabilized on the adapted grids for 90 flow through periods before making an average over the final 30 flow through periods to ensure a meaningful comparison.
    
    Fig.~{\ref{fig:bubble}}~(a), (b), (c) and (d) show the TKE field and the separation bubble captured on the $ IQ_{\eta-emp} $, $ IQ_{k-emp} $, $ IQ_{k-ke} $ and $ IQ_{k-tke} $ adapted grids. It can be observed that the level of captured TKE on $IQ_{k-emp}$ and $IQ_{k-tke}$ adapted grids is significantly reduced compared to the value on the coarse grid shown in Fig.~{\ref{fig:bubble_org}}~(a), while the $IQ_{\eta-emp}$ and $IQ_{k-ke}$ adapted grid still captures a high level of TKE. {\color{black}As shown in Table~\ref{tab:adapt_separation}, $IQ_{\eta-emp}$, $IQ_{k-emp}$ and $IQ_{k-tke}$ adapted grids improve the incorrect location of the separation point, while the $IQ_{k-ke}$ adapted grid still shows an incorrect separation location and provide for a larger separation bubble, whose reattachment point reaches $4.99$. The $IQ_{\eta-emp}$ and $IQ_{k-emp}$ adapted grids underestimate the length of the separation bubble by predicting an early reattachment similar to that on the original coarse grid.} The $IQ_{k-tke}$ adapted grid is able to predict the most accurate location of reattachment point at $X=4.55$. This confirms again that the size of the separation depends on the turbulent intensity captured on the top of the bubble.
    
    Figure.~{\ref{figure:adapt_reynolds_2}} shows the profiles of the Reynolds stress tensor components and the TKE compared with the reference LES and the experimental data at the slice $ X=2 $. The studied slice spans across the main flow and the bubble regions, and is characterized by a high TKE value in the mixing layer above the bubble. The coarse grid over-estimates all the components and wrongly captures the boundary layer while the fine grid shows great comparison against the reference data. It is observed that the $IQ_{\eta-emp}$ adapted grid does not provide a discernible improvement over the coarse grid; while, the $IQ_{k-emp}$ adapted grid shows an improvement of the captured TKE level in most of the computational domain, except for the lower wall region where an over-estimation of the TKE level is still observed.

    {\color{black}Due to the fact that the empirical formula estimates the $k_{\text{num}}$ value based on $k_{\text{sgs}}$, which tends to zero as it approaches the wall, the value of $k_{\text{num}}$ is low in the near-wall region, leading to a high value of estimated $IQ_{k-emp}$. Therefore, even the separation and reattachment regions are not targeted for refinement, which could explain the over-estimation of the TKE near the reattachment point on the $IQ_{k-emp}$ adapted grid.} $IQ_{k-ke}$ targets the high value of numerical TKE near the wall in the bubble region and leads to a refinement in the vicinity of the reattachment point, providing for an improvement in the captured TKE level near the wall, while showing an over-estimation of TKE in the mixing layer on the top region of the separation bubble.
    {The $IQ_{k-tke}$ adapted grid outperforms the other estimators by targeting all the essential regions including the mixing layer and the boundary layer, and by capturing the correct length of the separation bubble as well as a better level of TKE profile in the whole computational domain.}
    
    \subsection{SD7003 airfoil test case}
    \label{sec:sd7003}
    \begin{figure}[hbt!]
        \begin{center}
            \includegraphics[clip=true, trim= 0.1cm 4.1cm 0.1cm 5.1cm, width=.49\linewidth]{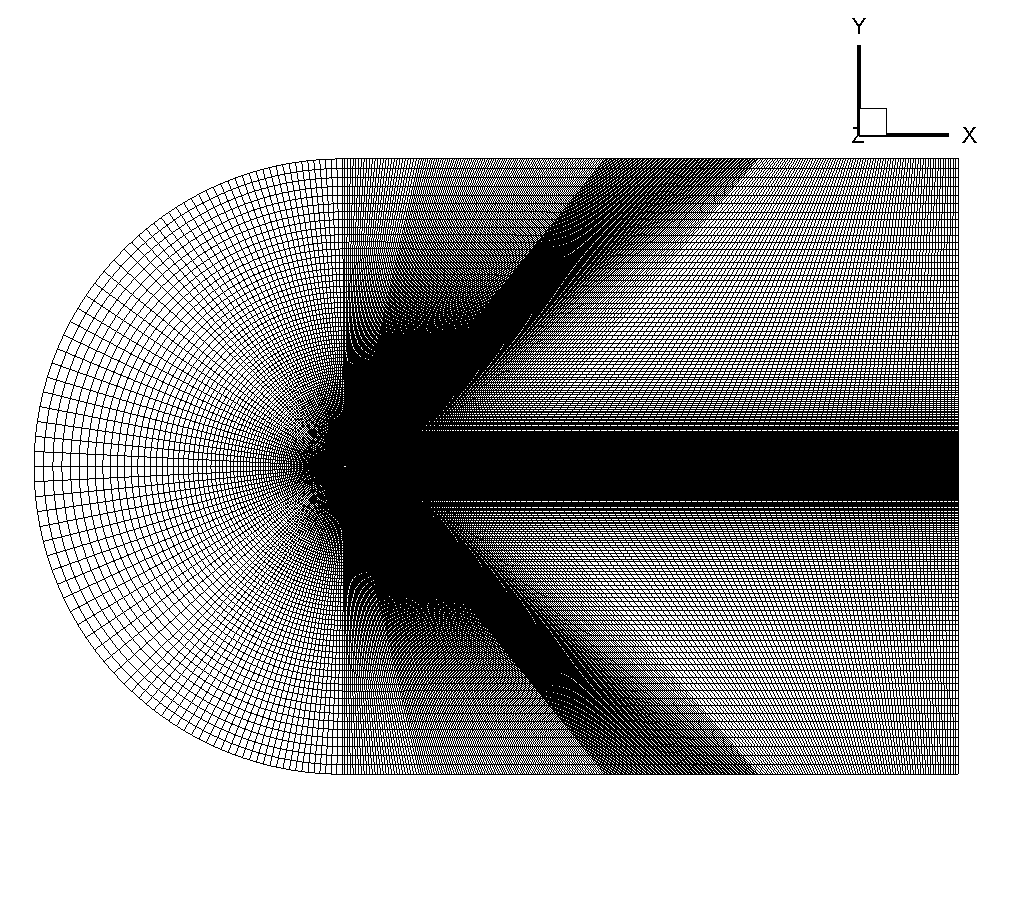}
            \includegraphics[clip=true, trim= 0.1cm 4.1cm 0.1cm 5.1cm, width=.49\linewidth]{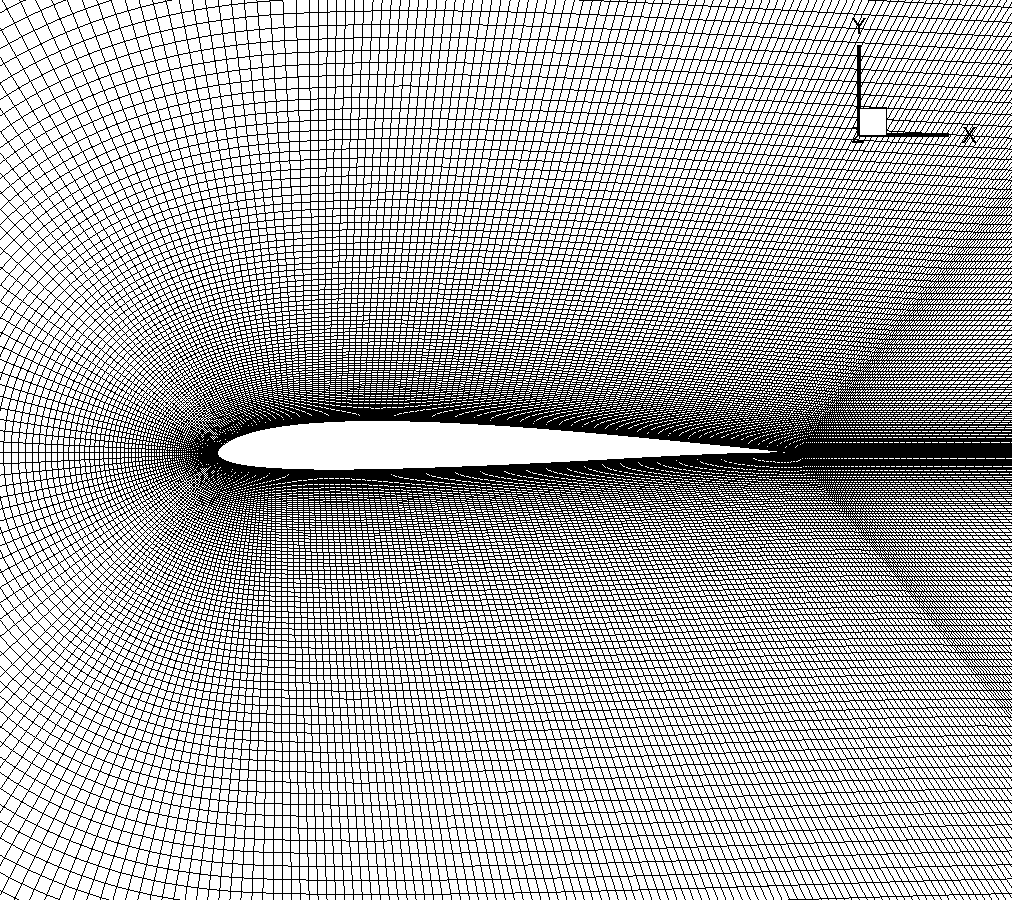}
        \end{center}
        \caption{SD7003 case grid.}
        \label{figure:d7003grid} 
    \end{figure}
    \subsubsection{Introduction}
The SD7003 belongs to the SD family of low-Reynolds number airfoils~\cite{selig1989airfoils} designed to minimize drag losses associated with the laminar separation bubble. The flow shows a laminar separation near the leading edge and transitions to turbulent flow before reattachment. The correct prediction of the separation bubble length relies on the correct prediction of the TKE level around the bubble. Our numerical test case is the SD7003 airfoil case with $Re=60000$ which has been widely studied in references~\cite{uranga2011implicit,garmann2013comparative}. We mainly focus on the case at an angle of attack $\alpha = 8^o$, where the flow shows more turbulent features. The C-grid extends to approximately $25$ chord lengths around the airfoil. The computational domain extends in the spanwise direction by $0.2$c following Uranga et al.~\cite{uranga2011implicit} who showed that a domain span of $0.2$c was sufficient for capturing spanwise structures. An original coarse grid of 2.16M cells (shown in Fig.~\ref{figure:d7003grid}) and the correspondingly refined grid with 7.99M cells are studied for comparison. To ensure the flow reaches statistically steady state, simulations are launched over nondimensionalized time $T=\frac{tU}{c}=120$, and the flow field is averaged in time over the final $T=30$ and in the spanwise direction to evaluate the error estimators and perform requisite flow analysis. The proposed error estimators are then applied to the coarse grid in order to generate higher quality grids. Only a small percentage ($5\%$) of the original grid cells are targeted and refined. The goal is to target critical regions where the grid refinement will lead to an improvement of the simulation quality, including the prediction of lift and drag coefficients and the structure of the separation bubble, with less additional computational cost compared to the fine grid.

\begin{table}
    \caption{Separation and reattachment points.}
    \centering
    \begin{tabular}{rrrrrcc}
        \hline
        Simulation & DOF & Sep.  & Reat.  & Length & $C_l$ & $C_d$\\
        \hline
        Garmann et al.\cite{garmann2013comparative}   & 54M & 0.031 & 0.303 & 0.272 & 0.917 & 0.0447\\
        Vermeire et al.\cite{vermeire2016implicit} & 2.3M & 0.031 & 0.345 & 0.314 & 0.946 & 0.0529\\
        Selig et al.~\cite{selig1989airfoils,selig1995summary} & Exp & - & - & - &  0.88 $\sim 0.95$ & 0.03 $\sim$ 0.04 \\
        Coarse   & 2.16M & 0.1368 & 0.4451 & 0.3083 & 0.9618 & 0.04261\\
        Fine   & 7.99M & 0.0201 & 0.3349 & 0.3148 & 0.9433 & 0.05273\\
        $IQ_{k-emp}$ adapted   & 2.76M & 0.1062 & 0.3473 & 0.2411 & 0.9246 & 0.05041\\
        $IQ_{k-ke}$ adapted   & 2.76M & 0.1214 & 0.4225 & 0.3011 & 0.9395 & 0.04681\\
        $IQ_{k-tke}$ adapted   & 2.76M & \textbf{0.0430} & \textbf{0.3026} & \textbf{0.2596} & 0.9193 & 0.04731\\
        \hline
    \end{tabular}%
    \label{tab:sd7003_summary}
\end{table}

\begin{table}
    \caption{\color{black} Grids used for the SD7003 airfoil test cases.}
    \centering
    \begin{tabular}{rrrrr}
        \hline
        Simulation & DOF & $\Delta s^{+}$ & $\Delta n^{+}$ & $\Delta z^{+}$\\
        \hline
        Garmann et al.\cite{garmann2013comparative}   & 54M & $\approx 10 $ & $\approx 0.2$~\footnotemark   & $\approx 10$\\
        Uranga et al.\cite{uranga2011implicit} & 1.8M & - & $\approx 2.0$ & $\approx 53$ \\
        Vermeire et al.\cite{vermeire2016implicit} & 2.3M & $\approx 22$ & $\approx 0.6$ & $\approx 20$\\
        Coarse   & 2.16M & $\approx 24$ & $\approx 3$ & $\approx 40$\\
        Fine   & 7.99M & $\approx 17$ & $\approx 1$ & $\approx 20$\\
        adapted   & 2.76M & $\approx 12$ & $\approx 1.5$ & $\approx 20$\\
        \hline
    \end{tabular}%
    \label{tab:sd7003_grid}
\end{table}
\footnotetext{Value from \cite{garmann2013comparative} is based on maximum $\Delta n^{+}$.}
\subsubsection{Flow solution}

\begin{figure}[hbt!]
        \centering
        \begin{subfigure}[b]{0.4\linewidth}            
            \includegraphics[clip=true, trim= 0.1cm 2.4cm 0.1cm 6.1cm, width=\textwidth]{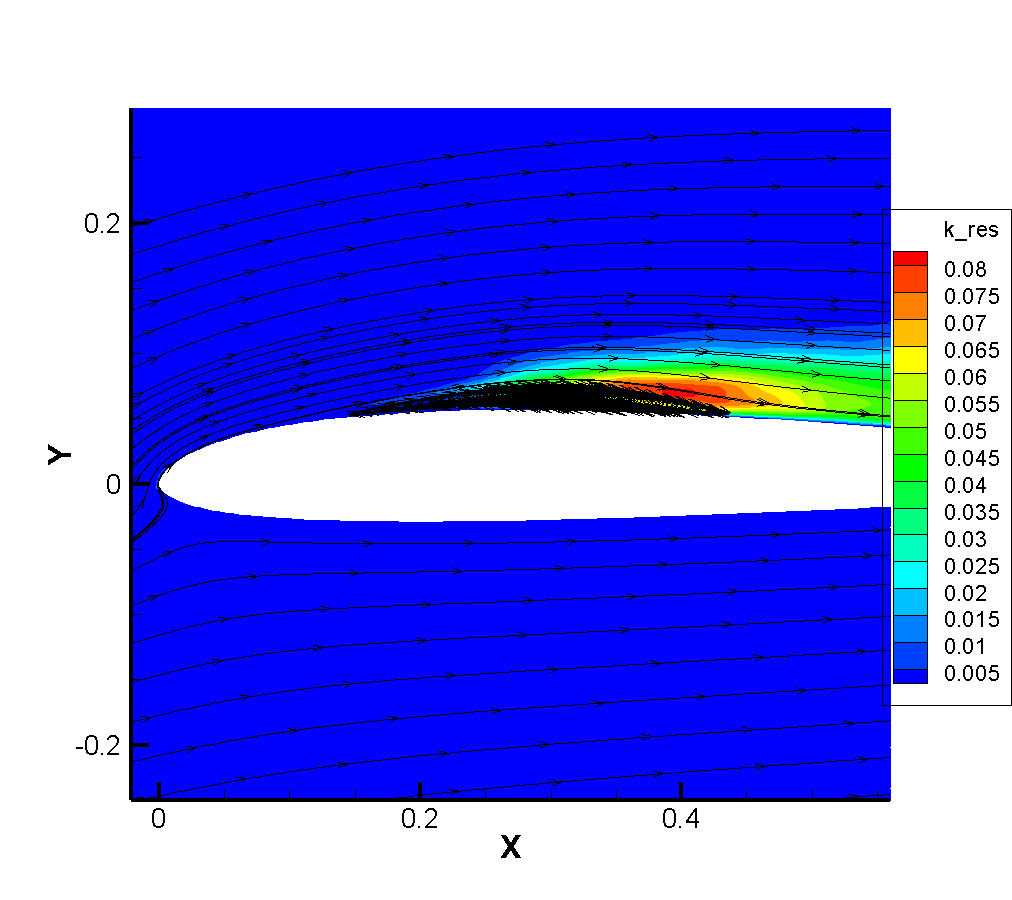}
            \caption{}
        \end{subfigure}
        \begin{subfigure}[b]{0.4\linewidth}            
            \includegraphics[clip=true, trim= 0.1cm 3.1cm 0.1cm 7.7cm, width=\textwidth]{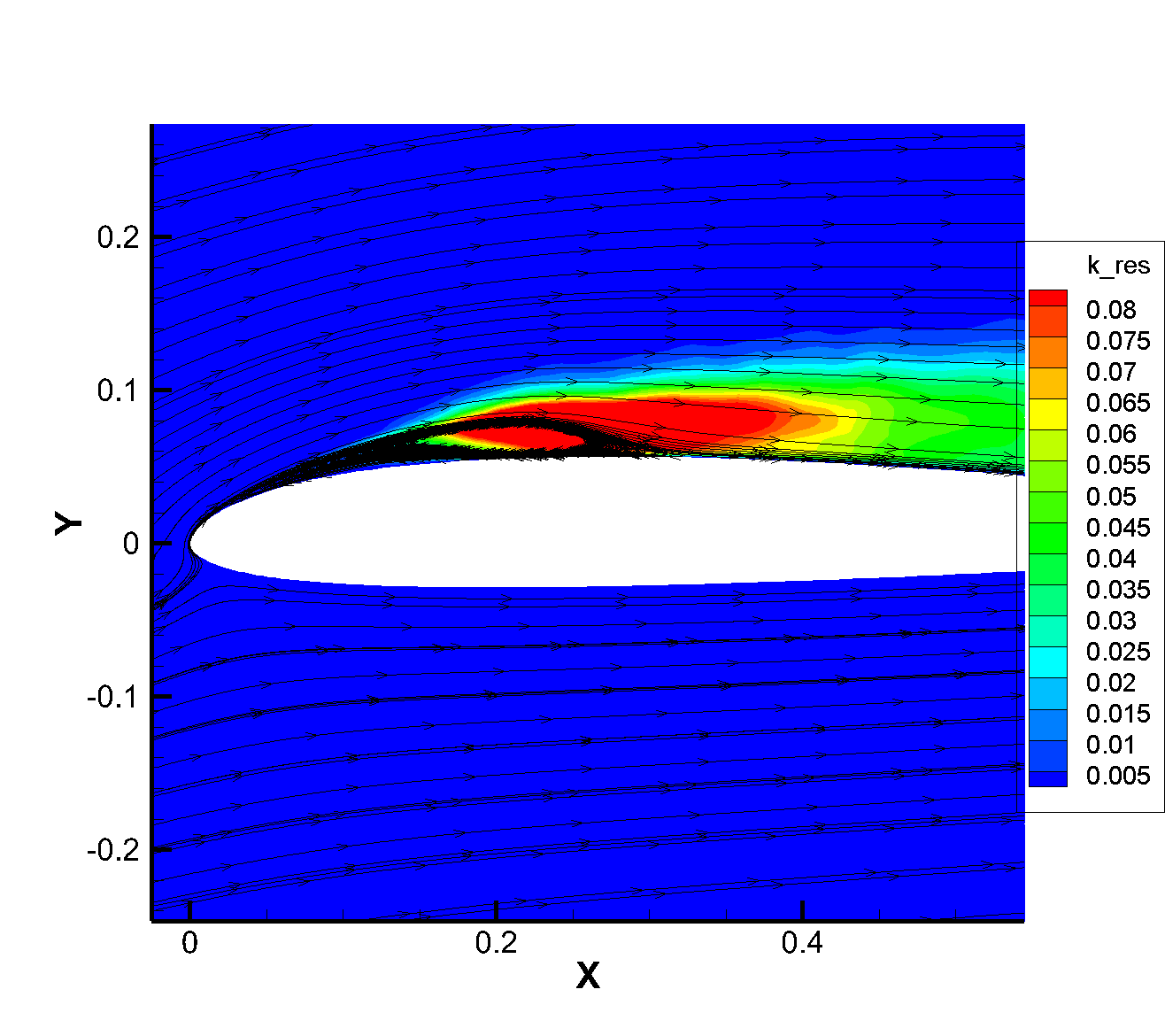}
            \caption{}
        \end{subfigure}
        \begin{subfigure}[b]{0.4\linewidth}            
            \includegraphics[clip=true, trim= 0.1cm 2.4cm 0.1cm 6.1cm, width=\textwidth]{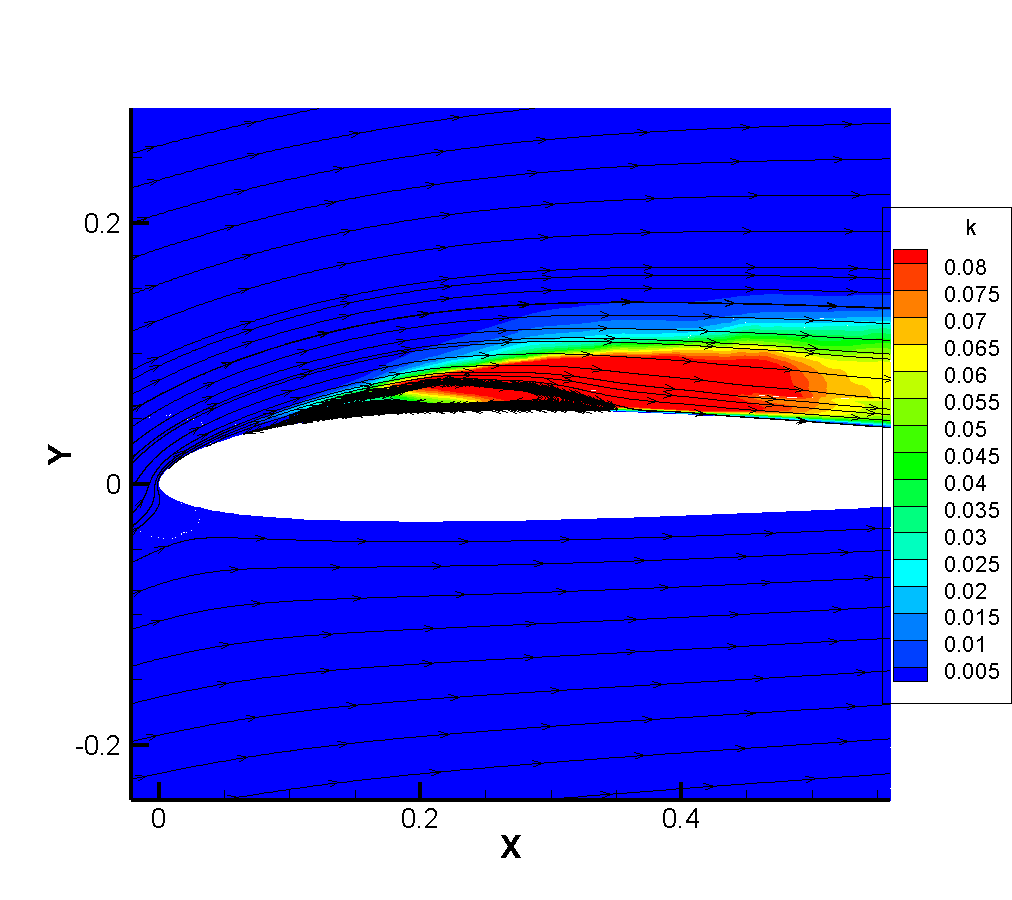}
            \caption{}
        \end{subfigure}
        \begin{subfigure}[b]{0.4\linewidth}
            \includegraphics[clip=true, trim= 0.1cm 2.4cm 0.1cm 6.1cm, width=\textwidth]{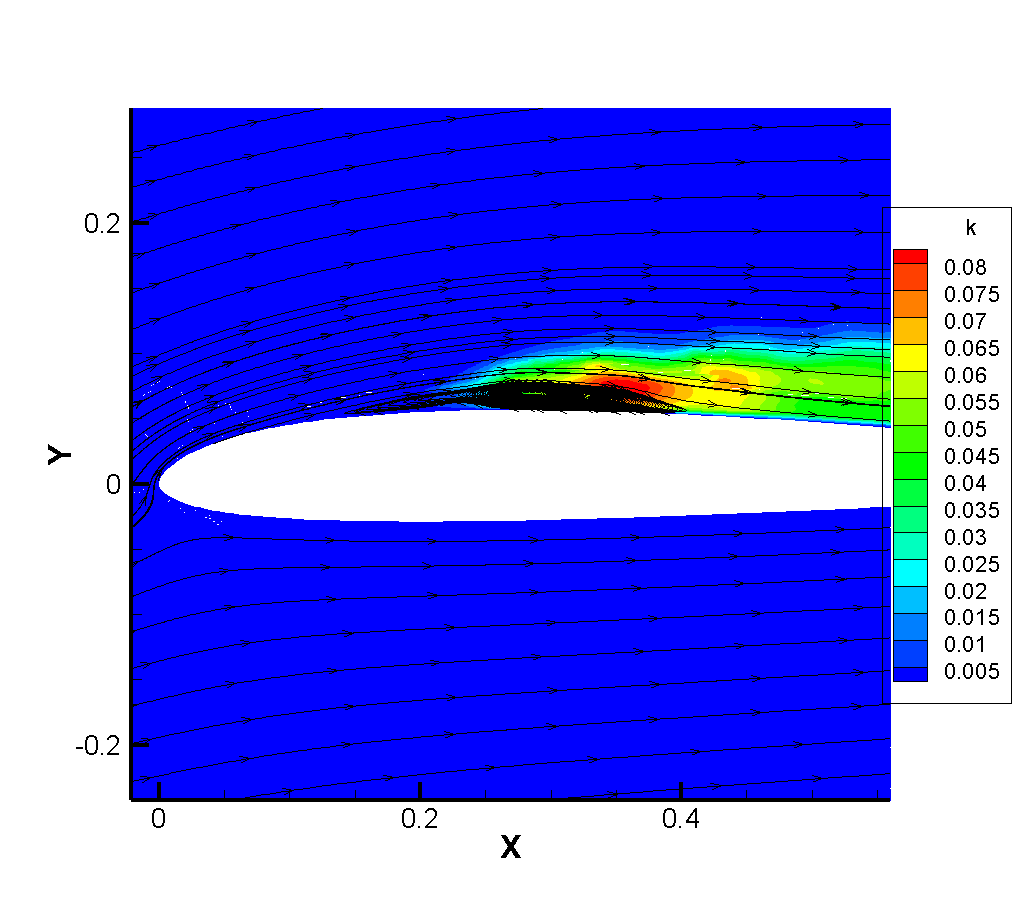}
            \caption{}
        \end{subfigure}
        \begin{subfigure}[b]{0.4\linewidth}
            \includegraphics[clip=true, trim= 0.1cm 2.4cm 0.1cm 6.1cm, width=\textwidth]{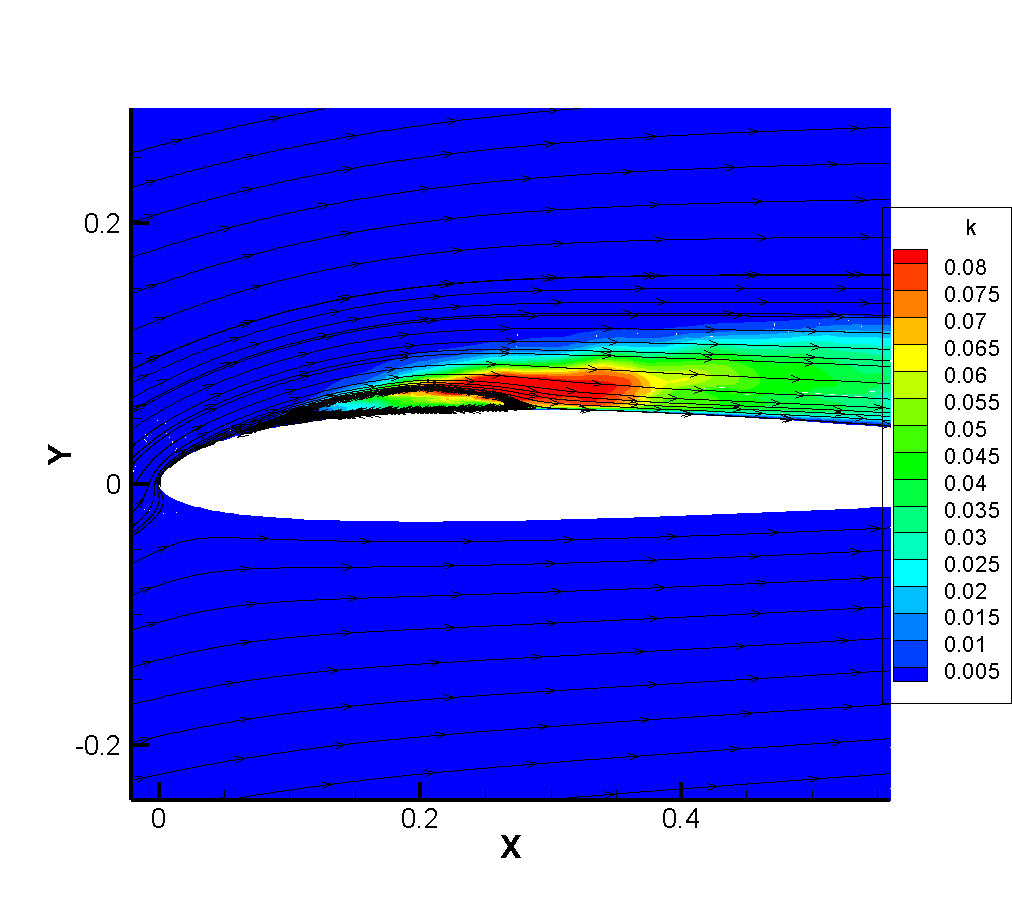}
            \caption{}
        \end{subfigure}
        \caption{Separation bubble and TKE contour on (a): coarse grid; (b): fine grid; (c): $IQ_{k-emp}$ adapted grid; (d): $IQ_{k-ke}$ adapted grid; (e): $IQ_{k-tke}$ adapted grid.}\label{fig:sd7003bubble}
    \end{figure}

{\color{black} Results on two levels of grids are available in Table~\ref{tab:sd7003_summary}, compared with the reference LES~\cite{garmann2013comparative,vermeire2016implicit} and experimental~\cite{selig1989airfoils,selig1995summary} results. There are two key observations. First, the noted range for the aerodynamic coefficients in the experimental data is due to differences in the spanwise boundary setup~\cite{selig1995summary}. Radespiel et al.~\cite{radespiel2007numerical} reported that the non two-dimensional wake flow pattern results in a strong influence of the spanwise measurement position on the drag coefficient. Hence, due to the uncertainty in the measured coefficients, two reliable LES results~\cite{garmann2013comparative,vermeire2016implicit} are used as benchmark values. Second, there is relatively good agreement between the two LES benchmark values~\cite{garmann2013comparative,vermeire2016implicit} in terms of the separation point but less so on the reattachment point, and lift and drag coefficients. We believe differences in the computational grid setup may have contributed to the observed differences. Garmann et al.~\cite{garmann2013comparative} reported that the grid resolution in the vicinity of the separation was a strong influence on the proper capture of the flow field and provided values for the minimum non-dimensional normal distance at the streamwise location $x/c = 0.1$. We have tabulated the data in Table~\ref{tab:sd7003_grid} and compared the values against other references as well as results achieved in this work. The normal distance of the first solution point off the surface of the airfoil is calculated as $\Delta n^+=u_{\tau}y/\nu$ where $u_{\tau}=\sqrt{C_f/2}U_{\infty}$ and $C_f \approx 0.01$ at the streamwise location x/c = 0.1~\cite{uranga2011implicit}. The measure of the cell resolution for the streamwise ($\Delta s^+$) and spanwise ($\Delta z^+$) directions are calculated in the same way and given in Table~\ref{tab:sd7003_grid}. Garmann et al.~\cite{garmann2013comparative} studied three levels of grids ranging from 6.7M to 12M and adopted a large computational domain extending to 100 chord lengths. As shown in Table~\ref{tab:sd7003_grid} the boundary layer is extremely refined for the grid used in~\cite{garmann2013comparative} with a $\Delta n^+ \approx 0.2$. However, Uranga et al.~\cite{uranga2011implicit} showed that a much coarser grid was able to provide sufficient resolution for the test case through a grid study using implicit LES, where they employed a grid with a computational domain that extended to 6 chord lengths with 1.8M degrees of freedom (DOF) and a $\Delta n^+ \approx 2$ at $x/c=0.1$. {\color{black}A validation of our in-house code for  $\alpha = 4^o$ resulted in $C_l$ and $C_d$ coefficients within 3\% of~\cite{uranga2011implicit}. Following~\cite{uranga2011implicit}, Vermeire et al.~\cite{vermeire2016implicit} adopted the same computational domain with a finer grid resulting to a $\Delta n^+ \approx 0.6$ at $x/c=0.1$.
Based on the computational grids adopted in~\cite{garmann2013comparative,vermeire2016implicit,uranga2011implicit}, we employed a computational domain that stretched to 25 chord lengths, 2.16M DOF, and a $\Delta n^+ \approx 3$ to serve as our coarse grid while the fine grid would be approximately 7.99M DOF and a $\Delta n^+ \approx 1$. As expected, our coarse grid estimates a delayed separation and reattachment of the flow, and also overestimates the lift coefficient $C_l$ and underestimates the drag coefficient $C_d$ compared to the fine grid and the reference LES results. The fine grid shows an overestimation of both $C_l$ and $C_d$ compared to~\cite{garmann2013comparative}, but show excellent agreement against~\cite{vermeire2016implicit}. The presented results highlight a strong correlation between the bubble length and the aerodynamic performance coefficients. Both the reported fine grid and~\cite{vermeire2016implicit} have equivalent bubble lengths and hence comparable $C_l$ and $C_d$ values. Fig.~{\ref{fig:sd7003bubble}}~(a) and (b) present the size of the separation bubble along with the contour of the captured TKE on the coarse and fine grids.  Fig.~{\ref{fig:sd7003bubble}}~(c) to (e) present flow solutions on adapted grids. These results  will be discussed in a subsequent  Sec.~\ref{sec:sd7003_sim_adapted_grids}. The delay and overestimation of the length of the separation bubble on the coarse grid are due to the incorrect delayed prediction of the transition point and an underestimation of the turbulence level above the bubble region.} As a result, the boundary layer is incorrectly captured on the coarse grid which leads to the observed errors in $C_l$ and $C_d$. The corresponding fine grid is able to better capture the separation bubble with a correct length as shown in Fig.~{\ref{fig:sd7003bubble}}~(b) and Table~\ref{tab:sd7003_summary}.

\subsubsection{Evaluation of numerical TKE}

    \begin{figure}[hbt!]
        \begin{center}
            \includegraphics[clip=true, trim= 0.1cm 0.1cm 0.1cm 0.1cm, width=.33\linewidth]{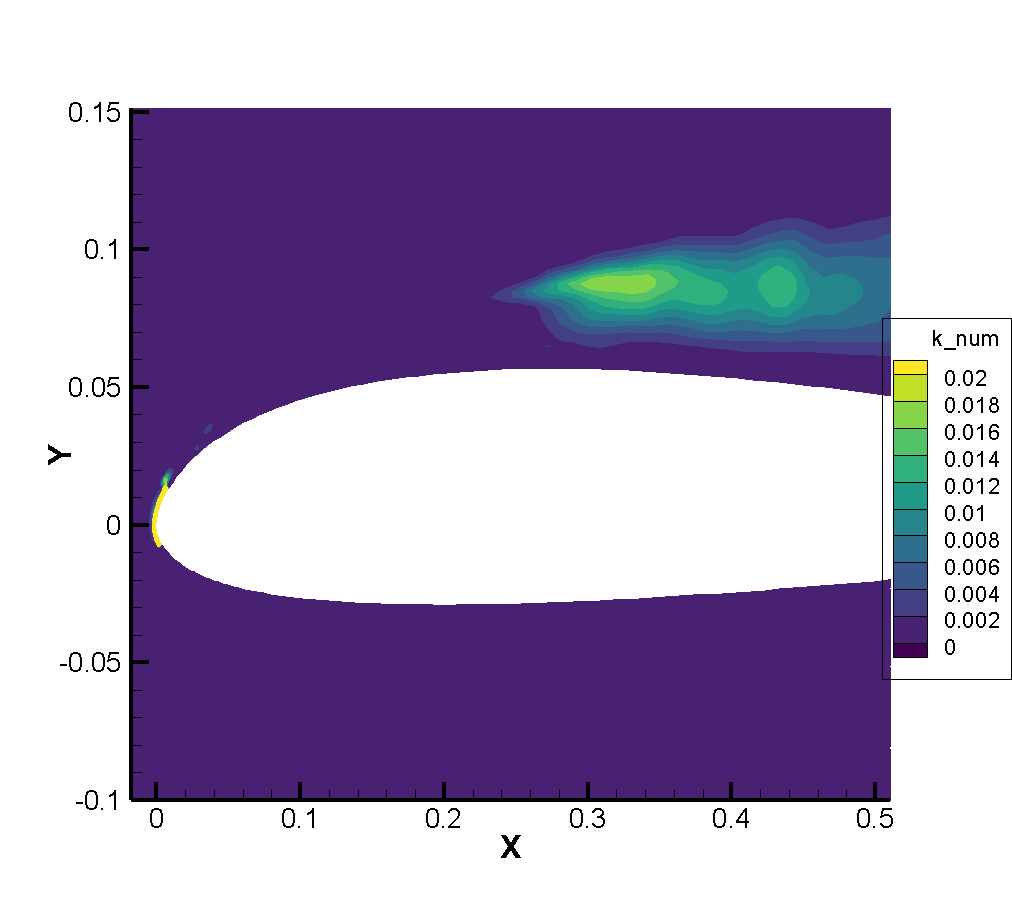}
            \includegraphics[clip=true, trim= 0.1cm 0.1cm 0.1cm 0.1cm, width=.33\linewidth]{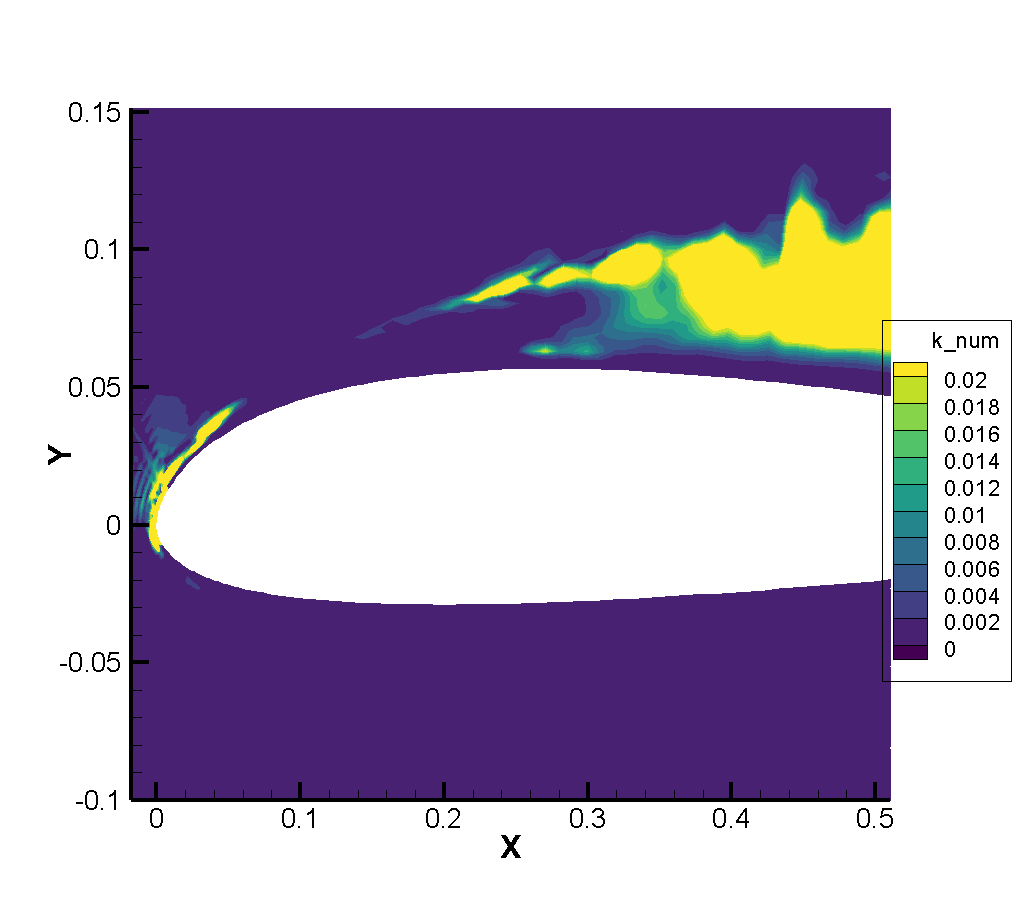}
            \includegraphics[clip=true, trim= 0.1cm 0.1cm 0.1cm 0.1cm, width=.33\linewidth]{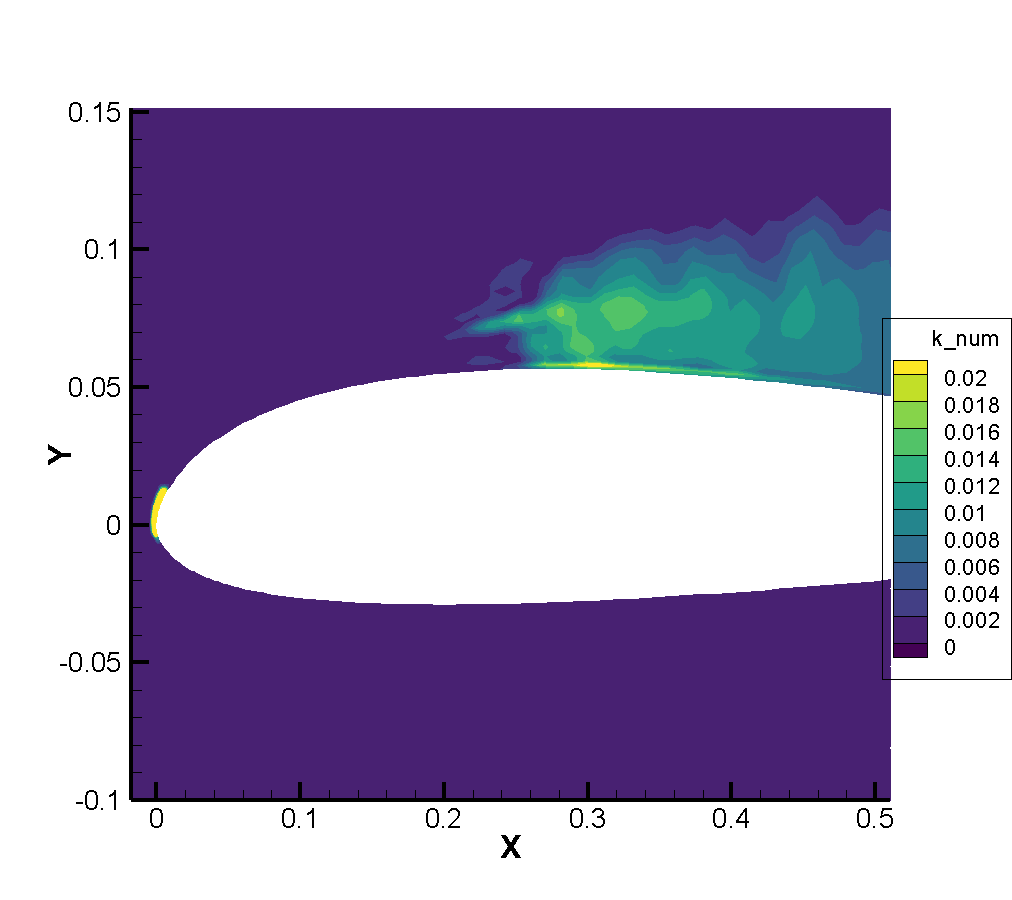}
        \end{center}
        \caption{Numerical TKE, $k_{\text{num}}$, based on empirical formula (left), KE numerical dissipation (middle) and TKE numerical dissipation (right).}
        \label{figure:sd7003_k_num}
    \end{figure}
In order to provide for the Index Quality on the coarse grid, the numerical TKE is estimated based on the three approaches presented in Section~\ref{sec:num_tke}. It is observed in Fig.~\ref{figure:sd7003_k_num} that all three approaches only show positive values of the numerical TKE and highlight similar regions that includes the leading edge area and the shear layer on top of the separation bubble starting at $X\approx 0.2$ where flow transitions to be turbulent. The main difference relies in the captured $k_{\text{num}}$ value level and the boundary layer. KE-based approach shows comparatively higher value of numerical TKE compared to other approaches. The TKE-based approach is able to highlight the near-wall region on the upper surface of the airfoil after the transition, especially targeting the boundary layer on the upper surface for a high value of numerical TKE even though the modeled TKE is low in the boundary layer; {\color{black} while both the empirical and KE-based approaches fail to flag the boundary layer.} Recalling that the {\color{black}KE-based approach estimates the numerical eddy viscosity $\nu_{\text{num}}$ based on the ratio between $\bar{\epsilon}_{n}$ and $\overline{\tau_{ij}\frac{\partial u_i}{\partial x_j}}$ as shown in Eq.~\ref{eq:ratio}}, the high value of $\overline{\tau_{ij}\frac{\partial u_i}{\partial x_j}}$ in the boundary layer due to a large velocity gradient actually leads to a low estimated value of numerical eddy viscosity $\nu_{\text{num}}$, {\color{black} resulting in a low value of estimated $k_{\text{num}}$ following Eq.~\ref{nu_num}.}}

\subsubsection{Error estimation}

    \begin{figure}[hbt!]
        \centering
        \begin{subfigure}[b]{0.33\textwidth}            
            \includegraphics[clip=true, trim= 0.1cm 0.1cm 0.1cm 0.1cm, width=\linewidth]{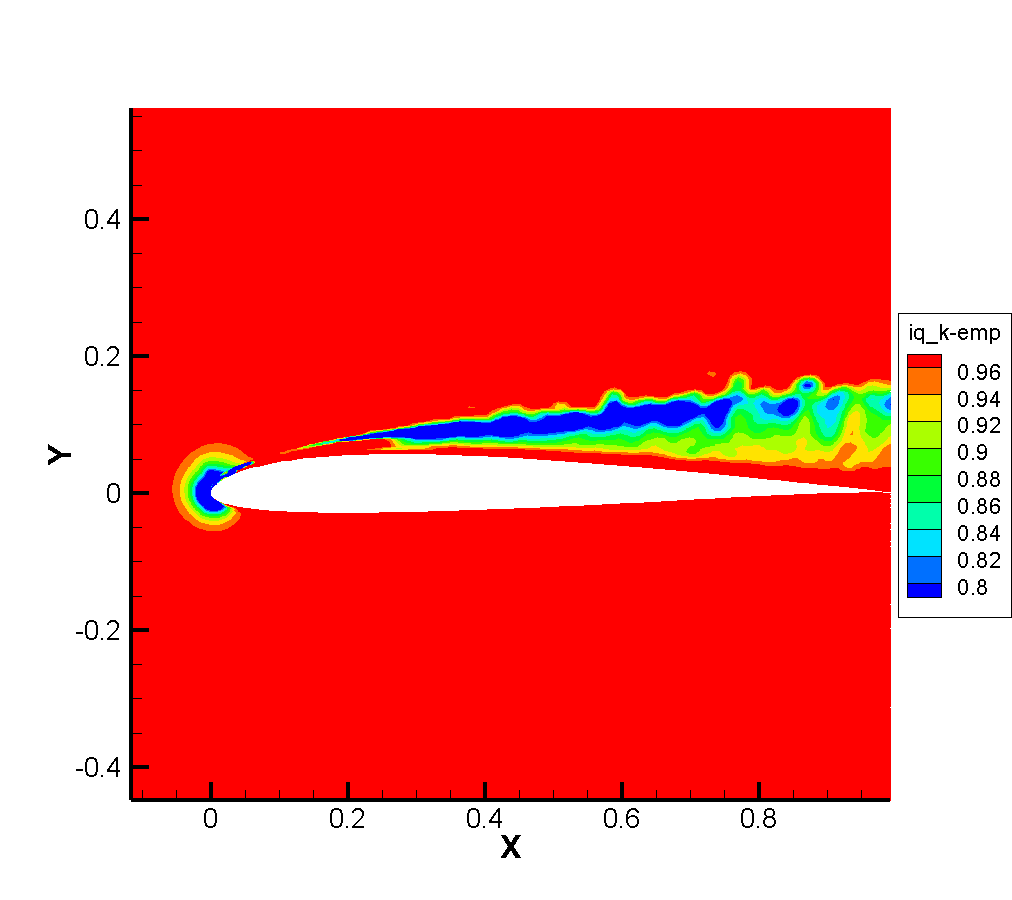}
            \caption{}
        \end{subfigure}
        \begin{subfigure}[b]{0.33\textwidth}            
            \includegraphics[clip=true, trim= 0.1cm 0.1cm 0.1cm 0.1cm, width=\linewidth]{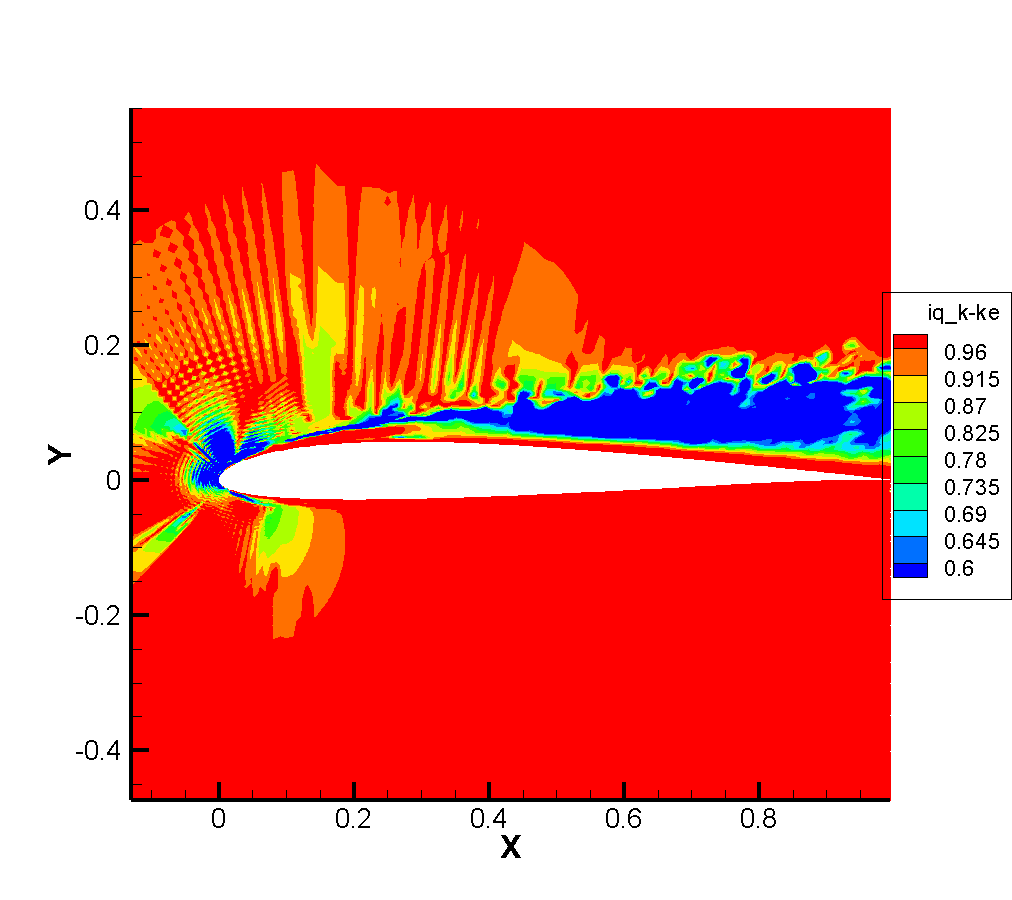}
            \caption{}
        \end{subfigure}
        \begin{subfigure}[b]{0.33\textwidth}            
            \includegraphics[clip=true, trim= 0.1cm 0.1cm 0.1cm 0.1cm, width=\linewidth]{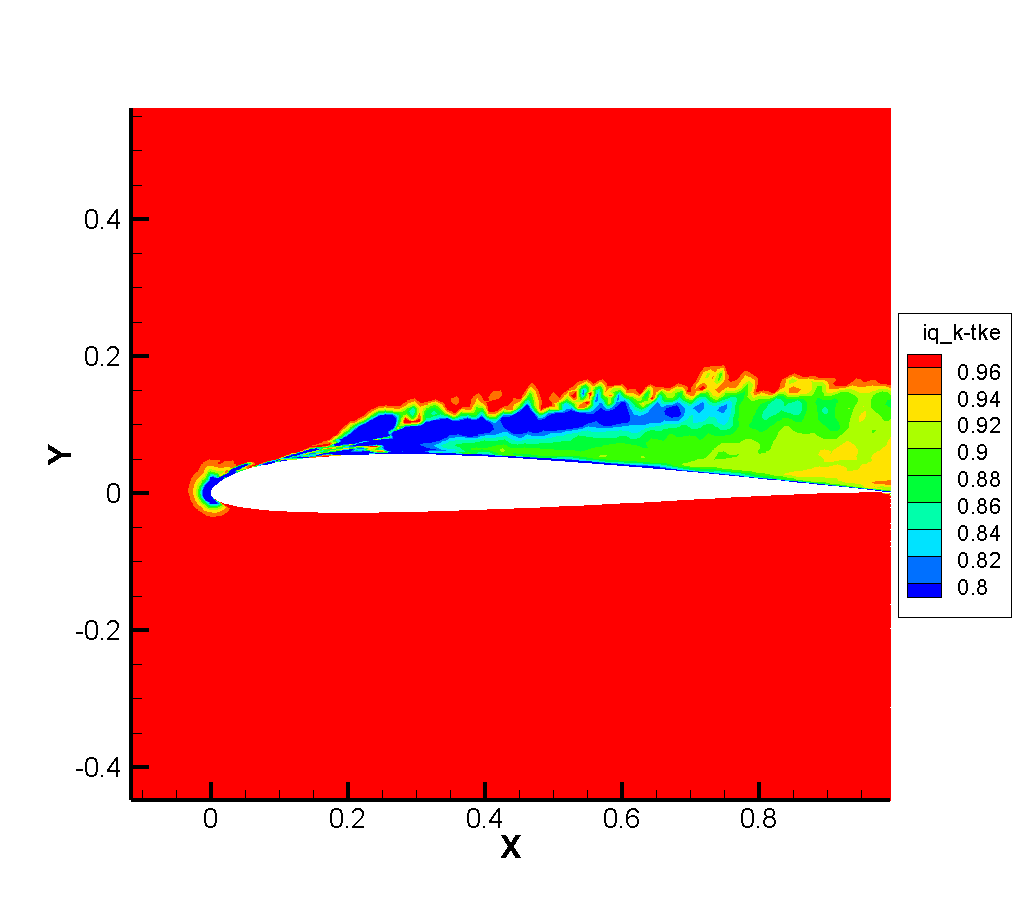}
            \caption{}
        \end{subfigure}
        \caption{Error estimation on coarse grid based on (a): $IQ_{k-emp}$; (b): $IQ_{k-ke}$; (c): $IQ_{k-tke}$.}\label{fig:sd7003error}
    \end{figure}
    
\begin{figure}[hbt!]
        \centering
        \begin{subfigure}[b]{0.33\textwidth}            
            \includegraphics[clip=true, trim= 0.1cm 0.1cm 0.1cm 0.1cm, width=\linewidth]{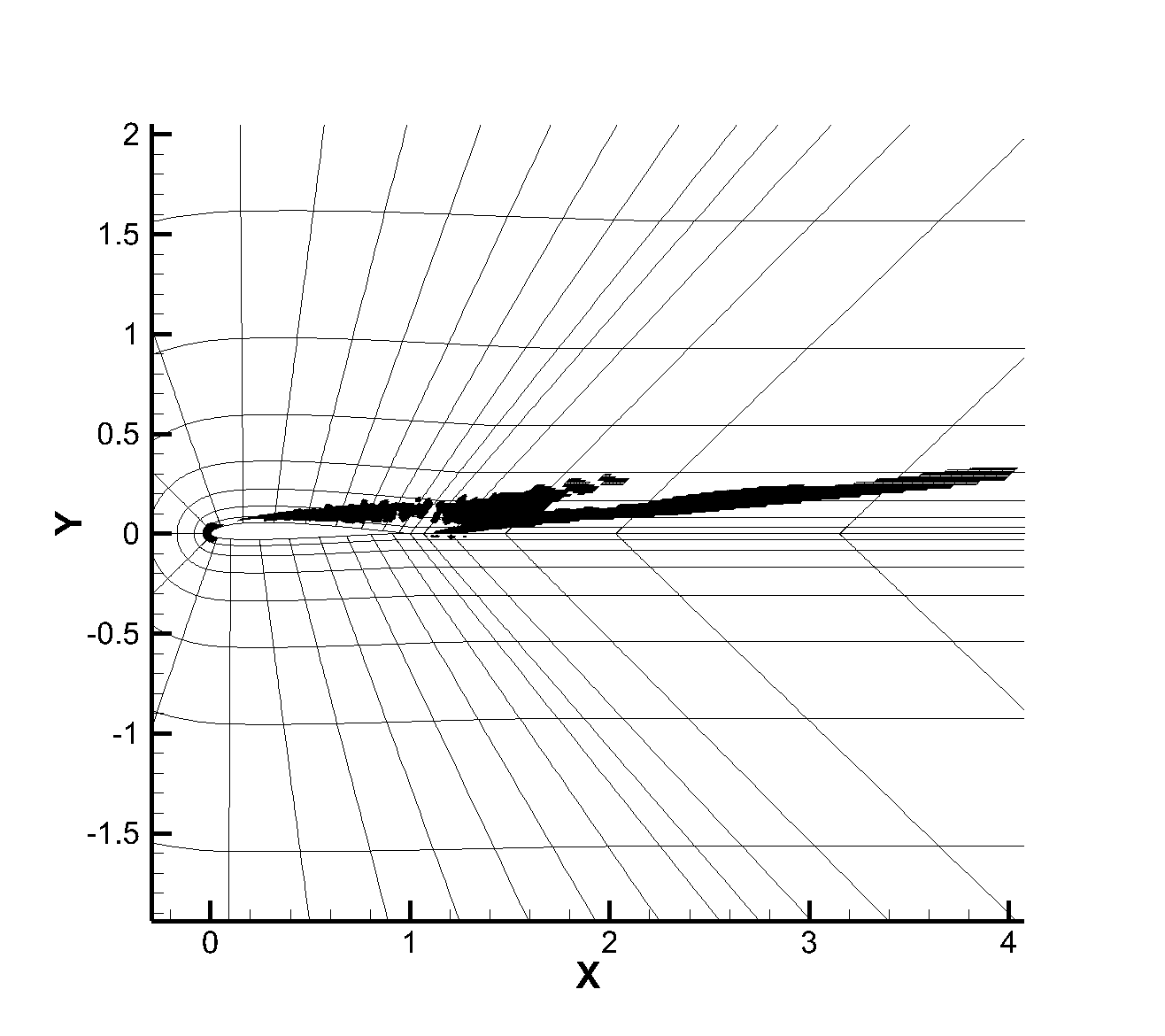}
            \caption{}
        \end{subfigure}
        \begin{subfigure}[b]{0.33\textwidth}            
            \includegraphics[clip=true, trim= 0.1cm 0.1cm 0.1cm 0.1cm, width=\linewidth]{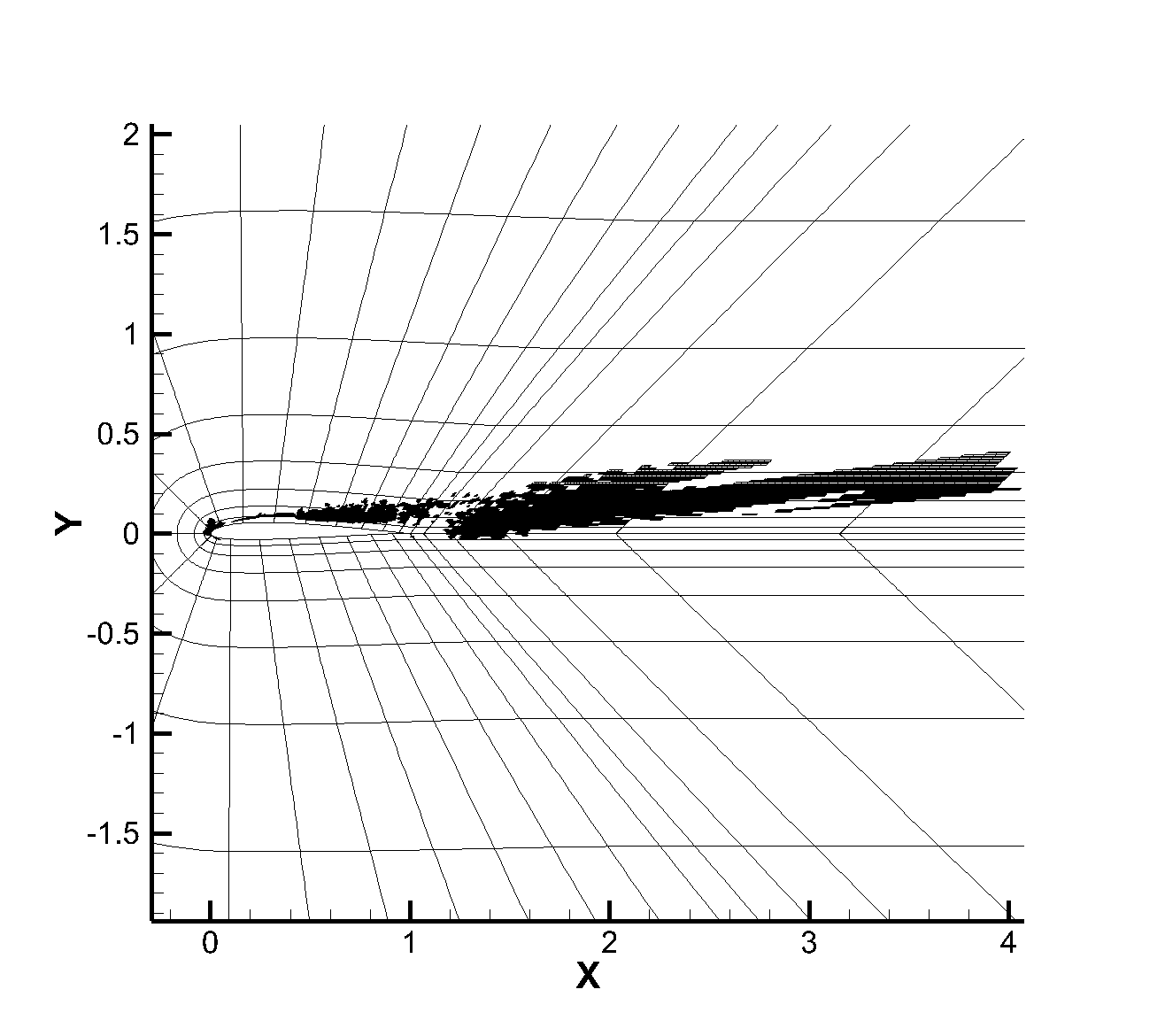}
            \caption{}
        \end{subfigure}
        \begin{subfigure}[b]{0.33\textwidth}            
            \includegraphics[clip=true, trim= 0.1cm 0.1cm 0.1cm 0.1cm, width=\linewidth]{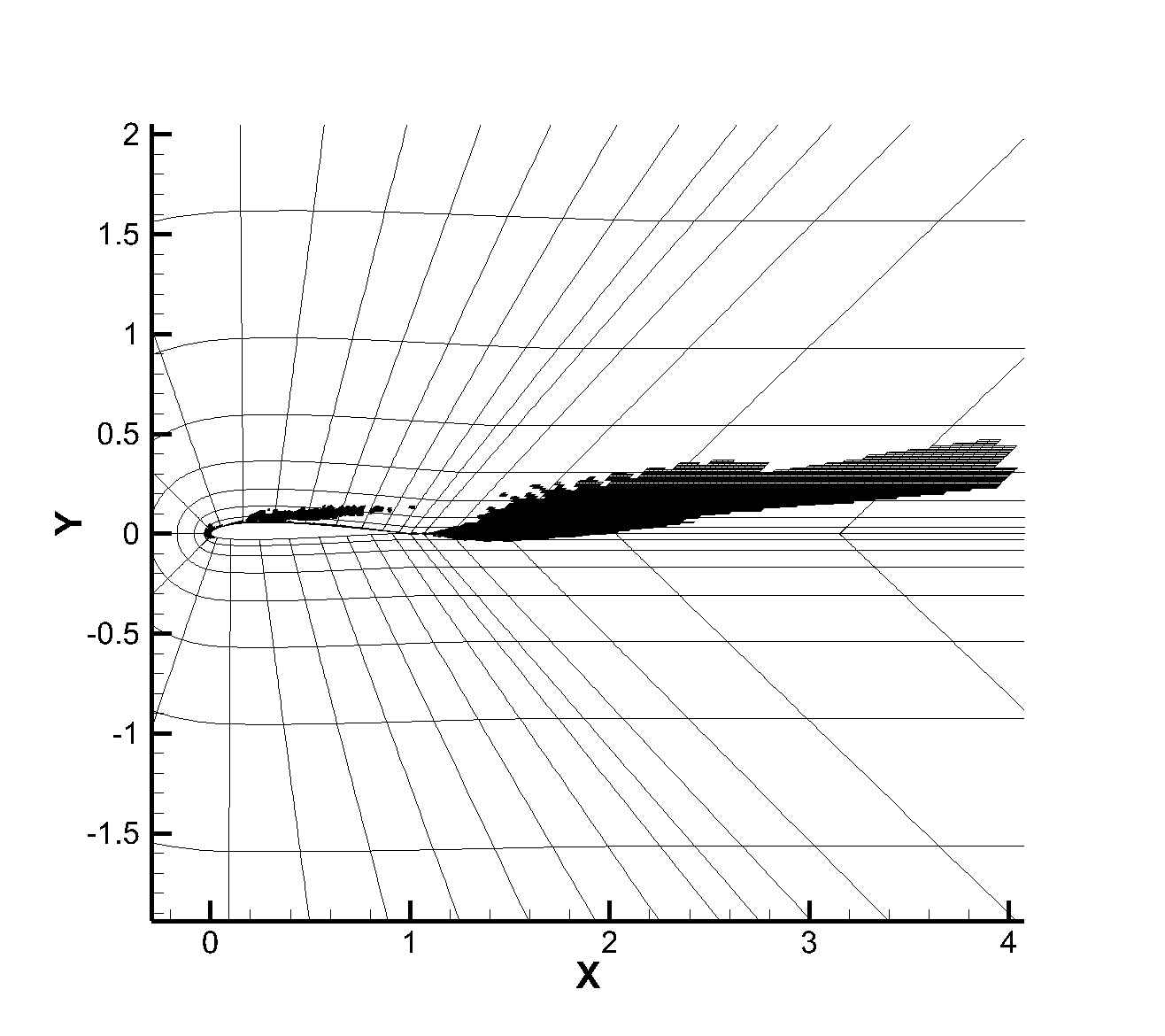}
            \caption{}
        \end{subfigure}
        \caption{Far view of adapted grids, (a): $ IQ_{k-emp}$ adapted grid; (c): $ IQ_{k-ke}$ adapted grid; (d): $ IQ_{k-tke}$ adapted grid.}\label{fig:sd7003adapted_grid_far}
    \end{figure}

\begin{figure}[hbt!]
        \centering
        \begin{subfigure}[b]{0.33\textwidth}            
            \includegraphics[clip=true, trim= 0.1cm 0.1cm 0.1cm 0.1cm, width=\linewidth]{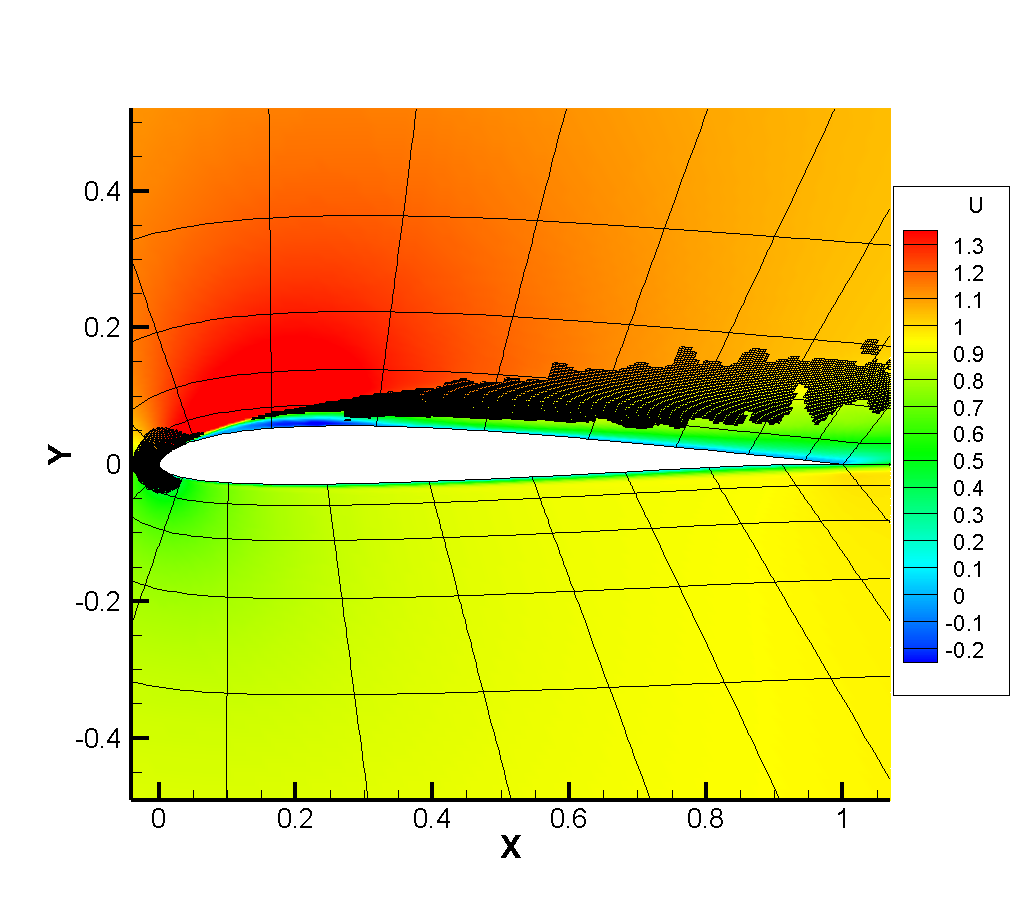}
            \caption{}
        \end{subfigure}
        \begin{subfigure}[b]{0.33\textwidth}            
            \includegraphics[clip=true, trim= 0.1cm 0.1cm 0.1cm 0.1cm, width=\linewidth]{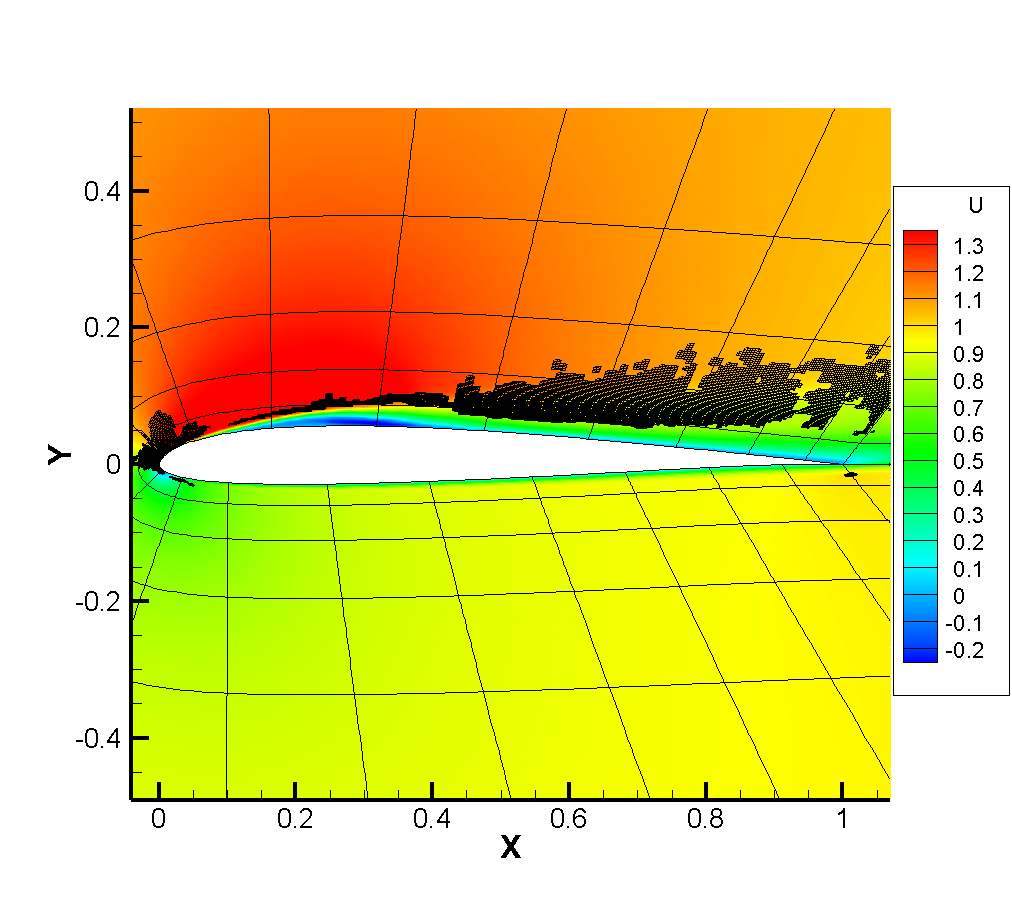}
            \caption{}
        \end{subfigure}
        \begin{subfigure}[b]{0.33\textwidth}            
            \includegraphics[clip=true, trim= 0.1cm 0.1cm 0.1cm 0.1cm, width=\linewidth]{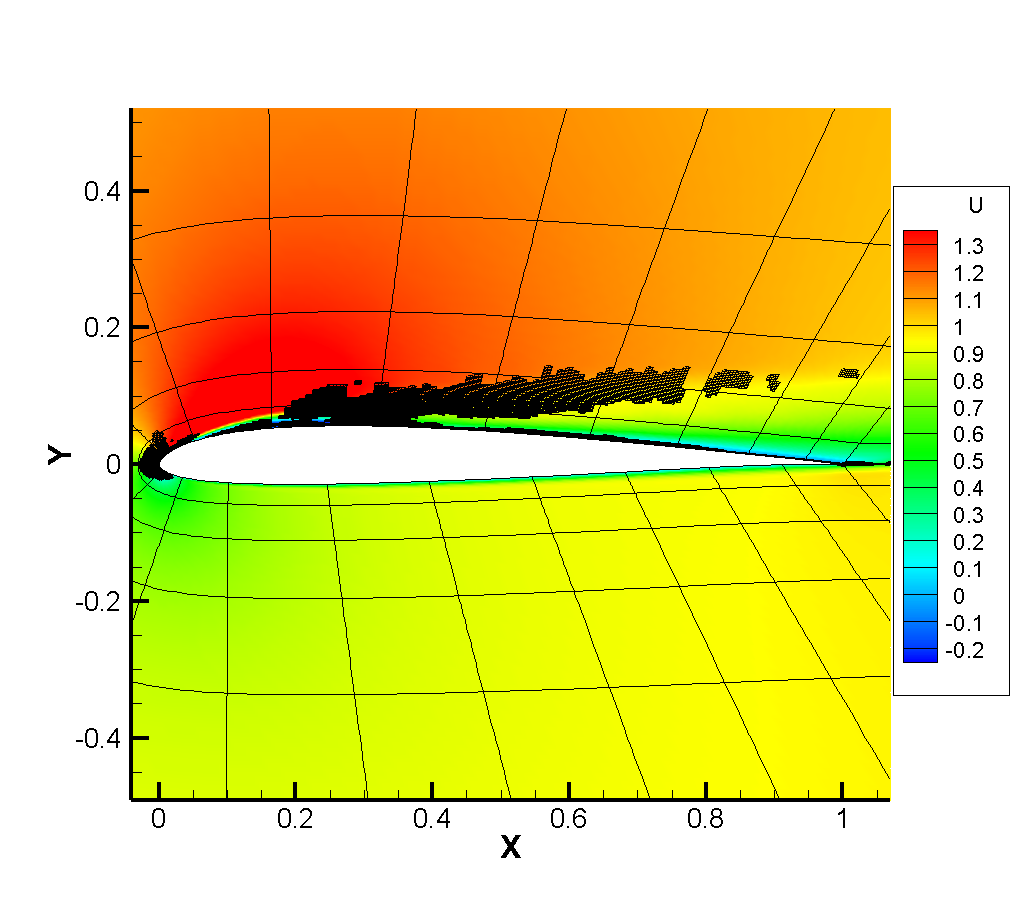}
            \caption{}
        \end{subfigure}
        \caption{Close view of adapted grids, (a): $ IQ_{k-emp}$ adapted grid; (c): $ IQ_{k-ke}$ adapted grid; (d): $ IQ_{k-tke}$ adapted grid.}\label{fig:sd7003adapted_grid}
    \end{figure}

As shown in Fig.~\ref{fig:sd7003error}, all three error estimators target the stagnation point, as well as the turbulent region in the shear layer above the separation bubble where high values of numerical TKE is estimated. $IQ_{k-ke}$ targets an even larger region above the separation bubble compared to the other error estimators; where, {\color{black} some of the cells along block interfaces are featured due to irregular transitioning between block faces. {\color{black} In the case of $IQ_{k-tke}$,} it targets the boundary layer on the upper surface due to high levels of numerical TKE as well as the entire separation bubble including the laminar separation point. It should be mentioned that the family of Index Quality error estimators is grounded on the assumption that the flow is in a fully turbulent regime~\cite{celik2009assessment}, while in the studied case, the region $X<0.2$ is laminar where both captured and modeled TKE show low values, such that the performance of the estimators in such regions primarily depends on the estimation of numerical TKE. The adapted grids are generated based on a $5\%$ refinement from the coarse grid and Fig.~\ref{fig:sd7003adapted_grid_far} gives an overview of the adapted grids based on three error estimators. As we can see the wake region is targeted by all the error estimators due to its turbulent nature. Since we are more interested in the performance of the error estimators in the near-wall region, we present the zoom-in view of the adapted grids in Fig.~\ref{fig:sd7003adapted_grid} for the three error estimators.}
    
\subsubsection{Simulation on adapted grids}
\label{sec:sd7003_sim_adapted_grids}
The simulation was carried out on the adapted grids and the results are presented in Table~\ref{tab:sd7003_summary} and Fig.~\ref{fig:sd7003bubble}. It is observed that the location of the separation and reattachment points are improved on all adapted grids, which led to earlier separation and reattachment compared to the coarse grid. Fig.~\ref{fig:sd7003bubble} shows that $IQ_{k-emp}$ and $IQ_{k-tke}$ perform better than $IQ_{k-ke}$ in terms of the prediction of the location, the length and the thickness of the separation bubble, as well as the location of the transition to turbulence, while the $IQ_{k-emp}$ adapted grid overestimates the level of TKE. The correct capture of the location and the size of the separation bubble highly depends on the correct estimation of the TKE level. {\color{black}The $C_l$ and $C_d$ values tended toward the fine grid results as expected; however, the $IQ_{k-tke}$ provided values that are in good agreement against~\cite{garmann2013comparative} instead of the fine grid. This is primarily due to a separation length that is more in line with~\cite{garmann2013comparative}.}

    \begin{figure}[hbt!]
        \begin{center}
            \includegraphics[clip=true, trim= 0.0cm 0.0cm 0.0cm 0.0cm,width=.7\columnwidth]{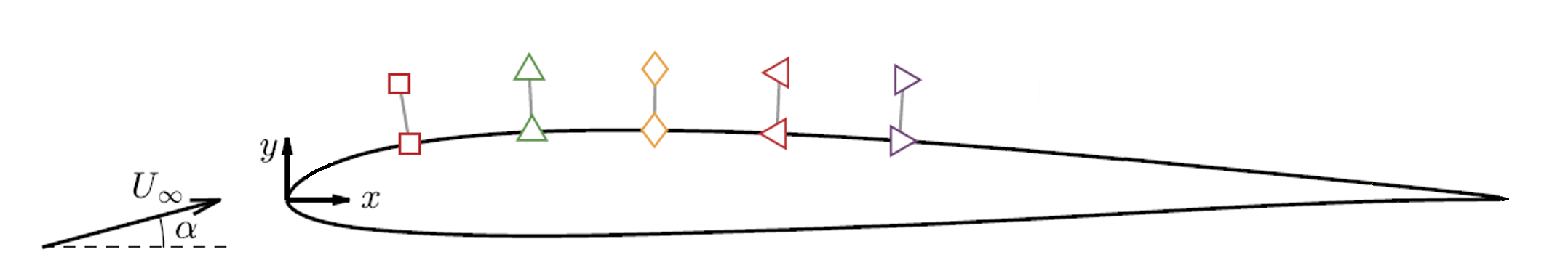}
        \end{center}
        \caption{Profile line location.}
        \label{figure:sd7003_location} 
    \end{figure}

\begin{figure}[hbt!]
        \centering
        \begin{subfigure}[b]{0.33\linewidth}            
            \includegraphics[clip=true, trim= 0.1cm 3.0cm 0.1cm 0.1cm, width=\textwidth]{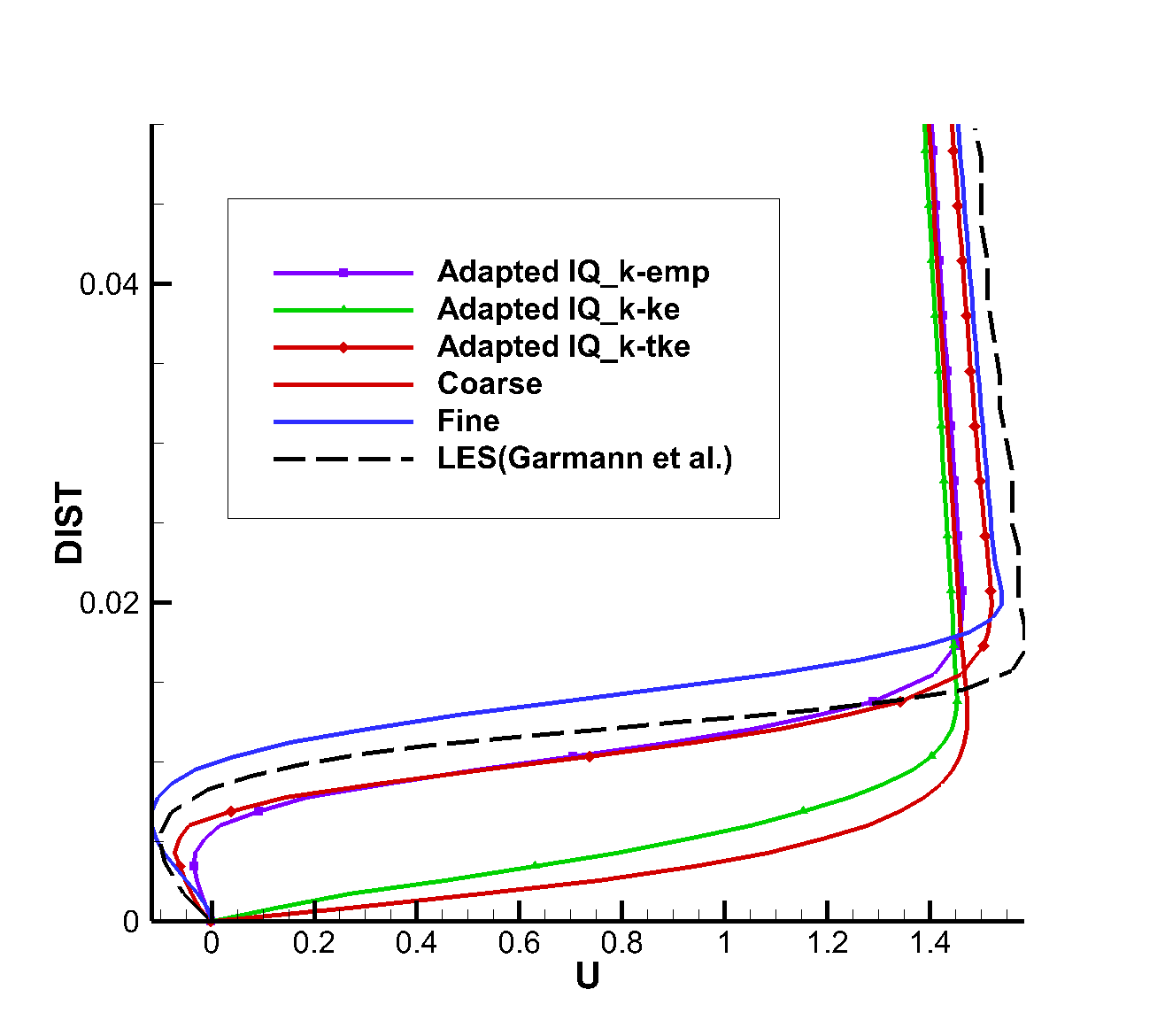}
            \label{fig:sd7003profileU01}
            \caption{}
        \end{subfigure}
        \begin{subfigure}[b]{0.33\linewidth}            
            \includegraphics[clip=true, trim= 0.1cm 3.0cm 0.1cm 0.1cm, width=\textwidth]{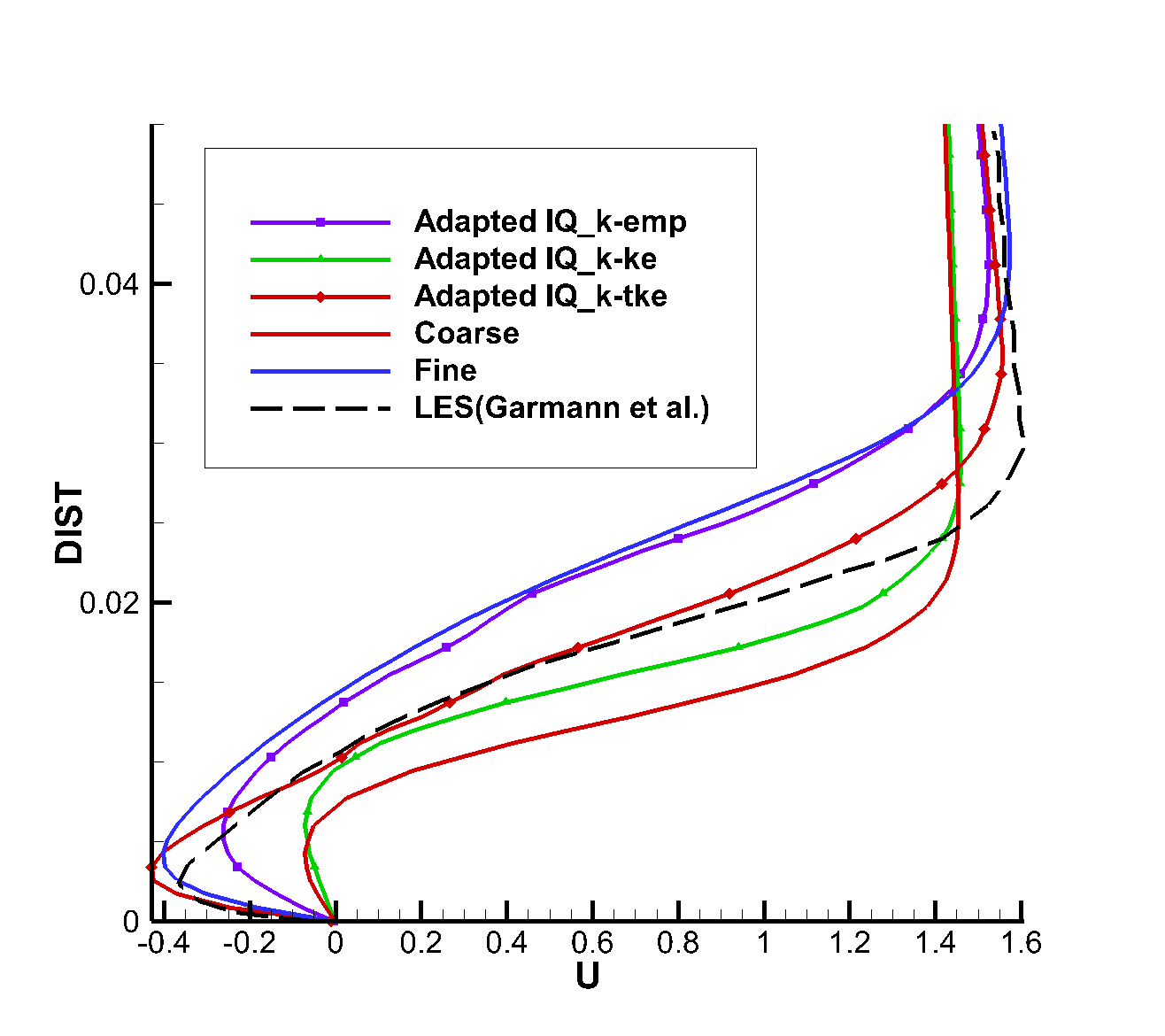}
            \label{fig:sd7003profileU02}
            \caption{}
        \end{subfigure}
        \begin{subfigure}[b]{0.33\linewidth}            
            \includegraphics[clip=true, trim= 0.1cm 3.0cm 0.1cm 0.1cm, width=\textwidth]{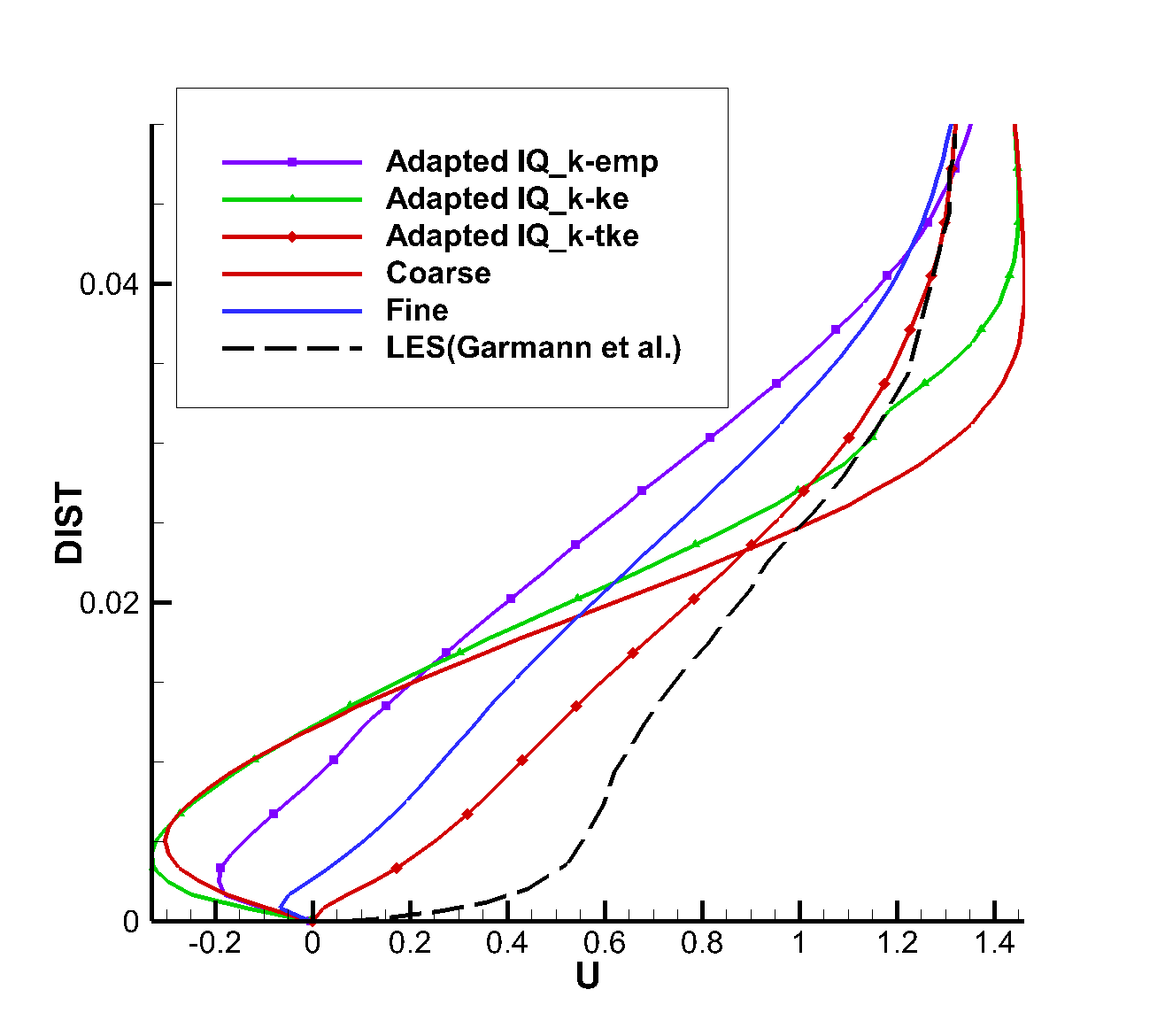}
            \label{fig:sd7003profileU03}
            \caption{}
        \end{subfigure}
        \begin{subfigure}[b]{0.33\linewidth}
            \includegraphics[clip=true, trim= 0.1cm 3.0cm 0.1cm 0.1cm, width=\textwidth]{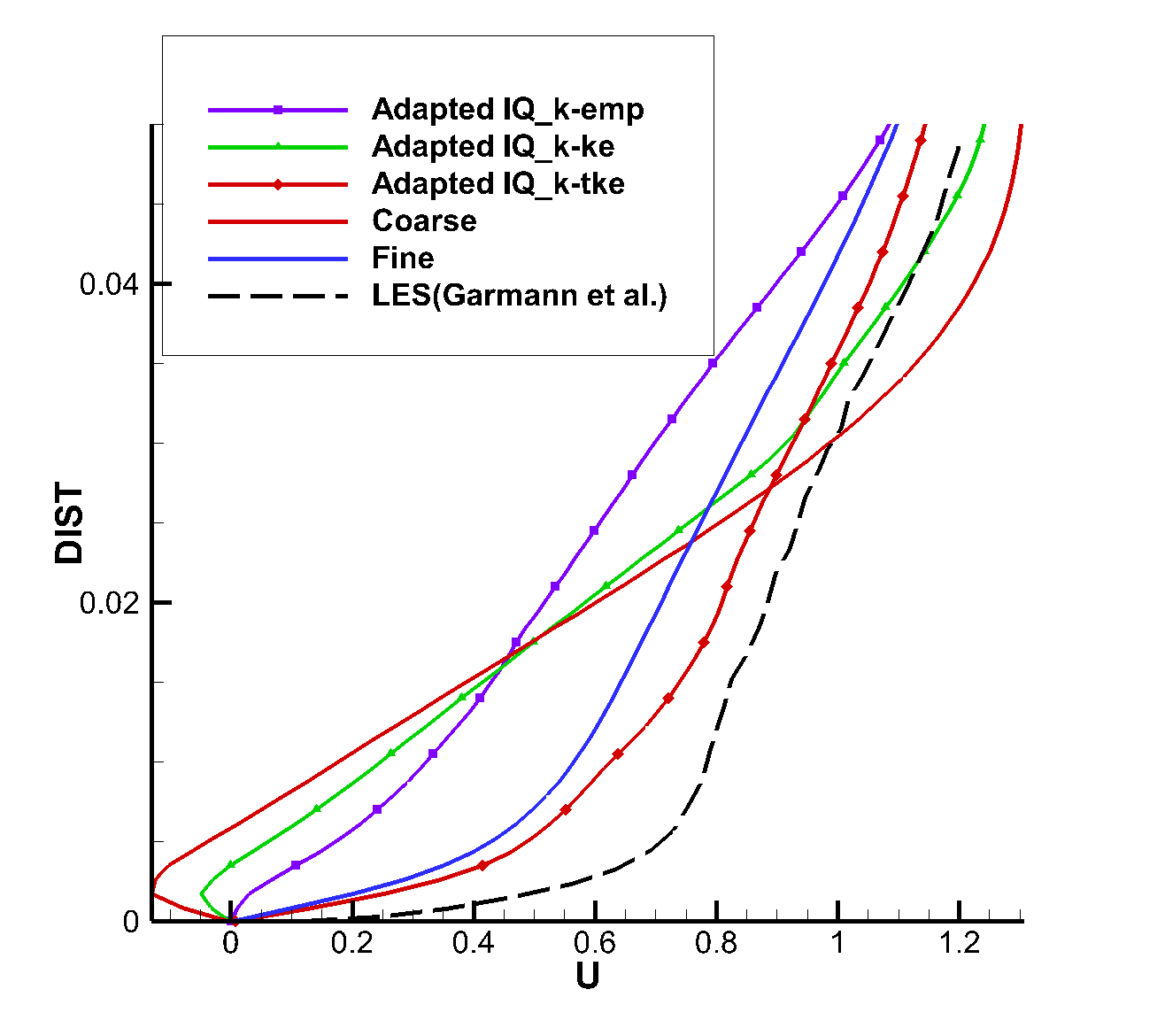}
            \label{fig:sd7003profileU04}
            \caption{}
        \end{subfigure}
        \begin{subfigure}[b]{0.33\linewidth}
            \includegraphics[clip=true, trim= 0.1cm 3.0cm 0.1cm 0.1cm, width=\textwidth]{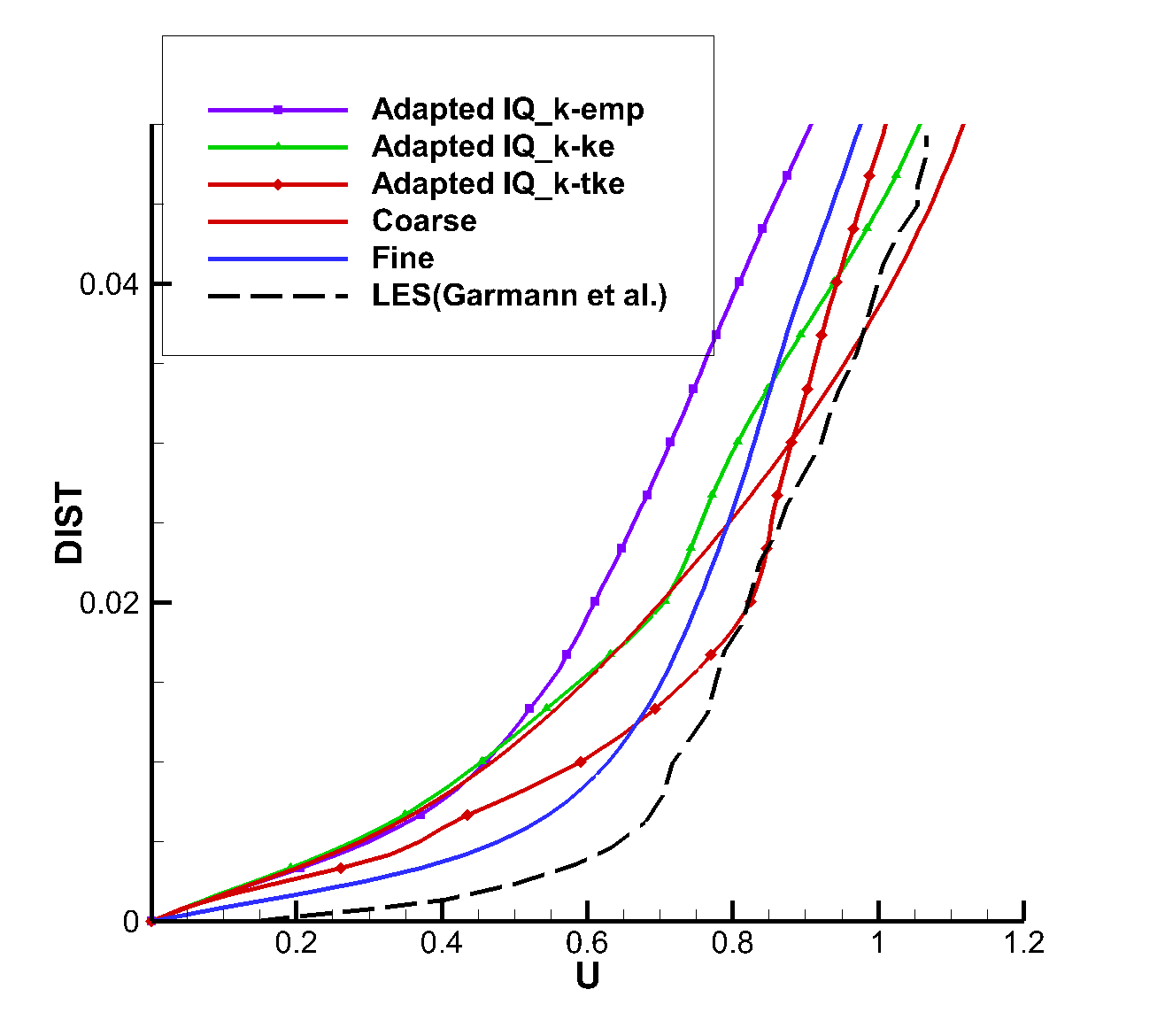}
            \label{fig:sd7003profileU05}
            \caption{}
        \end{subfigure}
        \caption{Streamwise mean velocity $U$ profiles at (a): $ X=0.1 $; (b): $ X=0.2 $; (c): $ X=0.3 $; (d): $ X=0.4 $; (e): $ X=0.5 $. }\label{fig:sd7003Uprofile}
\end{figure}

\begin{figure}[hbt!]
        \centering
        \begin{subfigure}[b]{0.33\linewidth}            
            \includegraphics[clip=true, trim= 0.1cm 3.0cm 0.1cm 2.1cm, width=\textwidth]{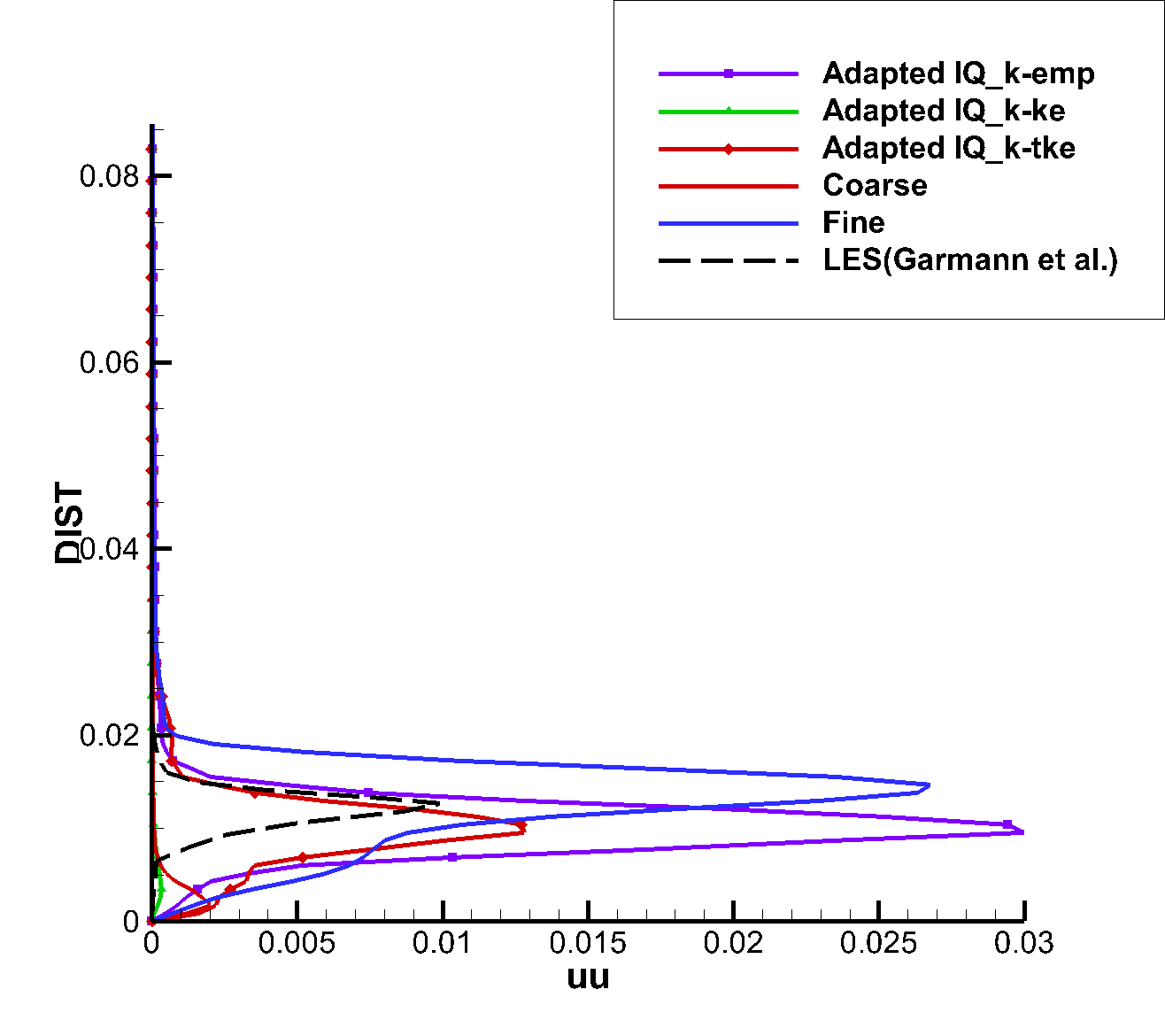}
            \label{fig:sd7003profileuu01}
            \caption{}
        \end{subfigure}
        \begin{subfigure}[b]{0.33\linewidth}            
            \includegraphics[clip=true, trim= 0.1cm 3.0cm 0.1cm 2.1cm, width=\textwidth]{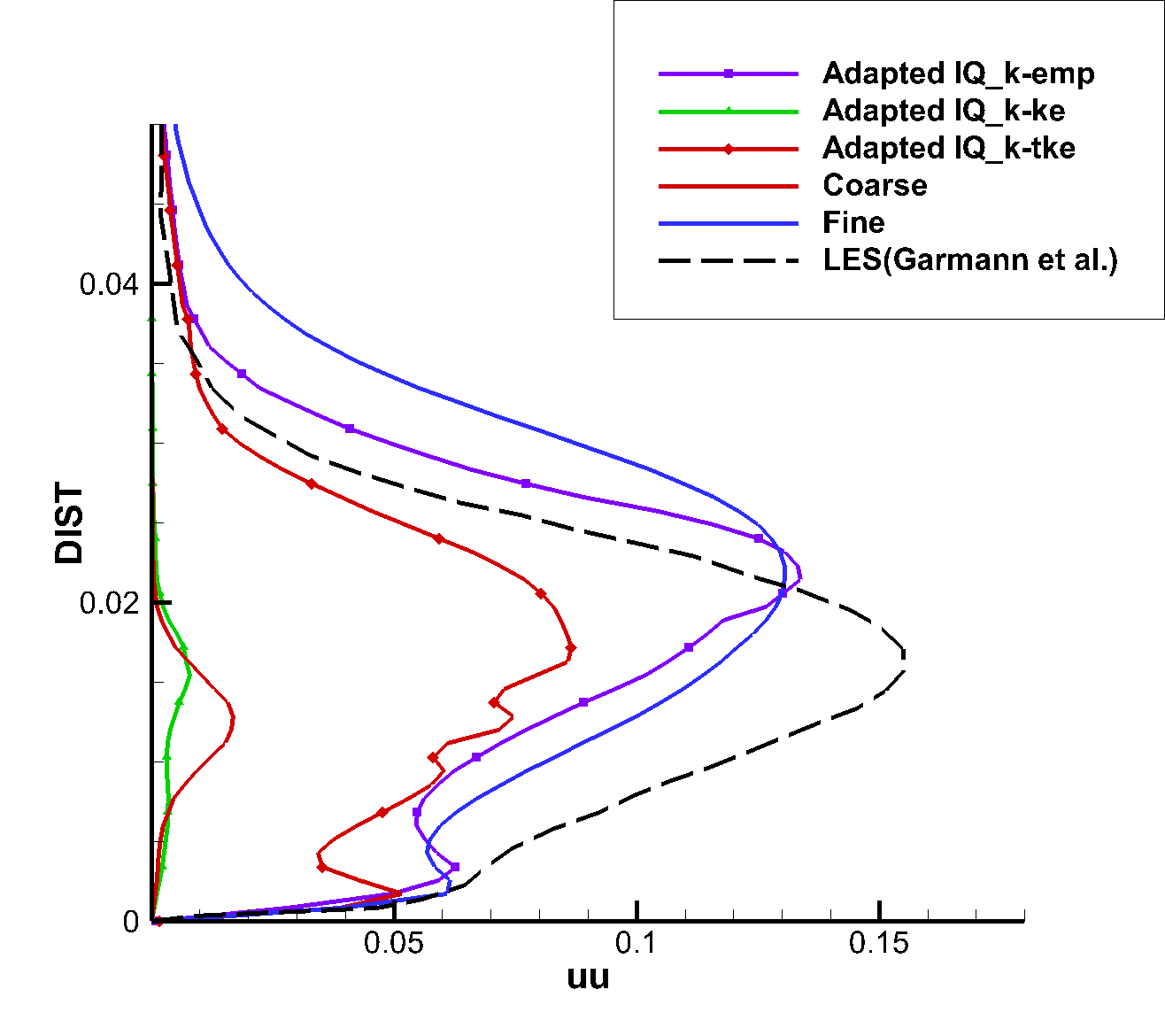}
            \label{fig:sd7003profileuu02}
            \caption{}
        \end{subfigure}
        \begin{subfigure}[b]{0.33\linewidth}            
            \includegraphics[clip=true, trim= 0.1cm 3.0cm 0.1cm 2.1cm, width=\textwidth]{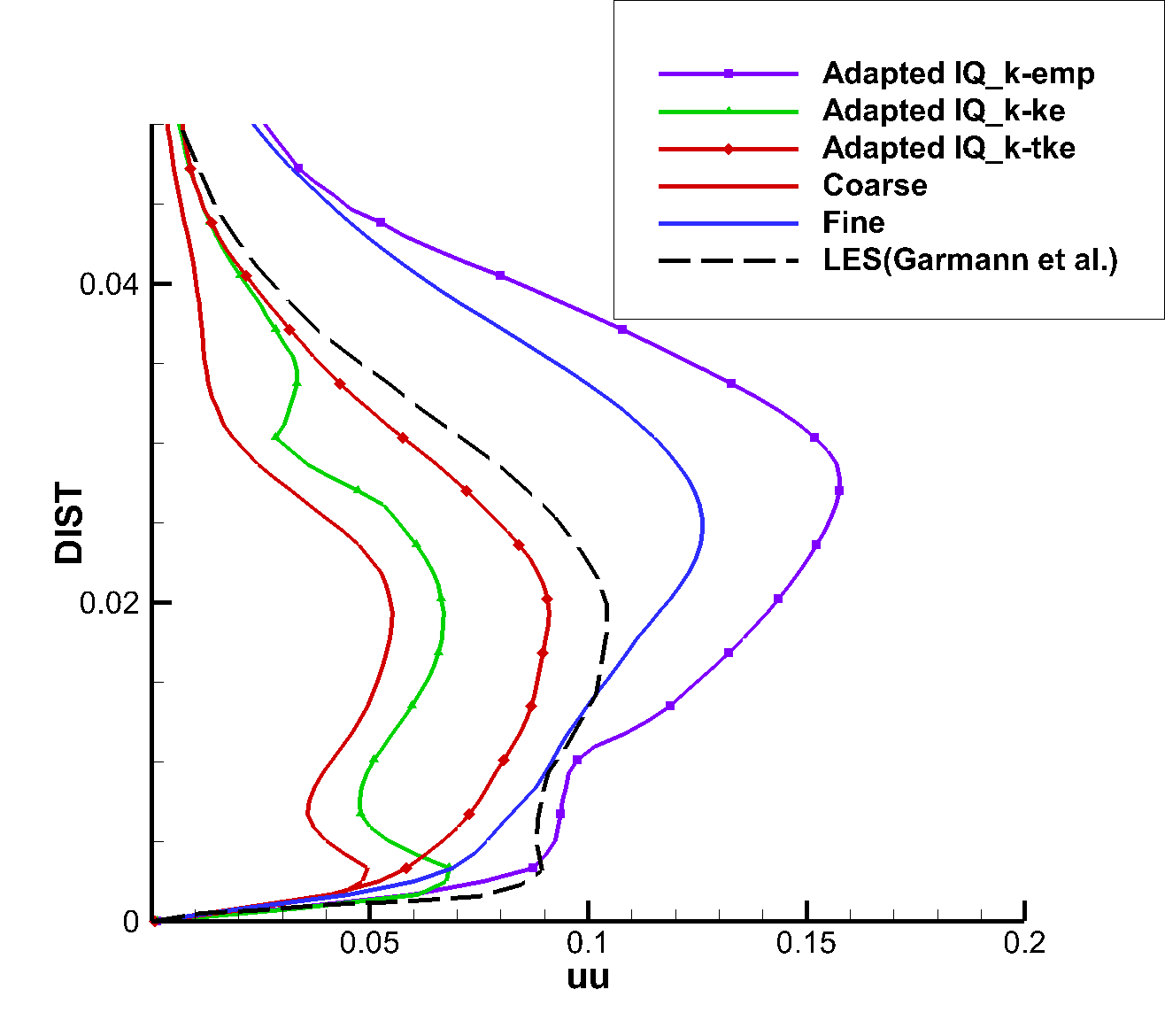}
            \label{fig:sd7003profileuu03}
            \caption{}
        \end{subfigure}
        \begin{subfigure}[b]{0.33\linewidth}
            \includegraphics[clip=true, trim= 0.1cm 3.0cm 0.1cm 2.1cm, width=\textwidth]{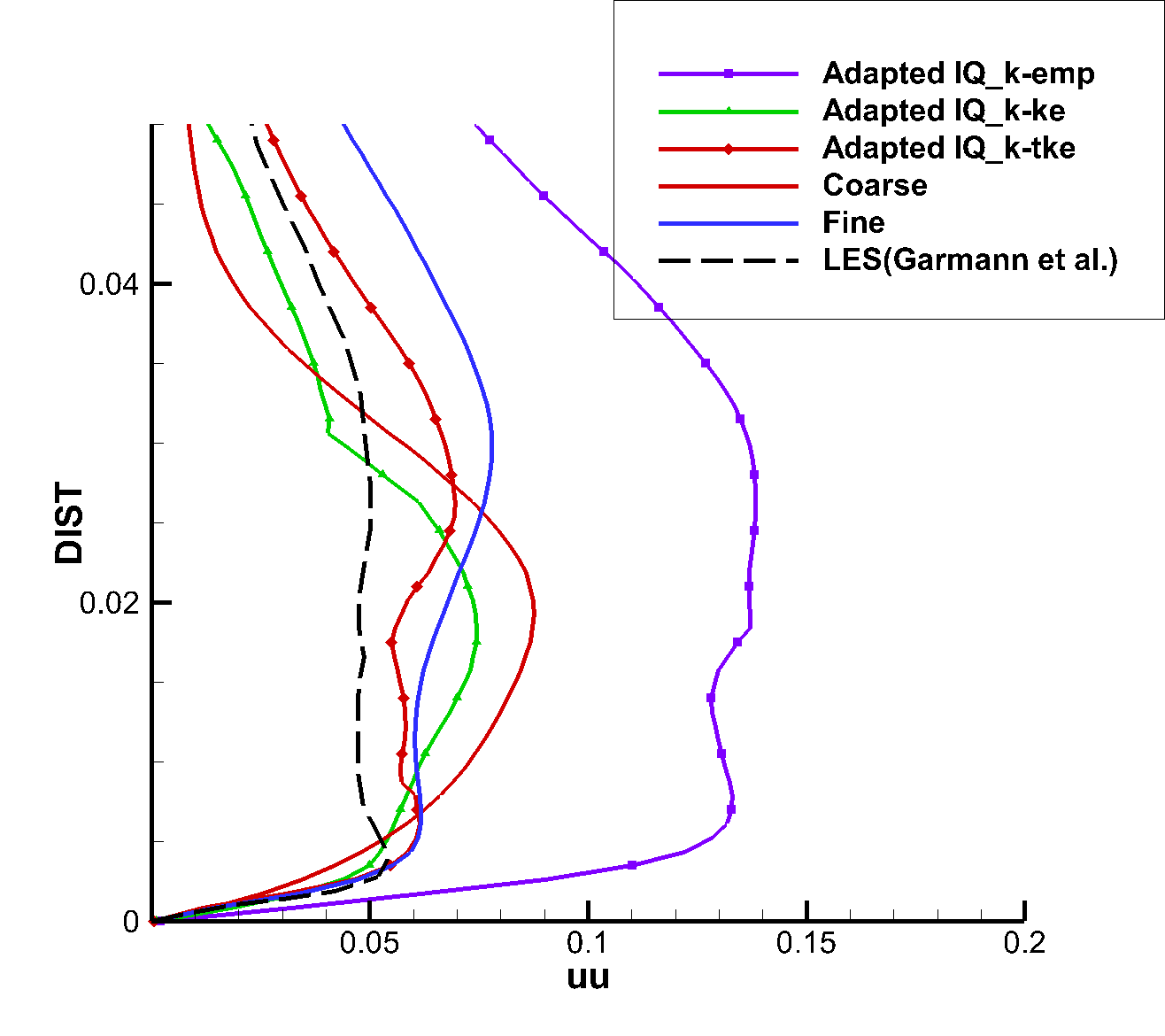}
            \label{fig:sd7003profileuu04}
            \caption{}
        \end{subfigure}
        \begin{subfigure}[b]{0.33\linewidth}
            \includegraphics[clip=true, trim= 0.1cm 3.0cm 0.1cm 2.1cm, width=\textwidth]{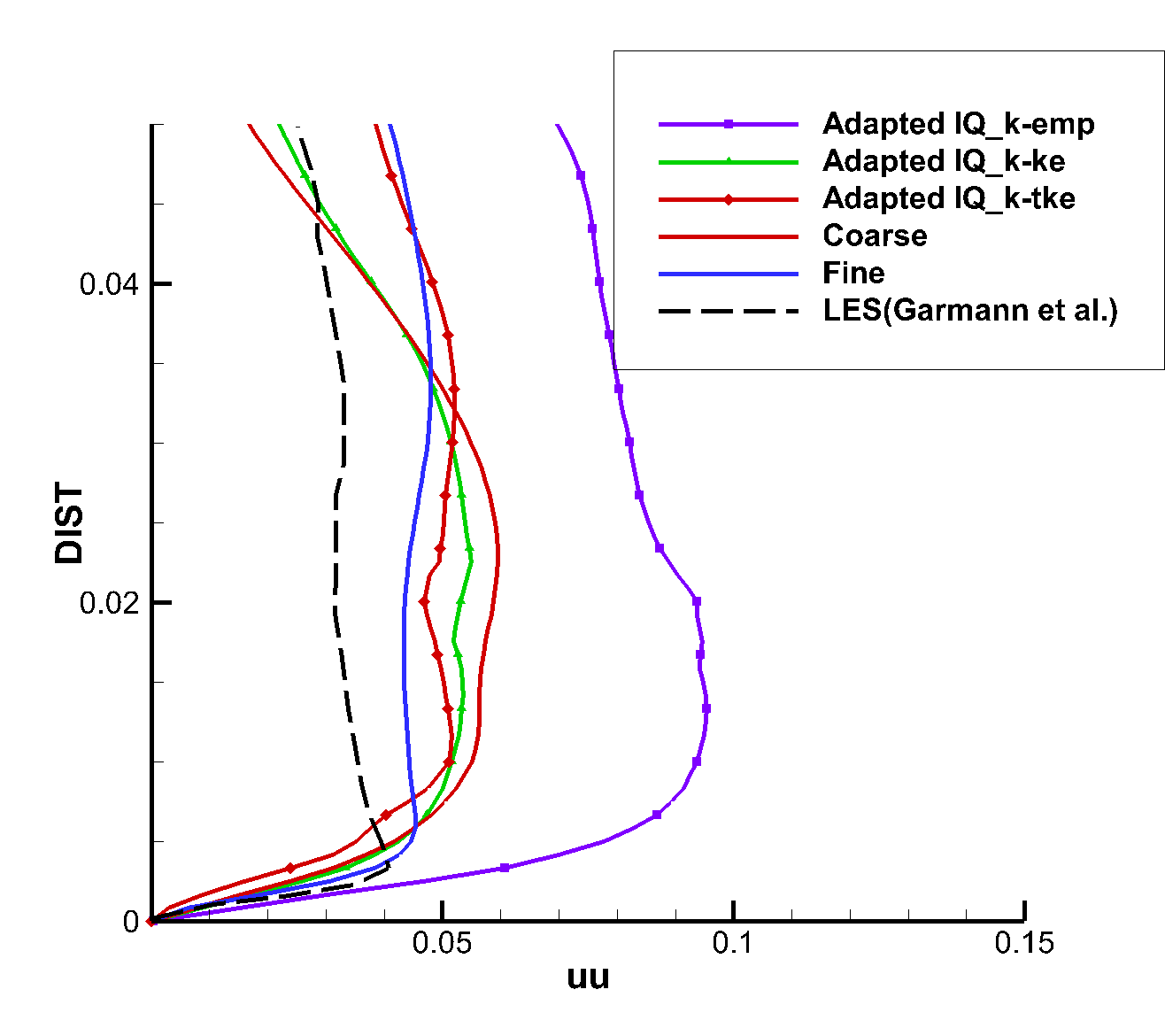}
            \label{fig:sd7003profileuu05}
            \caption{}
        \end{subfigure}
        \caption{$u'u'$ profiles at (a): $ X=0.1 $; (b): $ X=0.2 $; (c): $ X=0.3 $; (d): $ X=0.4 $; (e): $ X=0.5 $. }\label{fig:sd7003uuprofile}
\end{figure}

Profiles normal to the wall are studied in five different locations as shown in Fig.~\ref{figure:sd7003_location}. The profiles are taken every $0.1$c at $X= 0.1$ to $0.5$ and provide an insight into the size of the separated region. The streamwise mean velocity $U$ and the squared fluctuations $u'u'$ are studied at the stated locations and compared against the reference values. {\color{black}Fig.~\ref{fig:sd7003Uprofile} shows that the $IQ_{k-tke}$ adapted grid leads to the best prediction of the mean velocity at all five locations by providing results close to that on the fine grid and reference LES. Fig.~\ref{fig:sd7003Uprofile}~(c) shows that the flow is reattached at $X=0.3$ for the $IQ_{k-tke}$ adapted grid similar to the reference LES, while other grids still show backflow in the boundary layer. The octree data structure in the code results in the refinement of the targeted cells by a factor of 8, such that the adapted grids are in fact finer than the fine grid in the targeted regions. The $u'u'$ profiles are shown in Fig.~\ref{fig:sd7003uuprofile}. We observe that at $X=0.1$, only the $IQ_{k-tke}$ adapted grid is able to estimate the correct level of TKE. At a downstream location of $X=0.2$, the coarse grid and the $IQ_{k-ke}$ adapted grid show very low level of captured $u'u'$ and fail to predict the start of transition. The $IQ_{k-tke}$ adapted grid also underestimates the level of TKE at this location. Further downstream at $X=0.3$, which is near the flow reattachement, the $IQ_{k-tke}$ adapted grid captures $u'u'$ at the same level as reference LES~\cite{garmann2013comparative} while other grids show either an overestimation ($IQ_{k-ke}$ adapted grid and fine grid) or an underestimation ($IQ_{k-emp}$ adapted grid and coarse grid). At the downstream locations $X=0.4$ and $X=0.5$ the $IQ_{k-tke}$ adapted grid compares well against the reference LES and fine grid values, and outperforms $IQ_{k-emp}$ and $IQ_{k-ke}$ adapted grids.}
    

\section{Conclusion and Future Work}
    
In this work, two new approaches are proposed to evaluate the numerical TKE and estimate the Index Quality for LES grid adaptation. The first approach relies on the numerical dissipation of the kinetic energy equation and demonstrated promising outcomes but its propensity to result in negative values in laminar regions reduced its effectiveness. The second approach relies on the numerical dissipation of the turbulent kinetic energy equation for LES. The latter approach overcame the short-comings of the KE-based technique by stressing high numerical dissipation within boundary layers and turbulent regions. The proposed methods were compared against a classical empirical method for evaluating the numerical TKE. All three approaches were employed to develop a family of Index Quality error estimators: $IQ_{k-emp}$, $IQ_{k-ke}$, and $IQ_{k-tke}$ for LES. Its capabilities were demonstrated for the periodic hill and the SD7003 airfoil test cases, where we contrasted their differences by first studying the distribution of the numerical TKE then contours of Index Quality. All three IQ error estimators target the shear layer and turbulent mixing layer on the top of the separation bubble for the periodic hill and the separated region for the SD 7003 airfoil. However, the developed $IQ_{k-tke}$ added the boundary layer and reduced if not eliminated the tendency of the KE-based approach from producing negative numerical dissipation. 

One complete grid adaptation cycle is performed for both test cases. Generally, the adapted grids based on all error estimators are able to improve the quality of the LES results. The $IQ_{k-emp}$ estimator was able to target the turbulent mixing layer and the adapted grid improved the size of the separation bubble; however, the estimator failed to target the boundary layer. Unlike, the empirical approach, the $IQ_{k-ke}$ estimator was able to target the near-wall region however previously confirmed reports~\cite{castiglioni2015numerical} of negative values were observed. As such the estimator failed to target the laminar separation bubble in the SD7003 case. We found that the $IQ_{k-tke}$ estimator was able to balance the refinement between the mixing layer above the bubble and the boundary layer for both cases, which led to the best performance in terms of capturing correctly the Reynolds stress tensor component values and the length of the separation bubble on the corresponding adapted grid. Future work will see further validation of the technique for high Reynolds number cases with the presence of shocks and three-dimensional flows with greater inhomogeneous spanwise variations.

\bibliography{main}

\end{document}